\documentstyle[amscd,verbatim,leqno,12pt]{amsart}

\numberwithin{equation}{section}

\theoremstyle{plain}
\newtheorem{thm}{Theorem}[section]
\newtheorem{defn}[thm]{Definition}
\newtheorem{prop}[thm]{Proposition}
\newtheorem{lem}[thm]{Lemma}
\newtheorem{cor}[thm]{Corollary}
\newtheorem{conj}[thm]{Conjecture}
\newtheorem{propdefn}[thm]{Proposition-Definition}

\theoremstyle{definition}
\newtheorem{rem}[thm]{Remark}

\renewcommand{\b}{\bullet}
\newcommand{\beast}{\begin{eqnarray*}}
\newcommand{\east}{\end{eqnarray*}}

\newcommand{\N}{{\Bbb N}}
\newcommand{\Z}{{\Bbb Z}}
\newcommand{\Q}{{\Bbb Q}}

\newcommand{\F}{{\Bbb F}}

\newcommand{\M}{{\mathrm{M}}}
\newcommand{\MC}{{\mathrm{MC}}}

\newcommand{\Cone}{{\mathrm{Cone}}}
\newcommand{\Spec}{{\mathrm{Spec}}\,}

\newcommand{\Spf}{{\mathrm{Spf}}\,}
\newcommand{\Spm}{{\mathrm{Spm}}\,}
\newcommand{\Zar}{{\mathrm{Zar}}}
\newcommand{\lra}{\longrightarrow}
\newcommand{\ra}{\rightarrow}
\newcommand{\hra}{\hookrightarrow}

\newcommand{\lla}{\longleftarrow}

\newcommand{\isom}{\overset{\sim}{=}}
\newcommand{\ti}[1]{\widetilde{#1}}
\newcommand{\wt}[1]{\widetilde{#1}}

\newcommand{\ul}[1]{\underline{#1}}
\newcommand{\ol}[1]{\overline{#1}}
\newcommand{\os}{\overset}

\newcommand{\DR}{{\mathrm{DR}}}
\newcommand{\et}{{\mathrm{et}}}

\newcommand{\crys}{{\mathrm{crys}}}

\newcommand{\conv}{{\mathrm{conv}}}

\newcommand{\Hom}{{\mathrm{Hom}}}

\newcommand{\Ker}{{\mathrm{Ker}}}

\newcommand{\Coker}{{\mathrm{Coker}}}

\newcommand{\triv}{{\mathrm{triv}}}

\newcommand{\Mod}{{\mathrm{Mod}}}

\newcommand{\id}{{\mathrm{id}}}

\newcommand{\gp}{{\mathrm{gp}}}

\newcommand{\Coh}{{\mathrm{Coh}}}
\newcommand{\Strat}{{\mathrm{Str}}}

\newcommand{\red}{{\mathrm{red}}}

\newcommand{\HPD}{{\mathrm{HPD}}}
\newcommand{\HPDI}{{\mathrm{HPDI}}}
\newcommand{\pr}{{\mathrm{pr}}}

\newcommand{\an}{{\mathrm{an}}}

\newcommand{\D}{D}

\newcommand{\cA}{{\cal A}}
\newcommand{\cB}{{\cal B}}
\newcommand{\cC}{{\cal C}}
\newcommand{\cD}{{\cal D}}
\newcommand{\cE}{{\cal E}}
\newcommand{\cF}{{\cal F}}
\newcommand{\cG}{{\cal G}}
\newcommand{\cH}{{\cal H}}
\newcommand{\cI}{{\cal I}}
\newcommand{\cM}{{\cal M}}
\newcommand{\cN}{{\cal N}}
\newcommand{\cO}{{\cal O}}
\newcommand{\cP}{{\cal P}}
\newcommand{\cQ}{{\cal Q}}
\newcommand{\cK}{{\cal K}}
\newcommand{\cS}{{\cal S}}

\newcommand{\cU}{{\cal U}}

\newcommand{\cW}{{\cal W}}
\newcommand{\cX}{{\cal X}}
\newcommand{\cY}{{\cal Y}}
\newcommand{\cZ}{{\cal Z}}

\newcommand{\fP}{{\frak P}}

\renewcommand{\sp}{{\mathrm{sp}}}
\newcommand{\rig}{{\mathrm{rig}}}
\newcommand{\ex}{{\mathrm{ex}}}

\newcommand{\vpl}{\varprojlim}

\newcommand{\wh}{\widehat}
\renewcommand{\wt}{\widetilde}

\newcommand{\LF}{{({\mathrm{LFS}}/\cB)}}
\newcommand{\pLF}{{(p{\mathrm{LFS}}/\cB)}}
\renewcommand{\L}{{({\mathrm{LS}}/\cB)}}
\newcommand{\LB}{{({\mathrm{LS}}/B)}}

\renewcommand{\d}{\dagger}
\renewcommand{\dag}{\dagger}
\newcommand{\dd}{\ddagger}
\newcommand{\lam}{\lambda}
\newcommand{\Lam}{\Lambda}
\newcommand{\lr}{\langle \b \rangle}
\newcommand{\lmr}{\langle m \rangle}

\newcommand{\tls}{{\mathrm{triv. \allowbreak \, log 
\allowbreak \, str.}}}

\begin{document}

\title[Relative Log Convergent Cohomology I]
{Relative Log Convergent Cohomology and \\ Relative Rigid Cohomology I}

\author{Atsushi Shiho}
\address{Graduate School of Mathematical Sciences, 
University of Tokyo, 3-8-1 Komaba, Meguro-ku, Tokyo 153-8914, JAPAN.}
                          

\begin{abstract}
In this paper, we develop the theory of relative log convergent 
cohomology. 
We prove the coherence of relative log convergent 
cohomology in certain case 
by using the comparison theorem between 
relative log convergent cohomlogy and relative log crystalline 
cohomology, and we relates relative log convergent 
cohomology to relative rigid cohomology to 
show 
the validity of Berthelot's conjecture 
on the coherence and the 
overconvergence of relative rigid cohomology for proper smooth families 
when they admit nice proper log smooth compactification 
to which the coefficient extends logarithmically. 
\end{abstract}

\maketitle

\tableofcontents


\section*{Introduction}
Let $k$ be a perfect field of characteristic $p>0$ and let 
$V$ be a complete discrete valuation ring of mixed characteristic 
with residue field $k$. For a scheme $X$ of finite type over $k$, 
theory of convergent cohomology of $X$ over $V$ is developped in 
\cite{ogus1} and \cite{ogus2}. In \cite{ogus2}, Ogus proved 
two comparison theorems concerning convergent cohomology: 
One is the comparison theorem between convergent cohomology 
and crystalline cohomology for smooth schemes over $k$ 
and the other one is the comparison theorem 
between convergent cohomology and 
rigid cohomology for proper smooth schemes over $k$. As a consequence, he 
deduced the finiteness of rigid cohomology for proper smooth schemes 
over $k$ with coefficient (which is an overconvergent isocrystal). \par 
The theory of convergent cohomology is extended to the case of 
log schemes in \cite{ogus3}, \cite{shiho1}, \cite{shiho2}. 
In the paper \cite{shiho2}, we proved 
the comparison between log convergent 
cohomology and log crystalline cohomology for log smooth log schemes 
over $k$ and 
the comparison between log convergent cohomology and rigid cohomology 
for proper log smooth log schemes over $k$ 
under certain condition on coefficients. 
As a consequence, we proved the finiteness of 
rigid cohomology of a smooth $k$-scheme 
$X$ with coefficient $\cE$ (which is an overconvergent 
$F$-isocrystal) when $X$ admits a 
proper log smooth compactification to 
which $\cE$ extends logarithmically. 
(Note that, soon after the paper \cite{shiho2} appeared, 
the finiteness of rigid 
cohomology with coefficient is proved without any assumption by 
Kedlaya \cite{kedlaya1}. However, we think that 
the approach in \cite{shiho2} is 
still interesting. Recently Kedlaya has given a proof of 
a conjecture proposed in \cite{shiho2}, which also implies the finiteness of 
rigid cohomology in general case. See \cite{kedlaya2} and his 
subsequent papers for detail.) \par 
The purpose of this paper is to develop the theory 
of log convergent cohomology in relative situation. 
(We also remark here that the field $k$ in the main body of this paper 
is not necessarily perfect before the Frobenius structure comes in.) 
We prove the comparison theorem between relative log convergent cohomology and 
relative log crystalline cohomlogy and then we relates relative log convergent 
cohomology to the relative rigid cohomology (defined in \cite{chts}, 
see also \cite{berthelot1}). 
As a consequence of the former result, 
we deduce the coherence and the log convergence 
of the relative log convergent cohomology in certain case and 
by combining it with the latter result, we deduce the coherence and 
the overconvergence of relative rigid cohomology 
for proper smooth families when they admit nice proper log smooth 
compactification to which the coefficient (which is an overconvergent 
isocrystal) extends logarithmically. Note that, in relative situation, 
we need to develop the theory of `log' convergent cohomology 
even to treat only proper smooth families 
because we will need `log structure along boundary' of the given 
family to prove the coherence and the overconvergence. \par 
The coherence and the overconvergence of relative rigid cohomology 
for proper smooth families is conjectured by Berthelot (\cite{berthelot1}, 
see also \cite{tsuzuki3}). 
(However, we would like to note that 
we can formulate the conjecture in several ways. For the discussion on 
the statement of the conjecture, see Section 5.) 
Our result gives an affirmative answer to a version of this conjecture under 
certain assumption. This conjecture of Berthelot has been known 
in the case when the given family lifts to proper smooth families of 
formal schemes (\cite{berthelot1}, 
\cite{tsuzuki3}) and in the case where a given family is a family of 
abelian varieties and the coefficient is trivial (\cite{etesse}).
There is also an result (\cite{matsudatrihan}) treating the case 
when the base scheme of the family is a smooth curve and 
the coefficient is trivial, 
whose proof depends on an unpublished result of 
Kedlaya (\cite[6.1]{kedlaya0}). 
A weaker conjecture (generic coherence and 
generic overconvergence) is known in the case of 
proper smooth families of curves (\cite{tsuzuki3}). Our result seems 
to be applicable to many cases 
where the given family is not necessarily liftable 
(after any etale base change) 
to a family of formal schemes. \par 
Let us explain the content of this paper in more detail. 
In Section 1, we give a review of the relative log crystalline 
cohomology (defined in \cite{kkato}) and 
prove the finiteness and the base change property of it 
for proper log smooth integral morphisms which are slightly 
stronger than the known results in the literature. 
In Section 2, 
we give 
the definition and basic properties of relative log convergent 
cohomology 
and we prove the relative version of log convergent Poincar\'e lemma and 
the comparison theorem with relative log crystalline cohomology for 
log smooth morphisms. 
Some of the results in this section 
are also considered in \cite{nakkshiho}. Although the 
proofs are basically the same as the absolute case 
(given in \cite[Chapter 2]{shiho2}), 
we give (a sketch of) proofs by two reasons: One is for the completeness 
of this paper and the other is to point out that some parts of the proofs 
in \cite{shiho1}, \cite{shiho2} can be generalized and simplified. 
In particular, in this paper, there are no assumptions in the statement 
of theorems which require the log schemes to be `of Zariski type' 
(for definition, see \cite[1.1.1]{shiho2}). 
We also give correction to some errors in \cite{shiho2}. 
In Section 3, 
we prove the finiteness and the base change property of 
relative log convergent cohomology for 
proper log smooth integral morphisms `having log smooth parameter' (for 
definition, see Definition \ref{parameter}). We prove it by using 
the comparison theorem with relative log crystalline cohomology proved 
in the previous section, but we need a little argument here, because 
relative log crystalline cohomology is defined only 
in the case where the base log scheme admits a PD-closed immersion into 
a log (formal) scheme, while we do not need PD-structure in the case of 
relative log convergent cohomology. 
In Section 4, 
we introduce an analytic variant of 
relative log convergent cohomology (which we call relative log 
analytic cohomology)
and compare it with relative log convergent 
cohomology for proper log smooth integral morphisms having log smooth 
paramter. 
In relative case, relative log analytic cohomology is a 
sheaf of certain rigid analytic space, while relative log convergent 
cohomology is a sheaf on certain formal scheme. We prove that the latter 
is nothing but the specialization of the former in appropriate situation. 
From this result, we 
deduce the coherence and the log convergence 
of relative log analytic cohomology. 
In Section 5, after giving a review of relative rigid cohomology 
(defined in \cite{chts}), 
we prove the comparison 
theorem between relative log analytic cohomology and relative rigid 
cohomology for proper log smooth integral morphisms having log smooth 
parameter. This implies 
the coherence and the overconvergence of 
relative rigid cohomology for proper smooth families when they 
admit a proper log smooth integral 
compactification having log smooth parameter 
to which the given coefficient (which is an overconvergent isocrystal) 
extends logarithmically. \par 
We plan to write the following topics in forthcoming papers: 
First, we plan to establish variants of relative log convergent 
cohomology and relative rigid cohomology 
to prove (a version of) Berthelot's conjecture in general case. 
Second, we plan to discuss on the generic coherence and the 
generic overconvergence of relative rigid cohomology for 
non-smooth, non-proper families. \par 
The author would like to thank Pierre Berthelot 
for pointing out first that the construction of log tubular neighborhood 
in \cite{shiho2} can be done without `of Zariski type' hypothesis. 
He would like to thank Bruno Chiarellotto and Marianna 
Fornasiero for pointing out the same thing: Discussion with them on 
their paper \cite{chfo} was helpful also to him. He would like to 
thank to Yukiyoshi Nakkajima for allowing him to include 
some results which have overlaps with the joint work 
\cite{nakkshiho} in this paper. 
He would like to thank to Kazuya Kato for inviting him to give a talk 
at a conference held at Kyoto University, and to Nobuo Tsuzuki and 
Makoto Matsumoto for inviting him 
to give a talk at a conference held at Hiroshima University. 
The author is partly supported by Grant-in-Aid for Young Scientists (B), 
the Ministry of Education, Culture, Sports, Science and Technology of 
Japan, and JSPS Core-to-Core program 18005 whose representative is 
Makoto Matsumoto. 

\section*{Convention}
\hspace{-16pt}
(1) \,\, Throughout this paper, $k$ is a field of characteristic 
$p>0$, $W$ is a fixed Cohen ring of $k$ and $K$ is the fraction field 
of $W$. 
We fix a $p$-adic fine log formal scheme 
$(\cB,M_{\cB})$ separated and topologically of finite type 
over $\Spf W$ as a base log formal scheme. (Note that $\cB$ is 
Noetherian.) 
We denote the 
scheme $(\cB,M_{\cB}) \otimes_{\Z_p} \Z/p\Z$ by $(B,M_B)$. 
We denote the 
category of fine log (not necessarily $p$-adic) formal schemes 
which are separated and 
topologically of finite type over $(\cB,M_{\cB})$ by $\LF$ and denote 
the full subcategory of $\LF$ which consists of $p$-adic ones by $\pLF$. 
We denote the full subcategory of $\LF$ which consists of usual fine 
log schemes by $\L$ and the full subcategory which consists of 
fine log schemes over $(B,M_B)$ by $\LB$. 
We call an object in $\LF$ 
a fine log formal $\cB$-scheme, an object in $\pLF$ a $p$-adic fine 
log formal $\cB$-scheme, an object in $\L$ a fine log $\cB$-scheme and 
an object in $\LB$ a fine log $B$-scheme. (When log structure is trivial, 
we omit the term `fine log'.) Note that, in Section 5, 
we will impose more assumptions on $\cB$. \\
(2) \,\, For a formal $\cB$-scheme $T$, we denote the rigid analytic space 
assocated to $T$ by $T_K$. \\
(3) \,\, In this paper, we freely use terminologies concerning log structures 
defined in \cite{kkato}, \cite{shiho1} and \cite{shiho2}. 
For a fine log (formal) scheme $(X,M_X)$, $(X,M_X)_{\triv}$ denotes the 
maximal open sub (formal) scheme of $X$ on which the log structure $M_X$ is 
trivial. A morphism $f:(X,M_X) \lra (Y,M_Y)$ is said to be strict if 
$f^*M_Y = M_X$ holds. \\ 
(4) \,\, For a site ${\cal S}$, we denote the topos associated to
${\cal S}$ by ${\cal S}^{\sim}$. \\
(5) \,\, Fiber products of log formal schemes are completed unless 
otherwise stated. On the other hand, the completed tensor product of 
topological modules are denoted by $\hat{\otimes}$ to distinguish 
with the usual tensor product $\otimes$. \\


\section{Relative log crystalline cohomology}

In this section, first we give a review of the definition 
and the basic properties of log crystalline 
site and log crystalline cohomology (which are defined by 
Kato \cite{kkato}) and then we prove the 
finiteness and base change property of log crystalline cohomology when 
the coefficient is a locally free isocrystal. In the case where the 
coefficient is an isocrystal which comes from a locally free crystal, 
corresponding results are deduced from the results of Kato \cite{kkato} 
(the log version of the results of 
Berthelot \cite{berthelot}) fairly easily, but the results here are 
slightly 
stronger than them since a locally free isocrystal does not necessarily 
come from a locally free crystal. After that, we define the relative version 
of HPD-(iso)stratification (which are introduced in \cite[4.3.1]{shiho1} 
in the 
absolute case) and prove a relation with (iso)crystals. \par 
First we fix the notation concerning localization of categories and 
(system of) sheaves on $p$-adic formal schemes. 

\begin{defn}
For an additive category $\cC$, let us define the category 
$\cC_{\Q}$ in the following way$:$ An object of $\cC_{\Q}$ is the same as 
that of $\cC$ and the set of homomorphisms is defined by 
$$ \Hom_{\cC_{\Q}}(X,Y) := \Q \otimes_{\Z} \Hom_{\cC}(X,Y). $$
When we regard an object $X$ in $\cC$ as an object in $\cC_{\Q}$, 
we denote it by $\Q \otimes X$. 
\end{defn} 

\begin{defn}\label{m-coh}
\begin{enumerate}
\item 
Let $T$ be a $p$-adic formal scheme and let $T_n$ be the closed 
subscheme of $T$ defined by $p^n\cO_T$. Then we define the category 
$\M(T)$ as the category of projective systems $(M_n)_n$, where 
$M_n$ is an $\cO_{T_n}$-module of finite presentation satisfying 
$M_{n+1} \otimes_{\cO_{T_{n+1}}} \cO_{T_n} = M_n$. 
\item 
Let $T$ be a Noetherian $p$-adic formal scheme. 
Then 
we denote the category of coherent sheaves of 
${\cal O}_{T}$-modules by $\Coh({\cal O}_{T})$ and
let 
$\Coh(\Q \otimes {\cal O}_{T})$ be 
the category of sheaves of $\Q \otimes_{\Z} {\cal O}_{T}$-modules 
on $T$ which is isomorphic to $\Q \otimes_{\Z} F$
for some $F \in \Coh({\cal O}_{T})$. 
We call 
an object of $\Coh(\Q \otimes {\cal O}_{T})$ an isocoherent sheaf on 
$T$. \\
\end{enumerate}
\end{defn}

For a $p$-adic formal $\cB$-scheme $T$, it is known 
(\cite[10.10.3]{gd}, \cite[1.2]{ogus1}) 
that there are canonical equivalences of categories 
$\M(T) = \Coh(\cO_T), \Coh(\Q \otimes \cO_T) = \Coh (\cO_{T})_{\Q}$. 
So we have $\M(T)_{\Q} = \Coh(\Q \otimes \cO_T)$. 

\begin{rem} 
In \cite[\S 4.3]{shiho1}, we used the category $\Coh(\Q \otimes {\cal O}_T)$ 
for non-Noetherian $T$ at some points, but this is not correct: 
We should have used the category $\M(T)_{\Q}$ instead. For this reason, 
we will correct some of the arguments in \cite[\S 4.3]{shiho1} in 
Definition \ref{defhpdi} -- Proposition \ref{deszar-crys} below. 
\end{rem} 

Now we recall the definition of log crystalline site and 
relative log crystalline cohomology. 
Assume we are given a diagram 
\begin{equation}\label{diag-p}
(X,M_X) \os{f}{\lra} (Y,M_Y) \os{\iota}{\hra} (\cY,M_{\cY}), 
\end{equation}
where $f$ is a morphism in $\LB$, 
 $(\cY,M_{\cY})$ is an object in $\pLF$ and $\iota$ is the exact closed 
immersion defined by $p\cO_{\cY}$. Denote the canonical PD-structure 
on $p\cO_{\cY}$ by $\gamma$. (Note that $\L$ is a subcategory of 
$\pLF$. So we allow the case where $(\cY,M_{\cY})$ is a fine log 
$\cB$-scheme.) 

\begin{defn}\label{defcrys}
With the above notation, we define the log crystalline site 
$(X/\cY)^{\log}_{\crys}$ of $(X,M_X)/(\cY,M_{\cY})$ 
as follows$:$ An object is 
$T := ((U,M_U),(T,M_T),i,\delta)$, where $(U,M_U)$ is a fine log scheme 
strict etale over $(X,M_X)$ and $(T,M_T)$ is a fine log scheme over 
$(\cY,M_{\cY}) \otimes_{\Z_p} \Z/p^n\Z$ for some $n$. 
$i:(U,M_U) \hra (T,M_T)$ is an exact closed immersion over 
$(\cY,M_{\cY})$ and $\delta$ is a PD-structure 
on $\Ker(\cO_T \ra \cO_U)$ which is compatible 
with $\gamma$. The morphism is defined in natural way and the 
coverings are the ones induced from etale coverings of $T$. 
We denote the sheaf on $(X/\cY)^{\log}_{\crys}$ defined by 
$T \mapsto \Gamma(T,\cO_T)$ by ${\cal O}_{X/\cY}$. \par 
We denote the right derived functor $($resp. the $q$-th right derived 
functor$)$ of the functor 
$$ (X/\cY)^{\log,\sim}_{\crys} \lra \cY_{\Zar}^{\sim}; \,\,\,\,
\cF \mapsto (U \mapsto \Gamma((\cX \times_{\cY}U/U)^{\log}_{\crys},\cF)) $$ 
by $Rf_{X/\cY,\crys *}\cF$ $($resp. $R^qf_{X/\cY,\crys *}\cF)$. 
We call $R^qf_{X/\cY,\crys *}\cE$ the $q$-th relative log crystalline 
cohomology of $(X,M_X)/(\cY,M_{\cY})$ with coefficient $\cF$. 
\end{defn}

\begin{rem} 
The $q$-th relative log crystalline cohomology defined above is 
usually called the $q$-th log crystalline cohomology. We put the word 
`relative' just to emphasize that $\cY$ is not necessarily 
equal to $\Spf W$. 
\end{rem} 

\begin{rem}
$Rf_{X/\cY,\crys *}\cF$ (resp. $R^qf_{X/\cY,\crys *}\cF$) is usually 
denoted by $Rf_{X/\cY,*}\cF$ (resp. $R^qf_{X/\cY,*}\cF$). However, 
we decided to put the subscript `$\crys$' in this paper because 
we treat many other relative cohomology theories later in this paper. 
\end{rem} 

\begin{rem}
Let $(X^{(\b)},M_{X^{(\b)}}) \, (\b \in \Lambda)$ 
be a diagram (indexed by a small 
category $\Lambda$) of fine log $B$-schemes over 
$(Y,M_Y)$. Then we can define the log crystalline 
topos $(X^{(\b)}/\cY)^{\log,\sim}_{\crys}$ of 
$(X^{(\b)},M_{X^{(\b)}})/(\cY,M_{\cY})$ 
as the topos associated to the diagram of topoi 
$\{(X^{(\lambda)}/\cY)^{\log,\sim}_{\crys}\}_{\lambda \in \Lambda}$. 
(see \cite{saintdonat}, \cite[V.3.4]{berthelot}.)
\end{rem}

Next we recall the notion of (iso)crystals: 

\begin{defn} 
Let the notation be as in Definition \ref{defcrys}. 
\begin{enumerate}
\item 
A sheaf $\cF$ of ${\cal O}_{X/\cY}$-modules 
on $(X/\cY)^{\log}_{\crys}$ is called a crystal if, for 
any morphism $\varphi:T' \lra T$ in $(X/\cY)^{\log}_{\crys}$, the canonical 
homomorphism of sheaves $\varphi^*\cF_{T} \lra \cF_{T'}$ is an isomorphism. 
$($Here $\cF_T, \cF_{T'}$ denotes the sheaf on $T_{\Zar},T'_{\Zar}$ induced 
by $\cF$.$)$ 
\item 
A crystal $\cF$ is said to be of finite presentation 
$($resp. locally free of finite type$)$ 
if $\cF_T$ is an 
${\cal O}_T$-module of finite presentation 
$($resp. a locally free $\cO_T$-module of finite type$)$ 
for any $T \in (X/\cY)^{\log}_{\crys}$. 
We denote the category of crystals of finite presentation 
on $(X/\cY)^{\log}_{\crys}$ by $C_{\crys}((X/\cY)^{\log})$. 
\item 
A category of isocrystals $I_{\crys}((X/\cY)^{\log})$ on 
$(X/\cY)^{\log}_{\crys}$ is defined to be the category 
$C_{\crys}((X/\cY)^{\log})_{\Q}$. 
\item 
For $\cE = \Q \otimes \cF \in I_{\crys}((X/\cY)^{\log})$, 
we put 
$$ 
Rf_{X/\cY,\crys *}\cE := \Q \otimes_{\Z} Rf_{X/\cY,\crys *}\cF, 
\,\,\, 
R^qf_{X/\cY,\crys *}\cE := \Q \otimes_{\Z} R^qf_{X/\cY,\crys *}\cF $$
and we call $R^qf_{X/\cY,\crys *}\cE$ the $q$-th relative log crystalline 
cohomology of $(X,M_X) \allowbreak / \allowbreak 
(\cY,M_{\cY})$ with coefficient $\cE$. 
\end{enumerate}
\end{defn}

Next we define the notion of local freeness of isocrystals. 
To do this, first we define the notion of local freeness in the category 
$M(T)_{\Q}$. 

\begin{defn}\label{deflf1}
Let $T$ be a $p$-adic formal scheme. Then $M \in \M(T)_{\Q}$ is called 
locally free if $M$ is isomorphic to a direct summand of 
$\Q \otimes ({\cal O}_{T_n}^{\oplus r})_n$ $($where $T_n:=T \otimes_{\Z_p} \Z/p^n\Z)$ 
for some $r \in \N$ Zariski locally on $T$. 
When $T$ is a $p$-adic formal $\cB$-scheme, we say an object $M$ in 
$\Coh(\Q \otimes {\cal O}_T)$ a locally free $\Q\otimes_{\Z}{\cal O}_T$-module 
if the corresponding object in $M(T)_{\Q}$ is locally free. 
\end{defn} 

\begin{rem} 
When $T$ is an affine $p$-adic formal $\cB$-scheme, $M \in 
\Coh(\Q \otimes {\cal O}_T)$ is locally free if and only if it is 
associated to a finitely generated projective $\Q \otimes 
\Gamma(T,{\cal O}_T)$-module. Contrary to the name, a locally free 
$\Q \otimes_{\Z} {\cal O}_T$-module in the sense of Definition 
\ref{deflf1} does {\it not} necessarily has the form 
$\Q \otimes_{\Z} {\cal O}_{T}^{\oplus r}$ Zariski locally on $T$. 
$M \in \Coh(\Q\otimes {\cal O}_T)$ is a locally free 
$\Q \otimes_{\Z} {\cal O}_T$-module if and only if the corresponding 
coherent sheaf on the rigid analytic space $T_K$ 
is a locally free 
${\cal O}_{T_K}$-module with respect to the Grothendieck topology on 
$T_K$. 
\end{rem} 

Next we define the notion of a $p$-adic system in $(X/\cY)^{\log}_{\crys}$ 
as follows: 

\begin{defn} 
Let the notation be as above. A projective system of objects in 
$(X/\cY)^{\log}_{\crys}$ of the form 
$\{((U,M_U),(T_n,M_{T_n}),i_n,\delta_n)\}_n$ is called 
a $p$-adic system in $(X/\cY)^{\log}_{\crys}$ if 
$(T,M_T) \allowbreak := \varinjlim_n(T_n,M_{T_n})$ is a $p$-adic fine log 
formal $\cB$-scheme 
satisfying $(T,M_T) \otimes_{\Z_p} \Z/p^n\Z \allowbreak = (T_n,M_{T_n})$. 
\end{defn} 

Let $\cF$ be a crystal of finite type 
on $(X/\cY)^{\log}_{\crys}$ and let 
$T:=\{((U,M_U),(T_n,M_{T_n}),i_n, \allowbreak \delta_n)\}_n$ 
be a $p$-adic system 
in $(X/\cY)^{\log}_{\crys}$. Then $(\cF_{T_n})_n$ forms an object 
in $\M(T)$, which we denote by $\cF_T$. So, for an isocrystal 
$\cE := \Q \otimes \cF$ on $(X/\cY)^{\log}_{\crys}$ and a $p$-adic 
system $T$ in $(X/\cY)^{\log}_{\crys}$, we can associate the 
object $\Q \otimes \cF_T$ in $\M(T)_{\Q}$, which we denote by $\cE_T$. 

\begin{defn} 
An isocrystal $\cE=\Q \otimes \cF$ is called locally free if, 
for any $p$-adic system 
$T:=\{((U,M_U),(T_n,M_{T_n}),i_n,\delta_n)\}_n$ in $(X/\cY)^{\log}_{\crys}$, 
$\cE_T \in \M(T)_{\Q}$ is locally free in the sense of Definition 
\ref{deflf1}. 
\end{defn} 

\begin{rem} 
In the previous papers \cite{shiho1}, \cite{shiho2}, $M \in 
\Coh(\Q \otimes {\cal O}_T)$ is called locally free if it has the form 
$\Q \otimes_{\Z} {\cal O}_T^{\oplus r}$ Zariski locally on $T$ and an 
isocrystal $\cE$ is called locally free if each ${\cal E}_T$ has the form 
$\Q \otimes_{\Z} {\cal O}_T$ Zariski locally on $T$. 
But these definition are not good for this paper 
because the imposed condition is 
too strong. So we change the definition here. 
\end{rem}

\begin{rem} 
If $\cF$ is a locally free crystal of finite type, then 
$\cE = \Q \otimes \cF$ is a 
locally free isocrystal. However, we do not know if, for a locally 
free isocrystal $\cE$, there exists a locally free crystal of finite type 
$\cF$ satisfying $\cE \cong \Q \otimes \cF$. We are doubtful for this claim. 
\end{rem}

In the rest of this section, for a $p$-adic 
formal log scheme $(S,M_S)$, we denote the 
log scheme $(S,M_S) \otimes_{\Z_p} \Z/p^n\Z$ by $(S_n,M_{S_n})$. Now we 
recall several basic properties of log crystalline site and 
relative log crystalline cohomology. \par 
Assume we are given the diagram \eqref{diag-p} and assume for the moment that 
$(X,M_X)$ admits a closed immersion 
$(X,M_X) \hra (\cP,M_{\cP})$ into a $p$-adic fine log formal $\cB$-scheme 
which is formally log smooth over $(\cY,M_{\cY})$. 
Let $(D,M_D)$ be the $p$-adically completed log PD-envelope of 
$(X,M_X)$ in $(\cP,M_{\cP})$. 
Let $\MC_n$ be 
the category of $\cO_{D_n}$-modules of finite type 
with integrable quasi-nilpotent (in the sense of \cite[p.19]{ogus4}) 
log connection on $(\cP_n,M_{\cP_n})/(\cY_n,M_{\cY_n})$ which is compatible 
with the canonical connection on $\cO_{D_n}$, and 
let $\MC$ be the category of projective systems $\{(M_n,\nabla_n)\}_n$ of 
objects in $\MC_n$ such that $(M_{n+1},\nabla_{n+1}) 
\otimes_{\Z/p^{n+1}\Z} \Z/p^n\Z = (M_n,\nabla_n)$ holds. 
Then we have an equivalence of 
categories (\cite[6.2]{kkato}) of the following form: 
\begin{equation}\label{eqcrysconn}
C_{\crys}((X/\cY)^{\log}) \cong \MC, \,\,\, 
\cF \mapsto \{(\cF_{D_n}, \nabla_n: \cF_{D_n} \ra 
\cF_{D_n} \otimes_{\cO_{\cP_n}} \omega^1_{\cP_n/\cY_n})\}_n. 
\end{equation}
If we have $(D,M_D)=(\cP,M_{\cP})$ (this is the case when we have 
$(\cP,M_{\cP}) \times_{(\cY,M_{\cY})} (Y,M_Y) = (X,M_X)$), 
we have the equivalence of categories 
\begin{equation}\label{eqcrysconn2}
\MC \cong 
\left( 
\begin{aligned}
& \text{coherent sheaves with quasi-nilpotent integrable} \\
& \text{formal log connection on $(\cP,M_{\cP})/(\cY,M_{\cY})$} \\
\end{aligned}
\right), 
\end{equation}
which is given by $\{(M_n,\nabla_n)\}_n \mapsto 
(\varprojlim_n M_n, \varprojlim_n \nabla_n)$. 
Note that $\nabla_n: \cF_{D_n} \ra 
\cF_{D_n} \otimes_{\cO_{\cP_n}} \omega^1_{\cP_n/\cY_n}$ naturally 
extends to the complex of the form 
$\cF_{D_n} \otimes_{\cO_{\cP_n}} \omega^{\b}_{\cP_n/\cY_n}$ which is 
compatible with respect to $n$. 
We define the log de Rham complex $\DR(X/\cY,\cF)$ associated to 
$\cF$ as the complex 
$\vpl_n (\cF_{D_n} \otimes_{\cO_{\cP_n}} \omega^{\b}_{\cP_n/\cY_n})$. 
Then, by \cite[6.4]{kkato} (see also 
\cite[3.1.4]{shiho2}), we have 
$\DR(D/\cY,\cF) = R\ti{u}_* \cF$ 
(where $\ti{u}:(X/\cY)^{\log,\sim}_{\crys} \lra X^{\sim}_{\Zar}$ is 
the projection). Then we have the quasi-isomorphism 
$Rf_{X/\cY,\crys*}\cF = Rf_*\DR(D/\cY,\cF)$. 
For an isocrystal $\cE = \Q \otimes \cF$ on $(X/\cY)^{\log}_{\crys}$, 
we define the log de Rham complex $\DR(D/\cY,\cE)$ by 
$\DR(D/\cY,\cE):= \Q \otimes_{\Z} \DR(D/\cY,\cF)$. 
So we have the quasi-isomorphism 
$Rf_{X/\cY,\crys*}\cE = Rf_*\DR(D/\cY,\cE)$. \par 
Even if $(X,M_X)$ does not admit the closed immersion 
$(X,M_X) \hra (\cP,M_{\cP})$ as above, we always have an 
embedding system (\cite[2.18]{hyodokato})
\begin{equation}\label{emb}
(X,M_X) \os{g^{(\b)}}{\lla} (X^{(\b)},M_{X^{(\b)}}) 
\hra (\cP^{(\b)},M_{\cP^{(\b)}}).
\end{equation} 
For any abelian sheaf $\cF$ in $(X/\cY)^{\log,\sim}_{\crys}$, 
let us denote the 
pull-back of $\cF$ to $(X^{(\b)}/\cY)^{\log,\sim}_{\crys}$ by 
$\cF^{(\b)}$. Then we have the quasi-isomorphism (called cohomological 
descent) 
$\cF \os{\cong}{\lra} Rg^{(\b)}_*\cF^{(\b)}$, where the morphism of 
topoi $(X^{(\b)}/\cY)^{\log,\sim}_{\crys} \lra 
(X/\cY)^{\log,\sim}_{\crys}$ is also 
denoted by $g^{(\b)}$. Hence, for a crystal $\cF$ on $(X/\cY)^{\log}_{\crys}$, 
we have 
\begin{equation}\label{crys-D}
Rf_{X/\cY,\crys *}\cF = Rf_{X^{(\b)}/\cY,\crys *}\cF^{(\b)} = 
R(f \circ g^{(\b)})_* \DR(D^{(\b)}/\cY, \cF^{(\b)}), 
\end{equation}
where $(D^{(\b)},M_{D^{(\b)}})$ denotes the $p$-adically completed 
log PD envelope of $(X^{(\b)},M_{X^{(\b)}})$ in 
$(\cP^{(\b)},M_{\cP^{(\b)}})$. The same formula holds also in the case of 
isocrystals. \par 
We recall another construction which induces the formula like 
\eqref{crys-D}. Assume given a diagram \eqref{diag-p} with $f$ log smooth. 
Then, by \cite[3.14]{kkato}, we have an open covering 
$X = \bigcup_{j \in J}X_j$ by finite number of affine subschemes and 
exact closed immersions $i_j: (X_j,M_{X_j}) := 
(X_j,M_X |_{X_j}) \hra (\cP_j,M_{\cP_j}) \,(j \in J)$ 
into a fine log formal $\cB$-scheme $(\cP_j,M_{\cP_j})$ formally log 
smooth over $(\cY,M_{\cY})$ such that $\cP_j$ is affine and that 
$(\cP_j,M_{\cP_j}) \times_{(\cY,M_{\cY})} (Y,M_Y) = (X_j,M_{X}|_{X_j})$ 
holds. For a non-empty subset $L \subseteq J$, let 
$(X_L,M_{X_L})$ (resp. $(\cP_L,M_{\cP_L})$) be the fiber product of 
$(X_j,M_{X_j})$'s (resp. $(\cP_j,M_{\cP_j})$'s) for $j \in L$ over 
$(X,M_X)$ (resp. $(\cY,M_{\cY})$), and let $i_L:(X_L,M_{X_L}) \hra 
(\cP_L,M_{\cP_L})$ be the closed immersion induced by $i_j$'s ($j \in L$). 
For $m \in \N$, let $(X^{(m)},M_{X^{(m)}})$ (resp. 
$(\cP^{(m)},M_{\cP^{(m)}})$) be the disjoint union of $(X_L,M_{X_L})$'s 
(resp. $(\cP_L,M_{\cP_L})$'s) for $L \subseteq J$ with $|L|=m+1$, and 
let $i^{(m)}: (X^{(m)},M_{X^{(m)}}) \hra (\cP^{(m)},M_{\cP^{(m)}})$ be 
the disjoint union of $i_L$'s with $|L|=m+1$. 
(Note that $X^{(m)}$ is empty for $m \geq |J| =: N+1$.) 
Now let $\Delta^+_N$ be the category 
such that the 
objects are the sets $[m]:=\{0,1,2,\cdots, m\}$ with $0 \leq 
m \leq N$ and that 
$\Hom_{\Delta^+_N}([m],[m'])$ is the set of strictly increasing 
maps $[m] \ra [m']$. Then, if we fix a total order on $J$, 
we can regard $(X^{(\b)},M_{X^{(\b)}}) := 
\{(X^{(m)},M_{X^{(m)}})\}_{0 \leq m \leq N}$ 
(resp. $(\cP^{(\b)},M_{\cP^{(\b)}}) := 
\{(\cP^{(m)},M_{\cP^{(m)}})\}_{0 \leq m \leq N}$) naturally as a diagram of 
fine log $B$-schemes (resp. fine log formal $\cB$-schemes) 
indexed by $\Delta^+_N$ in the following way: To each $[m] \in \Delta^+_N$, 
we associate the fine log $B$-scheme $(X^{(m)}, \allowbreak 
M_{X^{(m)}})$ (resp. 
fine log formal $\cB$-scheme $(\cP^{(m)},M_{\cP^{(m)}})$). 
For a morphism $\alpha:[m'] \ra [m]$ in $\Delta^+_N$ and a subset 
$L \subseteq J$ with $|L|=m+1$, we define the subset $L'_{\alpha} \subseteq L$ 
to be the inverse image of $\alpha([m'])$ by the unique order-preserving 
isomorphism $L \os{\simeq}{\ra} [m]$, where the order on $[m]$ is the 
canonical one. Then we define the morphism 
$\alpha^*: (X^{(m)},M_{X^{(m)}}) \lra (X^{(m')},M_{X^{(m')}})$ 
(resp. $\alpha^*: (\cP^{(m)},M_{\cP^{(m)}}) \lra (\cP^{(m')},M_{\cP^{(m')}})$)
 as the composite 
$$ 
(X^{(m)},M_{X^{(m)}}) = 
\coprod_{L \subseteq J,|L|=m+1} (X_L,M_{X_L}) 
\os{\coprod \pr_{L'_{\alpha} \subseteq L}}{\lra} 
\coprod_{L \subseteq J,|L|=m'+1} (X_L,M_{X_L}) = (X^{(m')},M_{X^{(m')}})
$$
\begin{align*}
\text{(resp.} \,\, 
(\cP^{(m)},M_{\cP^{(m)}}) & = 
\coprod_{L \subseteq J,|L|=m+1} (\cP_L,M_{\cP_L}) \\ 
& \os{\coprod \pr_{L'_{\alpha} \subseteq L}}{\lra} 
\coprod_{L \subseteq J,|L|=m'+1} (\cP_L,M_{\cP_L}) = 
(\cP^{(m')},M_{\cP^{(m')}}) \,\text{),}
\end{align*} 
where 
$$ \pr_{L'_{\alpha}\subseteq L}: 
(X_L,M_{X_L}) = \prod_{j \in L, (X,M_X)}(X_j,M_{X_j}) \lra 
\prod_{j \in L'_{\alpha},(X,M_X)}(X_j,M_{X_j}) =(X_{L'_{\alpha}}, 
M_{X_{L'_{\alpha}}}) $$  
$$\text{(resp.} \,\, \pr_{L'_{\alpha}\subseteq L}: 
(\cP_L,M_{\cP_L}) = \prod_{j \in L, (\cY,M_{\cY})}(\cP_j,M_{\cP_j}) \lra 
\prod_{j \in L'_{\alpha},(\cY,M_{\cY})}(\cP_j,M_{\cP_j}) 
= (\cP_{L'_{\alpha}}, M_{\cP_{L'_{\alpha}}}) 
 \text{)}$$ 
is the 
projection associated to the inclusion $L'_{\alpha} \subseteq L$. 
Then we have a 
canonical diagram 
\begin{equation}\label{emb2}
(X,M_X) \os{g^{(\b)}}{\lla} (X^{(\b)},M_{X^{(\b)}}) 
\hra (\cP^{(\b)},M_{\cP^{(\b)}}) 
\end{equation} 
over $(\cY,M_{\cY})$ 
and by \cite[V 3.4.8]{berthelot}, we have again the formula 
\begin{equation}\label{crys-D2}
Rf_{X/\cY,\crys *}\cF = Rf_{X^{(\b)}/\cY,\crys *}\cF^{(\b)} = 
R(f \circ g^{(\b)})_* \DR(D^{(\b)}/\cY, \cF^{(\b)}), 
\end{equation}
where the notations are as in \eqref{crys-D}. The formula 
\eqref{crys-D2} has the advantage that the right hand side has only finitely 
many terms. \par 
Keep the situation of the diagram \ref{diag-p} with $f$ log smooth. 
Then, locally on $X$, we have an exact closed immersion 
$(X,M_X) \hra (\cP,M_{\cP})$ into a $p$-adic fine log formal $\cB$-scheme 
$(\cP,M_{\cP})$ formally smooth over $(\cY,M_{\cY})$ satisfying 
$(\cP,M_{\cP}) \times_{(\cY,M_{\cY})} (Y,M_Y) = (X,M_X)$. 
Then we have equivalences of categories \eqref{eqcrysconn}, 
\eqref{eqcrysconn2}. Hence $C_{\crys}((X/\cY)^{\log})$ is an abelian 
category locally on $X$, because so is the right hand side in 
\eqref{eqcrysconn2}. Since $C_{\crys}((X/\cY)^{\log})$ satisfies the 
descent property for Zariski covering of $X$, this implies that 
$C_{\crys}((X/\cY)^{\log})$ is an abelian category globally. \par 
Now we prove a finiteness property of relative log crystalline cohomology 
when the coefficient is a locally free isocrystal: 

\begin{thm}\label{finicrys}
Assume given a diagram \eqref{diag-p} and assume moreover that $f$ is 
proper and log smooth. Then, for an isocrystal $\cE$, 
the relative log crystalline cohomology $R^qf_{X/\cY,\crys *}\cE$ 
is an isocoherent sheaf on $\cY$ for any $q \in \N$ and it is zero 
for sufficiently large $q$. 
\end{thm} 

\begin{thm}\label{perfcrys}
With the above notation, assume that $\cE$ is a locally free isocrystal and 
assume moreover that either $f$ is integral or 
$Y$ is regular. Then $Rf_{X/\cY,\crys*}\cE$ 
is a perfect complex of $\Q \otimes_{\Z} \cO_{\cY}$-modules, that is, it is 
quasi-isomorphic to a bounded complex of locally free 
$\Q \otimes_{\Z} \cO_{\cY}$-modules $($in the sense of 
Definition \ref{deflf1}$)$ Zariski locally on $\cY$. 
\end{thm} 

\begin{rem}
When $\cE$ has the form $\Q \otimes \cF$ for a locally free crystal 
$\cF$ of finite type, Theorems \ref{finicrys} and \ref{perfcrys} are 
known by \cite[7.24.3]{berthelotogus}. 
\end{rem} 

\begin{pf*}{Proof of Theorem \ref{finicrys}}
To prove the theorem, 
we may assume that $\cY$ is affine. 
Let $\cY'$ be the closed subscheme 
of $\cY$ defined by the ideal $\{x \in \cO_{\cY} \,\vert\, p^nx=0 
\,\text{for some $n>0$}\}$, let $(X',M_{X'}) \lra (Y',M_{Y'}) 
\hra (\cY',M_{\cY}|_{\cY'})$ be the base change of the diagram 
\eqref{diag-p} by the exact closed immersion $(\cY',M_{\cY}|_{\cY'}) \hra 
(\cY,M_{\cY})$ and let $\cE'$ be the pull-back of $\cE$ to 
$(X'/\cY')^{\log}_{\crys}$. 
Then, by using a variant of \eqref{crys-D} for 
isocrystals, we can see the quasi-isomorphism 
$Rf_{X/\cY,\crys *}\cE = i_*Rf_{X'/\cY',\crys *}\cE'$, where 
$i$ denotes the closed immersion $\cY' \hra \cY$. Noting the 
fact that $i_*$ induces 
the equivalence of categories $\Coh(\Q \otimes \cO_{\cY'}) \cong 
\Coh(\Q \otimes \cO_{\cY})$, we see that we may replace $\cY$ by 
$\cY'$ to prove the theorem. So we may assume that 
${\cal O}_{\cY}$ has no $p$-torsion. \par 
Let us take a diagram 
$$ (X,M_X) \os{g^{(\b)}}{\lla} (X^{(\b)},M_{X^{(\b)}}) 
\hra (\cP^{(\b)},M_{\cP^{(\b)}}) $$
as \eqref{emb2} and take 
a crystal $\cF$ of finite presentation on $(X/\cY)^{\log}_{\crys}$ satisfying 
$\cE = \Q \otimes \cF$. Then $\cF$ induces an object 
$\cF_{\cP^{(0)}}$ in $\M(\cP^{(0)}) = \Coh(\cO_{\cP^{(0)}})$. Since 
$\cP^{(0)}$ is Noetherian, there exists a natural number $a$ satisfying 
$$ \{x \in \cF_{\cP^{(0)}} \,\vert\, p^nx=0\,\text{for some $n$}\} = 
\{x \in \cF_{\cP^{(0)}} \,\vert\, p^ax=0\}. $$ 
Then let us put $\cF' := \Coker(\Ker(p^a:\cF \ra \cF) \hra \cF)$, where 
cokernel and kernel are taken in the category 
$C_{\crys}((X/\cY)^{\log})$. Then, by using the equivalences of 
categories \eqref{eqcrysconn} and \eqref{eqcrysconn2}, we have 
$$ \cF'_{\cP^{(0)}} = 
\cF_{\cP^{(0)}}/\{x \in \cF_{\cP^{(0)}} \,\vert\, p^ax=0\} = 
\cF_{\cP^{(0)}}/\{x \in \cF_{\cP^{(0)}} 
\,\vert\, p^nx=0\,\text{for some $n$}\}.$$ 
Since we have $\Q \otimes \cF \cong \Q \otimes \cF'$, we may replace 
$\cF$ by $\cF'$, that is, we may assume that $\cF_{\cP^{(0)}}$ has 
no $p$-torsion. \par 
Let us denote the restriction of $\cF$ to 
$(X/\cY_n)^{\log,\sim}_{\crys}$ by $\cF_n$ and the restriction of 
$\cF_n$ to $(X^{(\b)}/\cY_n)^{\log,\sim}_{\crys}$ by $\cF^{(\b)}_n$. 
For $m \in \N$, take the open immersion 
$\cP^{(0,m)} \hra \cP^{(0)}$ of fine log formal $\cB$-schemes 
satisfying 
$\cP^{(0,m)} \times_{\cP^{(0)}} X^{(0)} = X^{(m)}$. Then 
$\cF_n^{(m)}$ induces the log de Rham complex 
$\DR(\cP^{(0,m)}_n/\cY_n,\cF^{(m)}_n)$. It is 
compatible with respect to $n$ and we have the quasi-isomorphism 
{\small{ \begin{equation}\label{drcrys}
Rf_{X^{(m)}/\cY_n,\crys *}\cF^{(m)}_n = 
R(f \circ g^{(m)})_{*}\DR(\cP^{(0,m)}_n/\cY_n,\cF^{(m)}_n) 
= (f \circ g^{(m)})_{*}\DR(\cP^{(0,m)}_n/\cY_n,\cF^{(m)}_n) 
\end{equation}}}
which is also compatible with respect to $n$. (Note however that 
the above quasi-isomorphism is not compatible with respect to $m$.) \par 
Now we prove the following claim: \\ 
\quad \\
{\bf claim.} \,\, Let the notation be as above. Then: \\ 
(1) \,\, We have the canonical quasi-isomorphism 
$$(Rf_{X/\cY_n,\crys *}\cF_n) \otimes^L_{\cO_{\cY_n}} \cO_{\cY_{n-1}} 
\os{\cong}{\lra} Rf_{X/\cY_{n-1},\crys *}\cF_{n-1}.$$ 
(2) \,\, We have the distinguished triangle 
$$ 
Rf_{X/Y,\crys *}\cF \os{p^{n-1}}{\lra} Rf_{X/\cY_n,\crys *}\cF 
\lra Rf_{X/\cY_{n-1}, \crys *}\cF_{n-1} \os{+1}{\lra}. $$ 

\begin{pf*}{Proof of claim}
First we prove (1). 
The construction of the map is standard and we omit it. To prove that the 
map is a quasi-isomorphism, we may replace $X, \cF$ by $X^{(\b)}, \cF^{(\b)}$ 
by the first equality in \eqref{crys-D2}
and then by $X^{(m)}, \cF^{(m)}$. So, by \eqref{drcrys}, 
we are reduced to showing that the map 
{\small{\begin{equation}\label{map1}
(f \circ g^{(m)})_{*}\DR(\cP^{(0,m)}_n/\cY_n,\cF^{(m)}_n) 
\otimes^L_{\cO_{\cY_n}} 
\cO_{\cY_{n-1}} \lra 
(f \circ g^{(m)})_{*}\DR(\cP^{(0,m)}_{n-1}/\cY_{n-1},\cF^{(m)}_{n-1})
\end{equation}}}
is a quasi-isomorphism. Since ${\cal O}_{\cY}$ 
has no $p$-torsion, it is flat over 
$\Z_p$. So $\cO_{\cY_n}$ is flat over $\Z/p^n\Z$ and it implies the 
quasi-isomorphism 
{\small{
$$(f \circ g^{(m)})_{*}\DR(\cP^{(0,m)}_n/\cY_n,\cF^{(m)}_n) 
\otimes^L_{\cO_{\cY_n}} \cO_{\cY_{n-1}} = 
(f \circ g^{(m)})_{*}\DR(\cP^{(0,m)}_n/\cY_n,\cF^{(m)}_n) 
\otimes^L_{\Z/p^n\Z} \Z/p^{n-1}\Z.$$}}
Next, since $\cF_{\cP^{(0)}}$ has no $p$-torsion, 
$(f \circ g^{(m)})_*\cF^{(m)}_{\cP^{(0,m)}_n}$ 
is flat over $\Z/p^n\Z$. So the each term of 
$(f \circ g^{(m)})_{*}\DR(\cP^{(0,m)}_n/\cY_n,\cF_n)$ is 
flat over $\Z/p^n\Z$ and then we have 
\begin{align*}
& \phantom{=} 
(f \circ g^{(m)})_{*}\DR(\cP^{(0,m)}_n/\cY_n,\cF^{(m)}_n) \otimes^L_{\Z/p^n\Z} 
\Z/p^{n-1}\Z \\ & = 
(f \circ g^{(m)})_{*}\DR(\cP^{(0,m)}_n/\cY_n,\cF^{(m)}_n) \otimes_{\Z/p^n\Z} 
\Z/p^{n-1}\Z \\ & = 
(f \circ g^{(m)})_{*}\DR(\cP^{(0,m)}_{n-1}/\cY_{n-1},\cF^{(m)}_{n-1}). 
\end{align*}
So the map \eqref{map1} is a quasi-isomorphism and we have proved 
the assertion (1). \par 
Next we prove (2). We can construct the first two maps easily, by using 
\eqref{crys-D}. To prove the assertion, we may replace 
$X, \cF$ by $X^{(m)},\cF^{(m)}$ as in the proof of (1). Then we are 
reduced to showing the existence of the distinguished triangle of the form 
\begin{align*}
(f \circ g^{(m)})_*\DR(\cP^{(0,m)}_1/Y,\cF^{(m)}_1) & \os{p^{n-1}}{\lra}
(f \circ g^{(m)})_*\DR(\cP^{(0,m)}_n/\cY_n,\cF^{(m)}_n) \\ & \lra 
(f \circ g^{(m)})_*\DR(\cP^{(0,m)}_{n-1}/\cY_{n-1},\cF^{(m)}_{n-1}) 
\os{+1}{\lra}. 
\end{align*}
We can construct this trangle by applying 
$(f \circ g^{(m)})_*\DR(\cP^{(0,m)}_n/\cY_{n},\cF^{(m)}_{n}) 
\allowbreak \otimes^L_{\Z/p^n\Z}$ to 
the triangle $\Z/p\Z \os{p^{n-1}}{\lra} \Z/p^n\Z \lra 
\Z/p^{n-1}\Z \os{+1}{\lra}$, 
since $(f \circ g^{(m)})_*\cF^{(m)}_{\cP^{(0,m)}_n}$ is flat over 
$\Z/p^n\Z$. So we have proved the assertion (2). 
\end{pf*}

By claim (2) and induction, we see that $Rf_{X/\cY_n,\crys *}\cF_n$ 
is bounded and it has 
finitely generated cohomologies for any $n$. 
(Note that, in the case $n=1$, we have $Rf_{X/Y,\crys *}\cF = 
Rf_*\DR(X/Y, \cF_X)$. So it is bounded and it has 
finitely generated cohomologies since $f$ is proper.) 
By this and claim (1), we 
see that the family $\{Rf_{X/\cY_n,\crys *}\cF_n\}_n$ forms 
a consistent system in the sense of \cite[B.4]{berthelotogus}. 
Hence, by \cite[B.9]{berthelotogus}, we see that 
$Rf_{X/\cY,\crys *}\cF = R\vpl Rf_{X/\cY_n,\crys *}\cF_n$ 
is bounded above (hence bounded) and has 
finitely generated cohomologies. 
Hence $Rf_{X/\cY,\crys *}\cE$ is bounded and has isocoherent cohomologies. 
So we are done. 
\end{pf*}

\begin{pf*}{Proof of Theorem \ref{perfcrys}}
First let us consider the case $Y$ is regular. 
We may assume that $Y$ is affine. 
In this case, 
$Rf_{X/Y,\crys *}\cF = Rf_*\DR(X/Y, \cF_X)$ is a bounded complex of 
quasi-coherent $\cO_Y$-modules with coherent cohomologies. So it is 
quasi-isomorphic to a bounded complex of coherent $\cO_Y$-modules. 
Now note that, since $Y$ is affine and regular, 
every coherent $\cO_Y$-module has 
a resolution of finite length by locally free $\cO_Y$-module of finite rank. 
So $Rf_{X/Y,\crys *}\cF$ is quasi-isomorphic to a bounded complex of 
locally free $\cO_{Y}$-modules of finite rank, that is, it is 
a perfect complex of $\cO_{Y}$-modules. From this fact and 
\cite[B.10]{berthelotogus}, we see that $Rf_{X/\cY,\crys *}\cF$ is 
a perfect complex of $\cO_{\cY}$-modules. So we are done. \par 
Next let us consider the case where $f$ is integral. 
By Theorem 
\ref{finicrys}, we know that $Rf_{X/\cY,\crys *}\cE$ is bounded and it 
has isocoherent cohomologies. So it is quasi-isomorphic to a bounded complex 
of isocoherent sheaves locally on $\cY$. Hence, to prove that 
$Rf_{X/\cY,\crys *}\cE$ is a perfect complex of 
$\Q \otimes_{\Z} \cO_{\cY}$-modules, it suffices to prove that 
$Rf_{X/\cY,\crys *}\cE$ is 
quasi-isomorphic to a bounded complex of flat 
$\Q \otimes_{\Z} \cO_{\cY}$-modules. \par 
Let the notation be as in the proof 
of Theorem \ref{finicrys} and put $\cE^{(\b)} := \Q \otimes \cF^{(\b)}$. 
Then, to prove the above assertion, we may replace $X,\cE$ by 
$X^{(m)}, \cE^{(m)}$. Then, since we have 
$Rf_{X^{(m)}/\cY,\crys *}\cE^{(m)} := 
\Q \otimes_{\Z} 
(f \circ g^{(m)})_*\DR(\cP^{(0,m)}/\cY, \cF^{(m)}) = 
(f \circ g^{(m)})_*(\Q \otimes_{\Z} \DR(\cP^{(0,m)}/\cY, \cF^{(m)}))$, 
it suffices to prove 
that 
$(f \circ g^{(m)})_* (\Q \otimes_{\Z} \cF^{(m)}_{\cP^{(0,m)}} \otimes 
\omega^i_{\cP^{(0,m)}/\cY}) \, (i \in \N)$ are flat 
$\Q \otimes_{\Z} \cO_{\cY}$-modules. Now let us note that 
$(\cP^{(0,m)},M_{\cP^{(0,m)}}) \lra (\cY,M_{\cY})$ is formally log smooth 
by definition and integral because 
 so is $f$. So $\cP^{(0,m)}$ is flat over $\cY$. Since 
$\Q \otimes_{\Z} \cF^{(m)}_{\cP^{(0,m)}} \otimes 
\omega^i_{\cP^{(0,m)}/\cY} \, (i \in \N)$ are locally free 
$\Q \otimes_{\Z} \cO_{\cP^{(0,m)}}$-modules, we can deduce that 
$(f \circ g^{(m)})_* (\Q \otimes_{\Z} \cF^{(m)}_{\cP^{(0,m)}} \otimes 
\omega^i_{\cP^{(0,m)}/\cY}) \, (i \in \N)$ are flat 
$\Q \otimes_{\Z} \cO_{\cY}$-modules. So we are done. 
\end{pf*}

\begin{rem} 
In the case 
$\cY = \Spf W$, 
Theorem \ref{finicrys} implies the finite-dimensionality 
of the (absolute) log crystalline cohomology $H^i((X/W)^{\log}_{\crys}, \cE)$. 
It is already used in \cite{ogus2} and \cite[3.1.5]{shiho2}, but we could not 
find a precise reference before. 
\end{rem} 

Next we prove the base change property. 

\begin{thm}\label{crysbc}
Assume we are given a diagram 
\begin{equation}
\begin{CD}
(X',M_{X'}) @>>> (Y',M_{Y'}) @>>> (\cY',M_{\cY'}) \\
@VVV @VVV @V{\varphi}VV \\
(X,M_X) @>f>> (Y,M_Y) @>{\iota}>> (\cY,M_{\cY}), 
\end{CD}
\end{equation}
where 
$f$ is proper 
log smooth integral, $\iota$ is the exact closed immersion defined by 
the ideal sheaf $p\cO_{\cY}$ and the squares are Cartesian. 
Then, for a locally free isocrystal $\cE$ on $(X/\cY)^{\log}_{\crys}$, 
we have the quasi-isomorphism 
$$ 
L\varphi^*
Rf_{X/\cY,\crys *}\cE \os{\sim}{\lra} 
Rf_{X'/\cY',\crys *} \varphi^*\cE. $$
\end{thm}

\begin{rem} 
When $\cE$ has the form $\Q \otimes \cF$ for a locally free crystal of 
finite type $\cF$, Theorem \ref{crysbc} is essentially known by 
\cite[7.8]{berthelotogus} plus limit argument. 
\end{rem}

\begin{pf} 
First, by the same argument as the proof of Theorem \ref{finicrys}, 
we may assume that $\cY,\cY'$ are affine and have no $p$-torsion. 
Let us take a diagram 
$$ (X,M_X) \os{g^{(\b)}}{\lla} (X^{(\b)},M_{X^{(\b)}}) \os{i^{(\b)}}{\hra} 
(\cP^{(\b)},M_{\cP^{(\b)}}) $$ 
as in \eqref{emb2} and 
take a crystal $\cF$ of finite presentation
 on $(X/\cY)^{\log}_{\crys}$ satisfying 
$\cE = \Q \otimes \cF$ such that $\cF_{\cP^{(0)}}$ has no $p$-torsion. 
(It is possible by the proof of Theorem \ref{finicrys}.) 
Let $\cF_n, \cF^{(\b)}_n$ be as in the proof of Theorem \ref{finicrys} and 
denote the map 
$\varphi \otimes_{\Z_p} \Z/p^n\Z: (\cY'_n,M_{\cY'_n}) \lra (\cY_n,M_{\cY_n})$ 
simply by $\varphi_n$. \par 
The map 
$L\varphi_{\b}^*Rf_{X/\cY_{\b},\crys *}\cF_{\b} \lra Rf_{X'/\cY'_{\b},\crys *} 
\varphi_{\b}^* \cF_{\b}$ is defined in a standard way 
(see \cite[7.8]{berthelotogus}, 
\cite[V 3.5.2]{berthelot}). This map induces the map 
\begin{equation}\label{aug1}
\Q \otimes_{\Z} (R \vpl L\varphi_{\b}^*Rf_{X/\cY_{\b},\crys *}\cF_{\b}) 
\os{\cong}{\lra} \Q \otimes_{\Z} R\vpl R f_{X'/\cY'_{\b}}\varphi_{\b}^*
\cF_{\b}. 
\end{equation} 
(Note that $L\varphi_{\b}^*Rf_{X/\cY_{\b},\crys *}\cF_{\b}$ is 
quasi-isomorphic to the complex of quasi-coherent 
${\cal O}_{\cY_{\b}}$-modules. So it admits a $\vpl$-acyclic resolution and so 
we can apply the functor $R \vpl$ to it.) The right hand side
 of \eqref{aug1} is quasi-isomorphic to 
to 
$\Q \otimes_{\Z} Rf_{X'/\cY',\crys *}\varphi^*\cF = 
Rf_{X'/\cY',\crys *}\varphi^*\cE$ by \cite[7.22]{berthelotogus} and 
the left hand side is calculated as follows: 
\begin{align*}
& \phantom{==} \Q \otimes_{\Z} (R \vpl 
L\varphi_{\b}^*Rf_{X/\cY_{\b},\crys *}\cF_{\b})
 \\ & = 
\Q \otimes_{\Z} (R \vpl L\varphi_{\b}^*(\cO_{\cY_{\b}} \otimes^L_{\cO_{\cY}} 
Rf_{X/\cY,\crys *}\cF)) \qquad \text{(\cite[B.5]{berthelotogus})} \\ 
& = 
\Q \otimes_{\Z} (R \vpl (\cO_{\cY'_{\b}} \otimes^L_{\cO_{\cY'}} 
L\varphi^*Rf_{X/\cY,\crys *}\cF)) \\ 
& = 
\Q \otimes_{\Z} (L\varphi^*Rf_{X/\cY,\crys *}\cF) 
\qquad \text{(\cite[B.9]{berthelotogus})} \\ 
& = L\varphi^*Rf_{X/\cY,\crys *}\cE. 
\end{align*}
So it suffices to show that the map \eqref{aug1} is a quasi-isomorphism. 
If we denote the cone 
$\Cone (L\varphi_{\b}^*Rf_{X/\cY_{\b},\crys *}\cF_{\b} \allowbreak 
\lra \allowbreak Rf_{X'/\cY'_{\b},\crys *} 
\varphi_{\b}^* \cF_{\b})$ by $C_{\b}$, it suffices to prove the 
equality $\Q \otimes_{\Z} R\vpl C_{\b} = 0$ in derived category. \par 
Let 
$({X'}^{(\b)}, M_{{X'}^{(\b)}}) := 
(X^{(\b)},M_{X^{(\b)}}) \times_{(X,M_X)} (X',M_{X'})$, 
$({\cP'}^{(\b)}, M_{{\cP'}^{(\b)}}) \allowbreak := \allowbreak 
(\cP^{(\b)},M_{\cP^{(\b)}}) \allowbreak 
\times_{(\cY,M_{\cY})} (\cY',M_{\cY'})$ and 
${g'}^{(\b)}:= g^{(\b)} \times_{(X,M_X)} (X',M_{X'})$. 
Then we have 
$$ C_{\b} = 
\Cone (L\varphi_{\b}^*Rf_{X^{(\b)}/\cY_{\b},\crys *}\cF^{(\b)}_{\b} 
\lra Rf_{{X'}^{(\b)}/\cY'_{\b},\crys *} \varphi_{\b}^* 
\cF^{(\b)}_{\b}).$$ 
For $m \in \N$, let us put 
$C_{\b}^{(m)} := 
\Cone (L\varphi_{\b}^*Rf_{X^{(m)}/\cY_n,\crys *}\cF^{(m)}_{\b} 
\lra Rf_{{X'}^{(m)}/\cY'_{\b},\crys *} \varphi_{\b}^* 
\cF^{(m)}_{\b})$. To prove the 
equality $\Q \otimes_{\Z} R\vpl C_{\b} = 0$, it suffices to prove 
the equality $\Q \otimes_{\Z} R\vpl C^{(m)}_{\b} = 0$ for all $m \in \N$, 
because we have $C^{(m)}_{\b}=0$ for sufficiently large $m$. 
Let us take a $p$-adic fine log formal 
$\cB$-scheme $(\cP^{(0,m)},M_{\cP^{(0,m)}})$ as in the proof of 
Theorem \ref{finicrys} and put 
$({\cP'}^{(0,m)},M_{{\cP'}^{(0,m)}}) := (\cP^{(0,m)},M_{\cP^{(0,m)}}) 
\times_{(\cY,M_{\cY})} (\cY',M_{\cY'}).$ 
To prove the claim, we may assume moreover that 
$\cP^{(0,m)}$ is affine, $\Q \otimes \cF_{\cP^{(0,m)}}$ is a free 
$\Q \otimes \cO_{\cP^{(0,m)}}$-module of finite rank 
and $\omega^1_{\cP^{(0,m)}/\cY}$ is a free $\cO_{\cP^{(0,m)}}$-module. Then 
we have 
{\scriptsize{
\begin{align*}
C_n^{(m)} & \cong \Cone(
L{\varphi_n}^*(f \circ g^{(m)})_*\DR(\cP^{(0,m)}_n/\cY_n, 
\cF_n^{(m)}) 
\lra 
(f \circ {g'}^{(m)})_*\DR({\cP'}^{(0,m)}/\cY'_n, 
{\varphi_n}^*\cF_n^{(m)})) \\ 
& \cong \Cone(
L{\varphi_n}^*(f \circ g^{(m)})_*\DR(\cP^{(0,m)}_n/\cY_n, 
\cF_n^{(m)}) \lra 
{\varphi_n}^*(f \circ g^{(m)})_*\DR(\cP^{(0,m)}_n/\cY_n, 
\cF_n^{(m)})). 
\end{align*}}}
(The second isomorphism follows from affine base change.) \par 
For a coherent $\cO_{\cP^{(0,m)}}$-module $\cM$, let us put 
$\cM_n := \cM \otimes {\Z_p} \Z/p^n\Z$ and put 
$$ 
C_{\b}^{(m)}(\cM) := \Cone(
L{\varphi_{\b}}^*(f \circ g^{(m)})_* \cM_{\b} \lra 
{\varphi_{\b}}^*(f \circ g^{(m)})_*\cM_{\b}). 
$$ 
Then, to prove the equality 
$\Q \otimes_{\Z} R\vpl C^{(m)}_{\b} = 0$, it suffices to prove 
the equality $\Q \otimes_{\Z} R\vpl C^{(m)}_{\b}(\cM) = 0$ for 
any $p$-torsion free coherent $\cO_{\cP^{(0,m)}}$-module 
$\cM$ such that $\Q \otimes_{\Z} \cM$ is locally free 
(in the sense of Definition \ref{deflf1}), because 
$\DR(\cP^{(0,m)}/\cY, \cF^{(m)})$ is a bounded complex consisting of 
such modules and 
we have 
$\DR(\cP^{(0,m)}_n/\cY_n, \cF_n^{(m)}) = 
\DR(\cP^{(0,m)}/\cY, \cF^{(m)}) \otimes_{\cO_{\cY}} \cO_{\cY_n}$. \par 
Let us fix a $p$-torsion free coherent $\cO_{\cP^{(0,m)}}$-module 
$\cM$ such that $\Q \otimes_{\Z} \cM$ is locally free and take a 
resolution $\cN^{\b} \lra \cM \lra 0$ of $\cM$ by finitely generated free 
$\cO_{\cP^{(0,m)}}$-modules. (It is possible since $\cP^{(0,m)}$ is 
affine.) Then, since $\cM$ is $p$-torsion free, the induced diagram 
$\cN^{\b}_n \lra \cM_n \lra 0$ gives a resolution of $\cM_n$ by 
finitely generated free $\cO_{\cP^{(0,m)}_n}$-modules. So 
$$ (f \circ g^{(m)})_*\cN^{\b}_n \lra (f \circ g^{(m)})_*\cM_n \lra 0 $$ 
is exact and each $(f \circ g^{(m)})_*\cN^q_n$ is a flat $\cO_{\cY_n}$-module. 
Applying $\varphi^*_n$ and noting the equality 
$\varphi^*_n (f \circ g^{(m)})_* = 
(f\circ {g'}^{(m)})_* {\varphi'_n}^*$, we obtain the diagram 
$$ (f\circ {g'}^{(m)})_* {\varphi'_n}^*\cN_n^{\b} \lra 
(f\circ {g'}^{(m)})_* {\varphi'_n}^*\cM_n $$ 
and the cone of it is a representative of $C_n^{(m)}(\cM)$, which is 
compatible with respect to $n$. 
So we have 
\begin{align*}
& \Q \otimes_{\Z} R\vpl C_{\b}^{(m)}(\cM) 
\\ = & 
\Q \otimes_{\Z} R\vpl 
\Cone((f\circ {g'}^{(m)})_* {\varphi'_{\b}}^*\cN_{\b}^{\b} \lra 
(f\circ {g'}^{(m)})_* {\varphi'_{\b}}^*\cM_{\b}) 
\\ = & 
\Q \otimes_{\Z} \vpl_n 
\Cone((f\circ {g'}^{(m)})_* {\varphi'_n}^*\cN_n^{\b} \lra 
(f\circ {g'}^{(m)})_* {\varphi'_n}^*\cM_n) 
\\ = & 
\Q \otimes_{\Z} (f\circ {g'}^{(m)})_* \Cone({\varphi'}^*\cN^{\b} \lra 
{\varphi'}^*\cM) 
\\ = & 
(f\circ {g'}^{(m)})_*\Cone({\varphi'}^*(\Q \otimes_{\Z}\cN^{\b}) \lra 
{\varphi'}^*(\Q \otimes_{\Z} \cM)) = 0. 
\end{align*}
(The second equality comes from the quasi-coherence and the surjectivity 
of the transition maps for each term and the final equality comes from 
the local freeness of $\Q \otimes_{\Z} \cM$.) 
So we have proved the equality 
$\Q \otimes_{\Z} R\vpl C_{\b}^{(m)}(\cM)=0$ and so we are done. 
\end{pf}

In the rest of this section, we introduce the relative version 
of HPD-(iso)stratification (\cite{ogus2}, \cite[4.3]{shiho1}) 
and prove a relation with (iso)crystals. \par 
First we define the notion of HPD-(iso)stratification 
(cf. \cite[4.3.1]{shiho1}): 

\begin{defn}\label{defhpdi}
Let $(X,M_X) \os{f}{\lra} (Y,M_Y) \os{\iota}{\hra} (\cY,M_{\cY})$ be as 
in \eqref{diag-p} and let $(X,M_X) \hra (\cP,M_{\cP})$ be a closed immersion 
into a $p$-adic fine log formal $\cB$-scheme over $(\cY,M_{\cY})$. 
Let $D(i) \, (i=0,1,2)$ be the $p$-adically completed log PD-envelope of 
$(X,M_X)$ in the $(i+1)$-fold fiber product of $(\cP,M_{\cP})$ over 
$(\cY,M_{\cY})$. Then we have the projections 
$p_i:D(1) \lra D(0) \,(i=1,2)$, $p_{ij}:D(2) \lra D(1)\,(1 \leq i < j \leq 3)$
 and the diagonal map $\Delta:D(0) \lra D(1)$. Then we define an 
HPD-strafinication $($resp. HPD-isostratification$)$ on an object 
$E$ in $\M(D(0))$ $($resp. $\M(D(0))_{\Q})$ as an isomorphism 
$\epsilon: p_2^*E \os{\cong}{\lra} p_1^*E$ satisfying 
$\Delta(\epsilon)=\id$ and 
$p_{12}^*(\epsilon) \circ p_{23}^*(\epsilon) = p_{13}^*(\epsilon)$. 
We denote the category of objects in $\M(D(0))$ $($resp. 
$\M(D(0))_{\Q})$ endowed with HPD-stratification $($resp. 
HPD-isostratification$)$ by 
$\HPD((X \hra \cP/\cY)^{\log})$ $($resp. 
$\HPDI((X \hra \cP/\cY)^{\log}))$. 
\end{defn} 

Let the notations be as above and assume that 
$(\cP,M_{\cP})$ is formally log smooth over $(\cY,M_{\cY})$. 
Then we have the canonical equivalence of categories (\cite[\S 6]{kkato})
$$ C_{\crys}((X/\cY)^{\log}) \lra \HPD((X \hra \cP/\cY)^{\log}) $$
and it naturally induces the fully-faithful functor 
$$ \Lambda:I_{\crys}((X/\cY)^{\log}) \lra \HPDI((X \hra \cP/\cY)^{\log}).$$
As in the absolute case (\cite[4.3.2]{shiho1}), we have the following: 

\begin{prop}\label{crys-hpdi}
If $(\cP,M_{\cP})$ is formally log smooth over $(\cY,M_{\cY})$ and 
$(\cP,M_{\cP}) \allowbreak 
\times_{(\cY,M_{\cY})} \allowbreak (Y,M_Y) = (X,M_X)$ holds, 
$\Lambda$ is an equivalence of categories. 
\end{prop}

The proposition follows from the following lemma 
(cf. \cite[4.3.3]{shiho1}): 

\begin{lem}\label{crys-hpdilem}
With the above notation, 
let $(\cU,M_{\cU}) \hra (\cP,M_{\cP})$ be a strict open immersion 
$(\cU$ may be empty$)$ and put 
$(U,M_U) := (X,M_X) \times_{(\cP,M_{\cP})} (\cU,M_{\cU})$. 
Let $(E,\epsilon)$ be an object in $\HPDI((X \hra \cP/\cY)^{\log})$, let 
$(L',\epsilon')$ be an object in $\HPD((U \hra \cU/\cY)^{\log}))$ 
such that $L' \in \M(\cU)=\Coh(\cO_{\cU})$ has no $p$-torsion, and 
assume given an isomorphism 
$\alpha: (E,\epsilon)|_{\cU} \os{\cong}{\lra} (\Q \otimes L',\id \otimes 
\epsilon')$ in $\HPDI((U \hra \cU/\cY)^{\log})$. Then there exists 
an object $(L'',\epsilon'')$ in $\HPD((X \hra \cP/\cY)^{\log})$ with the 
isomorphism $(\Q \otimes L'',\id \otimes \epsilon'') \cong (E,\epsilon)$ 
such that $L'' \in \M(\cP) = \Coh(\cO_{\cP})$ has no $p$-torsion and that 
$\alpha$ is induced by an isomorphism 
$(L'',\epsilon'')|_{\cU} \os{\cong}{\lra} (L',\epsilon)$ in 
$\HPD((U \hra \cU/\cY)^{\log})$. 
\end{lem} 

\begin{pf} 
The proof is similar to that of \cite[4.3.3]{shiho1}. 
First note the equivalence of categories 
$\M(\cZ) = \Coh(\cO_{\cZ}), \M(\cZ)_{\Q} = \Coh(\Q \otimes \cO_{\cZ})$ 
for $\cZ = \cP, \cU$. In this proof, we regard $L',E$ as an object
of $\Coh(\cO_{\cU}), \Coh(\Q \otimes \cO_{\cP})$, respectively. 
Take a coherent sheaf $L$ without $p$-torsion 
on $\cP$ such that $\Q \otimes L = E$ holds. 
By the argument in \cite[p.602 lines 16--23]{shiho1}, 
we may assume that $L |_{\cU} \cong L'$ 
holds via $\alpha$. \par 
Assume for the moment that $\cP$ is affine and put $\cP = \Spf A$. 
If we define $D(1)$ as in Definition \ref{defhpdi}, then $D(1)$ is affine 
over $\cP$ and so it is also affine. Put $D(1):=\Spf B(1)$. 
Then it is known that $B(1)$ is flat over $A$ (\cite{kkato}). 
Denote $\Gamma(\cP,E), \Gamma(\cP,L)$ by $E_A,L_A$. 
Then $\epsilon$ naturally induces the isomorphism 
$$ \epsilon_A: B(1) \otimes_A E_A \os{\cong}{\lra} E_A \otimes_A B(1). $$
(Note that, since $B(1)$ is flat over 
$A$ and $E_A$ is isocoherent, $B(1) \otimes_A E_A$ is isomorphic to 
$B(1) \hat{\otimes}_A E_A$ and so on. See for example 
\cite[7.1.6]{gr}.) Let 
$\theta_A: E_A \lra E_A \otimes_{A} B(1)$ 
be the map defined by $\theta_A(x):= \epsilon (1 \otimes x)$ and 
let us define $L''_A$ by $L''_A:=\theta^{-1}(L_A \otimes_A B(1))$. 
Then $L''_A$ is an $A$-module and by the argument in 
\cite[p.602 line 31 --p.604 line 16]{shiho1}, 
we have the following: 
\begin{enumerate}
\item 
$L''_A$ is a coherent sub $A$-module of $L_A$ satisfying 
$\Q \otimes L''_A = E_A$. 
\item 
$\epsilon_A(B(1) \otimes_A L''_A) \subseteq L''_A \otimes_A B(1)$ holds. 
\end{enumerate}
Now we prove that $\epsilon$ induces the isomorphism 
$\epsilon''_A: 
B(1) \otimes_A L''_A \os{\cong}{\lra} L''_A \otimes_A B(1)$. 
(It is implicitely used in \cite{shiho1} and \cite{ogus2}, but the proof 
was omitted there.) To see this, it suffices to prove that 
$C:=\Coker(\epsilon''_A)$
is zero. Since we have 
$B(1) \otimes_A L''_A = B(1) \hat{\otimes}_A L''_A$ and 
$L''_A \otimes_A B(1) = L''_A \hat{\otimes}_A B(1)$, 
it suffices to prove 
$C/p^nC=0$ for any $n$. Note that $A \otimes_{\Delta^*,B(1)} \epsilon''_A$ 
is, by definition, the identity map $L''_A \lra L''_A$. So we have 
$A \otimes_{\Delta^*,B(1)} C = 0$, hence 
$(A/p^nA) \otimes_{\Delta^*, B(1)/p^nB(1)} (C/p^nC)=0$. Since 
$\Ker(\Delta^*: B(1)/p^nB(1) \lra A/p^nA)$ is a nil-ideal and $C/p^nC$ is 
finitely generated, it implies 
$C/p^nC = 0$. So $\epsilon''_A$ is an isomorphism. The pair 
$(L''_A,\epsilon''_A)$ naturally induces an object in 
$\HPD((X \hra \cP/\cY)^{\log})$ (denoted also by $(L''_A,\epsilon''_A)$) 
satisfying $\Q \otimes (L''_A, \epsilon''_A) = (E,\epsilon)$. \par 
Now let us remove the condition that $\cP$ is affine. 
For an affine open sub formal scheme 
$\Spf A$ of $\cP$, we can 
define $(L''_A, \epsilon''_A)$ by the method of the previous paragraph. 
Now let $\Spf A' \subseteq \Spf A \subseteq \cP$ be open immersions. 
Then we have 
{\tiny{
{\allowdisplaybreaks{
\begin{align*} 
L''_{A'} & = \theta_{A'}^{-1}(L_{A'} \otimes_{A'} (B(1) \hat{\otimes}_A A')) \\
& = 
\Ker(E_{A'} \os{\theta_{A'}}{\lra} 
E_{A'} \otimes_{A'} (B(1) \hat{\otimes}_A A') \lra 
(E_{A'}/L_{A'}) \otimes_{A'} (B(1) \hat{\otimes}_A A')) \\ 
& = 
\Ker(E_{A'} \os{\theta_{A'}}{\lra} 
E_{A'} \otimes_{A'} (B(1) \hat{\otimes}_A A') \lra 
(E_{A'}/L_{A'}) \otimes_{A'} (B(1) \otimes_A A'))  \,\,\, 
\text{($E_{A'}/L_{A'}$ is $p$-torsion)} \\ 
& = 
\Ker(E_A \otimes_A A' \os{\theta_A \otimes_A A'}{\lra} 
E_{A} \otimes_{A} B(1) \otimes_A A' \lra 
(E_{A}/L_{A}) \otimes_{A} B(1) \otimes_A A') \\ 
& = 
\Ker(E_A \os{\theta_A}{\lra} E_A \otimes_A B(1) \lra 
(E_A/L_A) \otimes_A B(1)) \otimes_A A' \\ 
& = 
\theta_A^{-1}(L_A \otimes_A B(1)) \otimes_A A' 
= L''_A \otimes_A A'. 
\end{align*}}}}}
So the construction of $(L''_A,\epsilon''_A)$ in the previous paragraph is 
compatible with respect to the open immersion 
$\Spf A' \hra \Spf A$. Moreover, if we have 
$\Spf A \subseteq \cU$, we have 
$$ 
L_A \subseteq \theta^{-1}(L_A \otimes_A B(1)) = L''_A, 
$$ 
since $L |_{\cU} \cong L'$ is stable under $\epsilon$. So we have 
$L_A = L''_A$ in this case. Therefore, 
we can glue $(L''_A, \epsilon''_A)$'s to define 
globally an object $(L'',\epsilon'')$ in $\HPD((X \hra \cP/\cY)^{\log})$ 
which satisfies $\Q \otimes (L'',\epsilon'') = (E, \epsilon)$ and 
$(L'',\epsilon'')|_{\cU} \cong (L',\epsilon')$ via $\alpha$. 
So we are done. 
\end{pf} 

By the same argument as in \cite[0.7.5]{ogus2} and \cite[4.3.4]{shiho1}, 
we have the following (we omit the proof): 

\begin{prop}\label{deszar-crys}
Let $(X,M_X) \os{f}{\lra} (Y,M_Y) \os{\iota}{\hra} (\cY,M_{\cY})$ be as in 
\eqref{diag-p} with $f$ log smooth. Then the descent for finite Zariski 
open coverings holds for the category $I_{\crys}((X/\allowbreak \cY)^{\log})$. 
\end{prop}


\section{Relative log convergent cohomology (I)}

In this section, first we give the definition and some basic 
properties of relative log convergent site. Then we give a proof 
of the relative version of log convergent Poincar\'e lemma. After that, 
we give a proof of the comparison theorem between relative log 
convergent cohomology and relative log crystalline cohomology. 
Some of the results in this section are studied also 
in \cite{nakkshiho}. 
The proofs are basically 
the same as the absolute case (given in \cite{shiho2}), but 
we would like to point out that some of the proofs are 
slightly generalized or slightly simplified than those given in 
\cite{shiho2}. Also, we correct some errors in \cite{shiho2}. \par 
First we give a definition of pre-widenings, widenings and enlargements,
following the absolute case (\cite[2.1.1, 2.1.8, 2.1.9, 
2.1.11]{shiho2}): 

\begin{defn}\label{quaddef}
Let $f:(\cX,M_{\cX}) \lra (\cY,M_{\cY})$ be a morphism in $\pLF$. 
Define the category 
${\cal Q}((\cX/\cY)^{\log})$ 
of quadruples on
$(\cX,M_{\cX})/(\cY, \allowbreak 
M_{\cY})$ as follows$:$ The objects are 
the data 
 $((\cZ,M_{\cZ}),(Z,M_Z),i,z),$ where $(\cZ,M_{\cZ})$ is an object in $\LF$
over $(\cY,M_{\cY})$, 
$(Z,M_Z)$ is an object in $\L$ over $(\cY,M_{\cY})$, 
$i$ is a closed immersion $(Z,M_Z) \hra (\cZ,M_{\cZ})$ over
 $(\cY,M_{\cY})$ and $z$ is a morphism $(Z,M_Z) \lra (\cX,M_{\cX})$ 
 in $\LF$ over $(\cY,M_{\cY})$. 
We define a morphism of quadruples on 
$(\cX,M_{\cX})/(\cY,M_{\cY})$ 
in an obvious way. $($See 
\cite[2.1.8]{shiho2}.$)$ 
\end{defn}

\begin{defn}\label{widedef}
Let $f:(\cX,M_{\cX}) \lra (\cY,M_{\cY})$ be as above. \\
$(1)$ \,
A quadruple $(({\cZ},M_{\cZ}),(Z,M_Z),i,z)$ on $(\cX,M_{\cX})/(\cY,M_{\cY})$ 
is called a pre-widening on $(\cX,M_{\cX})/(\cY,M_{\cY})$ 
if $(\cZ,M_{\cZ})$ is in $\pLF$. \\
$(2)$ \,
A quadruple $((\cZ,M_{\cZ}),(Z,M_Z),i,z)$ 
on $(\cX,M_{\cX})/(\cY,M_{\cY})$ 
is called a widening on $(\cX,M_{\cX})/(\cY,M_{\cY})$ 
if $i$ is a homeomorphic exact closed immersion. $($That is, 
$Z$ is a scheme of definition of $\cZ$ via $i$.$)$ \\
$(3)$ \,
A $($pre$)$-widening $((\cZ,M_{\cZ}),(Z,M_Z),i,z)$ 
on $(\cX,M_{\cX})/(\cY,M_{\cY})$ 
is said to be exact
if $i$ is exact. It is said to be affine if $\cZ$, $z$ and 
the structure morphism $\cZ \lra \cY$ are affine. \\
$(4)$ \, A quadruple $((\cZ,M_{\cZ}),(Z,M_Z),i,z)$ on 
$(\cX,M_{\cX})/(\cY,M_{\cY})$ 
is called an 
enlargement if it is both pre-widening and widening, it is exact and 
$\cZ$ is flat over $\Spf W$. 
\end{defn}

\begin{rem}
In \cite[2.1.11]{shiho2}, we should have added the condition that $T$ is 
flat over $\Spf V$. 
\end{rem}

We often denote a pre-widening, a widening or an enlargement 
$((\cZ,M_{\cZ}),(Z,M_Z), \allowbreak i,z)$ simply by 
$((\cZ,M_{\cZ}), (Z, M_Z))$ or $\cZ$. For a pre-widening 
$\cZ:= ((\cZ,M_{\cZ}),(Z,M_Z), \allowbreak i,z)$, we define
the associated widening by the quadruple 
$((\hat{\cZ},M_{\cZ} \vert_{\hat{\cZ}}), (Z,M_Z),\hat{i},z),$ 
where $\hat{\cZ}$ 
is the completion of
$\cZ$ along $Z$ and $\hat{i}$ is the closed immersion 
$(Z,M_Z) \hra (\hat{\cZ},M_{\cZ} \vert_{\hat{\cZ}})$ induced by $i$. 

We define relative log convergent site and isocrystals on it 
as follows: 

\begin{defn}
Let $\tau$ be one of the words $\{ \Zar(={\mathrm{Zariski}}), 
\et(={\mathrm{etale}})\}$.  
For a morphism $(\cX,M_{\cX}) \lra (\cY,M_{\cY})$ in $\pLF$, 
we define the log convergent site $(\cX/\cY)^{\log}_{\conv,\tau}$ 
of 
$(\cX,M_{\cX})/(\cY,M_{\cY})$ with respect to $\tau$-topology as follows$:$ 
The objects 
 are the enlargements $\cZ$ on $(\cX,M_{\cX})/(\cY,M_{\cY})$ and 
 the morphisms are the morphism of enlargements. 
A family of morphisms 
$$ \{ 
((\cZ_{\lambda}, M_{\cZ_{\lambda}}), 
 (Z_{\lambda}, M_{Z_{\lambda}}), i_{\lambda}, z_{\lambda}) \lra 
 ((\cZ,M_{\cZ}),(Z,M_Z),i,z) \}_{\lambda \in \Lambda} $$
is a covering if the morphisms $(\cZ_{\lambda}, M_{\cZ_{\lambda}})
\ra (\cZ,M_{\cZ})$ are strict, form a covering of $\cZ$ with respect to 
$\tau$-topology and $(Z_{\lambda}, M_{Z_{\lambda}})$ is canonically 
isomorphic to 
$(\cZ_{\lambda}, M_{\cZ_{\lambda}}) \times_{(\cZ,M_{\cZ}),i} 
\allowbreak (Z, \allowbreak M_Z)$. 
When the log structures are 
 trivial, we omit the superscript $^{\log}$ in 
$(\cX/\cY)^{\log}_{\conv,\tau}$. 
We denote 
the right derived functor 
$($resp. the 
$q$-th right derived functor$)$ of the functor 
$$ 
(\cX/\cY)^{\log,\sim}_{\conv,\tau} \lra \cY^{\sim}_{\Zar}; 
\,\, \cE \mapsto (U \mapsto 
\Gamma((\cX \times_{\cY} U/U)^{\log}_{\conv,\tau},\cE)) 
$$ 
by $Rf_{\cX/\cY,\conv *}\cE$ 
$($resp. $R^qf_{\cX/\cY,\conv *}\cE)$. We call 
$R^qf_{\cX/\cY,\conv *}\cE$ 
the $q$-th relative log convergent cohomology of $(\cX,M_{\cX})/
(\cY,M_{\cY})$ with coefficient $\cE$. 
We denote the sheaf on 
$(\cX/\cY)^{\log}_{\conv,\tau}$ defined by 
$\cZ \mapsto \Gamma(\cZ, {\cal O}_{\cZ})$ 
$($resp. $\Q \otimes_{\Z} \Gamma(\cZ, {\cal O}_{\cZ}))$
by $\cO_{\cX/\cY}$ $($resp. $\cK_{\cX/\cY})$. 
\end{defn}

\begin{defn}
Let the notations be as above. An isocrystal on the log 
convergent site $(\cX/\cY)^{\log}_{\conv,\tau}$ is a sheaf 
${\cal E}$ on $(\cX/\cY)^{\log}_{\conv,\tau}$ satisfying 
the following conditions$:$ 
\begin{enumerate}
\item
For any enlargement $\cZ$, the sheaf ${\cal E}_{\cZ}$ on $\cZ$ 
induced by ${\cal E}$ is an isocoherent sheaf. 
\item 
For any morphism $f:\cZ' \lra \cZ$ of enlargements, the 
homomorphism $f^*{\cal E}_{\cZ} \lra {\cal E}_{\cZ'}$ of 
sheaves on $\cZ'$ induced by ${\cal E}$ is an 
isomorphism. 
\end{enumerate}
We denote the category of isocrystals on 
$(\cX/\cY)^{\log}_{\conv,\tau}$ 
by 
$I_{\conv,\tau}((\cX/\cY)^{\log})$. 
When the log structures are trivial, we omit the superscript 
$^{\log}$. 
\end{defn}

\begin{defn} 
Let the notations be as above. Then an isocrystal 
$\cE$ is said to be locally free if, 
for any enlargement $\cZ$, the sheaf ${\cal E}_{\cZ}$ on $\cZ$ 
induced by 
${\cal E}$ is a locally free $\Q \otimes_{\Z} {\cal O}_{\cZ}$-module 
in the sense of Definition \ref{deflf1}. 
\end{defn}

For a morphism $(\cX,M_{\cX}) \lra (\cY,M_{\cY})$ and 
a (pre-)widening $\cZ$ on $(\cX,M_{\cX})/(\cY,M_{\cY})$, we can define 
the notion of localized log convergent site 
$(\cX/\cY)^{\log}_{\conv,\tau} \vert_{\cZ}$ and the category of isocrystals 
$I_{\conv,\tau}((\cX \allowbreak /\cY)^{\log}|_{\cZ})$ on 
it in the same way as 
\cite[2.1.17, 2.1.19]{shiho2}. \par 
Note that, 
for a simplicial object $(\cX^{(\b)},M_{\cX^{(\b)}})$ in $\pLF$ and a 
morphism $(\cX^{(\b)},\allowbreak M_{\cX^{(\b)}}) \allowbreak 
\lra (\cY,M_{\cY})$ in $\pLF$, 
we can define the log convergent topos 
$(\cX^{(\b)}/\cY)^{\log,\sim}_{\conv,\tau}$ of 
$(\cX^{(\b)},M_{\cX^{(\b)}})/(\cY,M_{\cY})$ 
(cf. \cite[p.46--47]{shiho2}). 
We can prove the following proposition in the same way as 
\cite[2.1.20]{shiho2} (we omit the proof): 

\begin{prop}\label{cohdes}
Let $(\cX,M_{\cX}) \lra (\cY,M_{\cY})$ be a morphism in $\pLF$ and 
let $\tau$ be one of the words 
$\{ \Zar(={\mathrm{Zariski}}), \et(={\mathrm{etale}})\}$.
Let $g^{(\b)}:(\cX^{(\bullet)},M_{\cX^{(\b)}}) \lra 
(\cX,M_{\cX})$ be a strict $\tau$-hypercovering. 
Let $\theta:= (\theta_*,\theta^{-1}): 
({\cX}^{(\bullet)}/\cY)^{\log,\sim}_{\conv,\tau} \lra 
(\cX/\cY)^{\log,\sim}_{\conv,\tau} $
be the morphism of topoi characterized by 
$\theta^{-1}(E)^{(i)} := g^{(i),-1}(E).$ 
Then, for any abelian sheaf $E$ on $(\cX/\cY)^{\log}_{\conv,\tau}$, 
the canonical homomorphism 
$E \lra R\theta_* \theta^{-1} E$
is a quasi-isomorphism. 
\end{prop}

\begin{rem}
The description of $\theta^{-1}$ in \cite[2.1.20]{shiho2} was
wrong. The correct one is given above. 
\end{rem}

The equivalence of the category of isocrystals 
on $(\cX/\cY)^{\log}_{\conv,\et}$ and that on 
$(\cX/\cY)^{\log}_{\conv,\Zar}$ holds also in relative case 
(\cite[2.1.21]{shiho2}, we omit the proof): 

\begin{prop}\label{zid}
Let $(\cX,M_{\cX}) \lra (\cY,M_{\cY})$ 
be as above and let $\cZ$ be a 
$($pre-$)$widening of $(\cX,M_{\cX})/(\cY,M_{\cY})$. Let us denote the 
canonical morphism of topoi 
$$ (\cX/\cY)^{\log,\sim}_{\conv,\et} \lra 
(\cX/\cY)^{\log,\sim}_{\conv,\Zar} \quad \text{$($resp.} \,\,\,
(\cX/\cY)^{\log,\sim}_{\conv,\et} \vert_{\cZ} 
\lra (\cX/\cY)^{\log,\sim}_{\conv,\Zar} \vert_{\cZ} \,\,\, \text{$)$} $$ 
by $\epsilon$. Then$:$ 
\begin{enumerate}
\item
for any $\cE \in I_{\conv,\et}((\cX/\cY)^{\log})$ 
$($resp. $\cE \in I_{\conv,\et}((\cX/\cY)^{\log} \vert_{\cZ}))$, we have 
$R\epsilon_* \cE = \epsilon_*\cE$. 
\item 
The functor $\cE \mapsto \epsilon_* \cE$ induces the equivalence 
of categories 
\begin{align*}
I_{\conv,\et}((\cX/\cY)^{\log}) & \overset{\sim}{\lra} 
I_{\conv,\Zar}((\cX/\cY)^{\log}) \\ 
\text{$($resp.} \,\,\, 
I_{\conv,\et}((\cX/\cY)^{\log} \vert_{\cZ}) & \overset{\sim}{\lra} 
I_{\conv,\Zar}((\cX/\cY)^{\log} \vert_{\cZ}) 
\,\,\, \text{$).$} 
\end{align*}
\end{enumerate}
\end{prop} 

In the rest of this paper, we denote the category 
$I_{\conv,\et}((\cX/\cY)^{\log})=
I_{\conv,\Zar}((\cX/\cY \allowbreak )^{\log})$ (resp. 
$I_{\conv,\et}((\cX/\cY)^{\log} \vert_{\cZ})= 
I_{\conv,\Zar}((\cX/\cY)^{\log} \vert_{\cZ})$) simply by 
$I_{\conv}((\cX/\cY)^{\log})$ (resp. 
$I_{\conv}((\cX/\cY)^{\log} \vert_{\cZ})$) and call it an isocrystal 
on $(\cX/\cY)^{\log}_{\conv}$ (resp. $(\cX/\cY)^{\log}_{\conv} |_{\cZ}$), 
by abuse of terminology. \par 
Next we introduce the notion of the exactification of a closed immersion 
of fine log formal $\cB$-schemes. (This is essentially 
introduced in \cite[0.9]{chfo} as 
log formal tube. We will not use this terminology because we would like to 
use the word `tube' for the rigid analytic one.) 

\begin{propdefn}
Let $\cC$ be the category of homeomorphic 
exact closed immersions $i: (Z,M_Z) \hra 
(\cZ,M_{\cZ})$ in $\LF$ and let 
$\cD$ be the category of closed immersions in $\LF$. Then the 
canonical inclusion functor $\cC \lra \cD$ has a right adjoint of 
the form $((Z,M_Z) \hra (\cZ,M_{\cZ})) \mapsto ((Z,M_Z) \hra (\cZ^{\ex},
M_{\cZ^{\ex}}))$. We call this functor the exactification. 
\end{propdefn}

\begin{pf} 
Let $i: (Z,M_Z) \hra (\cZ,M_{\cZ})$ be a closed immersion in $\LF$ and 
we prove the existence of $(Z,M_Z) \hra (\cZ^{\ex},
M_{\cZ^{\ex}})$. To prove it, 
we may replace $\cZ$ by the formal completion 
of $\cZ$ along $Z$. Moreover, since it suffices to prove the claim 
etale locally 
(by etale descent), we may assume that $i$ admits a chart 
$(P_{\cZ} \ra M_{\cZ}, Q_Z \ra M_Z, P \os{\alpha}{\ra} Q)$ such that 
$\alpha^{\gp}$ is surjective. If we put 
$P':= \alpha^{\gp,-1}(Q), \cZ' := \cZ \hat{\times}_{\Spf \Z_p\{P\}} 
\Spf \Z_p\{P'\}, M_{\cZ'} :=$ log structure associated to $P'_{\cZ'} 
\ra {\cal O}_{\cZ'}$, we obtain a factorization 
$$ (Z,M_Z) \os{i'}{\lra} (\cZ',M_{\cZ'}) \os{f}{\lra} (\cZ,M_{\cZ}) $$ 
such that $i'$ is an exact closed immersion and $f$ is formally log etale. 
(Here we define the formal log etaleness of a morphism in $\LF$ 
by infinitesimal lifting property.) If we define 
$\cZ^{\ex}$ to be the completion of $\cZ'$ along $Z$ and $M_{\cZ^{\ex}}$ 
to be the pull-back of $M_{\cZ'}$ to $\cZ^{\ex}$, we have the diagram 
$$ (Z,M_Z) \os{i^{\ex}}{\lra} (\cZ^{\ex},M_{\cZ^{\ex}}) 
\os{g}{\lra} (\cZ,M_{\cZ}) $$ 
such that $i^{\ex}$ is a homeomorphic exact closed immersion and 
that $g$ is formally log etale. We can check that 
$i \mapsto i^{\ex}$ is indeed the desired functor by this diagram. 
\end{pf}

\begin{cor}\label{affine}
Let $i: (Z,M_Z) \hra (\cZ,M_{\cZ})$ be a closed immersion in $\LF$ and 
let $(Z,M_Z) \hra (\cZ^{\ex},M_{\cZ^{\ex}})$ be its exactification. Then 
the canonical morphism $(\cZ^{\ex},M_{\cZ^{\ex}}) \allowbreak 
\lra (\cZ,M_{\cZ})$ defined 
by adjointness is formally log etale and affine. 
\end{cor}

For a (pre-)widening $\cZ=((\cZ,M_{\cZ}),(Z,M_Z),i,z)$, the quadruple 
$\cZ^{ex} := ((\cZ^{\ex}, \allowbreak M_{\cZ^{\ex}}), \allowbreak 
(Z,M_Z),i^{\ex},z)$ is a widening. 
It is called the exactification of $\cZ$. Note that 
the widening obtained as 
the exactification of a 
pre-widening is equal to the exactification of the widening associated to 
the  given pre-widening. 
The notion of the exactification was defined only 
under the existence of a chart 
in \cite{shiho2} and this forced us to do complicated arguments at several 
parts of the paper \cite{shiho2}. With this general definition of the 
exactification, we can 
slightly generalize some notions and 
simplify some arguments in \cite{shiho2}. 
For example, we can define the notion of the system of universal enlargements 
of (not necessarily exact) (pre-)widenings as follows: 

\begin{defn}\label{defsysenl}
Let $\cZ := ((\cZ,M_{\cZ}),(Z,M_Z),i,z)$ be a $($pre-$)$widening. 
Let $((\cZ^{\ex}, \allowbreak M_{\cZ^{\ex}}), \allowbreak 
(Z,M_Z),i^{\ex},z)$ be 
the exactification of $\cZ$ and let ${\cal I}$ be the 
ideal $\Ker ({\cal O}_{\cZ^{\ex}} \lra {\cal O}_{Z})$. For 
$n \in \N$, let 
$B_n(\cZ)$ be the 
formal blow-up of $\cZ^{\ex}$ with respect 
to the ideal $p {\cal O}_{\cZ^{\ex}} + {\cal I}^n$, 
let 
$T'_n(\cZ)$ be the open sub formal scheme
$$ \{x \in B_{n}(\cZ) \,\vert\, (p {\cal O}_{\cZ^{\ex}} + 
{\cal I}^n) \cdot {\cal O}_{B_{n}(\cZ),x} = 
p{\cal O}_{B_{n}(\cZ),x} \} $$ 
and 
let $T_n(\cZ)$ be 
the closed subscheme of $T'_n(\cZ)$ defined by the ideal 
$\{x \in {\cal O}_{T'_n(\cZ)}\,\vert\, \exists n, \allowbreak p^nx=0\}$. 
Let $\lambda_n:T_{n}(\cZ) \lra \cZ$ be the canonical 
morphism and let $Z_n := \lambda_n^{-1}(Z)$. 
Then the quadruple 
{\footnotesize{ $$ T_{n}(\cZ) := ((T_{n}(\cZ),M_{\cZ^{\ex}} 
|_{T_{n}(\cZ)}),(Z_n, M_Z|_{Z_n}), 
Z_n \hra T_{n}(\cZ), (Z_n,M_Z) \overset{\lambda_n}{\lra} 
(Z,M_Z) \overset{z}{\lra} (X,M))  $$ }}
is an enlargement for each $n$ and the family 
$\{T_{n}(\cZ)\}_{n \in \N}$ 
forms an inductive system of enlargements. The morphisms 
$\lambda_n$'s define the morphisms of $($pre-$)$widenings 
$T_{n}(\cZ) \lra \cZ \,(n \in \N)$ which is compatible with 
transition morphisms. We call this inductive system 
the system of universal enlargements of $\cZ$. 
\end{defn} 

We can see that the analogues of 
\cite[2.1.23--27]{shiho2} are true also in relative situation. 
That is, we have the following (we omit the proof): 
 
\begin{lem}\label{2.1.23}
Let $\cZ$ be a widening and let $\cZ'$ be an enlargement. Then 
a morphism $\cZ' \lra \cZ$ as quadruples factors through $T_n(\cZ)$ for some 
$n$. Moreover, such a factorization is unique 
as a morphism to the inductive system 
$\{T_{n}(\cZ)\}_{n \in \N}$. 
\end{lem} 

\begin{lem}\label{2.1.24}
With the above notation, the morphism of sheaves $h_{T_n(\cZ)} \lra h_{\cZ}$ 
is injective. $($Here $h_?$ denotes the sheaf on relative log convergent 
site associated to $?.)$ 
\end{lem}

\begin{lem}\label{2.1.25}
Let $\cZ$ be a pre-widening and let $\hat{\cZ}$ be the associated widening. 
Then we have the canonical isomorphism of enlargements 
$T_n(\hat{\cZ}) \os{=}{\lra} T_n(\cZ)$. 
\end{lem}

\begin{lem}\label{2.1.26}
Let 
$$g:((\cZ',M_{\cZ'}),(Z',M_{Z'}),i',z') \lra ((\cZ,M_{\cZ}),(Z,M_Z),i,z) $$
be a morphism of $($pre-$)$widenings and assume that 
$(\cZ',M_{\cZ'}) \times_{(\cZ,M_{\cZ})} (Z,M_Z) = (Z',M_{Z'})$ 
holds naturally and that $\cZ' \lra \cZ$ is flat.  
Then $g$ induces the natural isomorphism of enlargements 
$$ T_{n}(\cZ') \overset{\sim}{\lra} T_{n}(\cZ) \times_{\cZ} {\cZ'}. $$
\end{lem}

\begin{lem}\label{zure}
Let 
$$ g: \cZ := ((\cZ,M_{\cZ}),(Z,M_Z),i,z) \lra 
\cZ':= ((\cZ,M_{\cZ}),(Z',M_{Z'}),i',z') $$ 
be a morphism of $($pre-$)$widenings such that $g$ is 
identity on $\cZ$ and $(Z,M_Z) \lra (Z',M_{Z'})$ is an exact closed 
immersion. Assume that there exists 
an ideal $J$ of ${\cal O}_{\cZ}$ and $m \in \N$ such that 
$J^{m+1} \subseteq \pi {\cal O}_{\cZ}$ and 
$I_Z \subseteq I_{Z'} + J$ holds, where $I_Z, I_{Z'}$ are 
the defining ideals of $Z$, $Z'$ in $\cZ$, respectively. 
Denote the morphism of enlargements 
$T_{n}(\cZ) \lra T_{n}(\cZ')$ induced by $g$ by $g_n$. 
Then$:$ 
\begin{enumerate}
\item 
There exists a homomorphism of formal schemes 
$h_n: T_{n}(\cZ') \lra T_{m+n}(\cZ)$ such that the composites 
$$ T_{n}(\cZ') \overset{h_n}{\lra} T_{m+n}(\cZ) \overset{g_{m+n}}{\lra} 
T_{m+n}(\cZ'), \,\,\,\, 
T_{n}(\cZ) \overset{g_n}{\lra} T_{n}(\cZ') \overset{h_n}{\lra} 
T_{m+n}(\cZ) $$ 
coincide with the canonical transition morphisms. 
\item 
For an isocrystal $\cE$, we have the natural isomorphism 
$\varphi_n: 
h_n^*\cE_{T_{m+n}(\cZ)} \overset{\sim}{\lra} \cE_{T_{n}(\cZ')}$ 
such that the composites 
$$ 
h_n^* \circ g_{m+n}^* \cE_{T_{m+n}(\cZ')} 
\overset{h_n^*\cE(g_{m+n})}{\lra} 
h_n^* \cE_{T_{m+n}(\cZ)} 
\overset{\varphi_n}{\lra} 
\cE_{T_{n}(\cZ')}, 
$$
$$
g_n^* \circ h_{n}^* \cE_{T_{m+n}(\cZ)} 
\overset{g_n^*\varphi_n}{\lra} 
g_n^* \cE_{T_{n}(\cZ')} 
\overset{E(g_n)}{\lra} 
\cE_{T_{n}(\cZ)}, 
$$
coincide with the isomorphisms induced by the canonical 
morphism of enlargements $T_{n}(\cZ') \lra T_{m+n}(\cZ')$ and 
$T_{n}(\cZ) \lra T_{m+n}(\cZ)$, respectively. 
\end{enumerate}
\end{lem}

\begin{rem} 
In \cite[2.1.25, 2.1.26]{shiho2}, we forgot to assume that $T$ 
is exact. However, now we do not need the exactness hypothesis 
(as written in Lemmas \ref{2.1.25}, \ref{2.1.26}), 
because now we can define the system of universal enlargements for 
(pre-)widenings which are not necessarily exact. 
\end{rem} 

As another application of the exactification of a closed immersion 
of fine log formal $\cB$-schemes, we define the notion of tubular 
neighborhood of fine log formal $\cB$-schemes as follows 
(In \cite[2.2.4]{shiho2}, 
it was defined under certain assumption, but we can define it without 
any assumption here): 

\begin{defn}
For a closed immersion $(Z,M_Z) \hra (\cZ,M_{\cZ})$ of 
fine log formal $\cB$-schemes, we define the tubular neighborhood 
$]Z[^{\log}_{\cZ}$ of $(Z,M_Z)$ in $(\cZ,M_{\cZ})$ as the rigid 
analytic space $\cZ^{\ex}_K$. We define the specialization map 
$$ 
\sp:]Z[^{\log}_{\cZ} \lra \hat{\cZ} (:=\text{$($the completion of 
$\cZ$ along $Z)$} \simeq Z) 
$$
by the composite $]Z[^{\log}_{\cZ} = \cZ^{\ex}_K \os{\sp}{\lra} \cZ^{\ex} \lra 
\hat{\cZ}$, where the first map is the usual specialization map and the 
second map is the one induced by the canonical morphism $\cZ^{\ex} \lra 
\cZ$. 
\end{defn}

\begin{rem}
$]Z[^{\log}_{\cZ}$ and $\{T_n(\cZ)\}_n$ are related by 
$]Z[^{\log}_{\cZ} = \bigcup_{n=1}^{\infty} T_n(\cZ)_K$. 
\end{rem}

Next, we recall the relation between the category of isocrystals and 
certain categories of stratifications and 
give the definition of the log de Rham complex on tubular neighborhood 
associated to an isocrystal. 
Let us consider the situation
\begin{equation}
\begin{CD}
(X,M_{X}) @>i>> (\cP,M_{\cP}) \\
@VfVV @VgVV \\
(\cY,M_{\cY}) @= (\cY,M_{\cY}), 
\end{CD}
\end{equation}
where $(X,M_X)$ is an object in $\L$, 
$f$ is a morphism in $\pLF$, $i$ is a closed immersion in 
$\pLF$ and $g$ is a formally log smooth morphism 
in $\pLF$. 
For $n \in \N$, 
let $(\cP(n),M_{\cP(n)})$ be the $(n+1)$-fold 
fiber product of $(\cP,M_{\cP})$ over $(\cY,M_{\cY})$. Then, we have a 
closed immersion $i(n): (X,M_{X}) \hra (\cP(n),M_{\cP(n)})$ and 
the quadruple $\cP(n) := ((\cP(n),M_{\cP(n)}),(X,M_X),i(n),\id)$ is a 
pre-widening of $(X,M_{X})/(\cY,M_{\cY})$. So we have a system of 
universal enlargements $\{T_m(\cP(n))\}_m$. Moreover, we have 
the `projections' 
{\small{ $$ p_{i,m}: T_{m}(\cP(1)) \lra T_{m}(\cP) 
\,\,\,\, (i=1,2), \quad 
p_{ij,m}: T_{m}(\cP(2)) \lra T_{X,m}(\cP(1)) 
\,\,\,\, (1 \leq i < j \leq 3) $$ }}
and the `diagonal morphism' $\Delta_m: T_{m}(\cP) \lra T_{m}(\cP(1)),$ 
which are morphisms of enlargements and compatible with respect to $m$. 
Let $\Strat'((X \hra \cP/\cY)^{\log})$ 
be the category 
of compatible family of isocoherent sheaves $E_m$ on 
$T_{m}(\cP) \, (m \in \N)$ endowed with compatible isomorphisms 
$\epsilon_m: p^*_{2,m}E \overset{\sim}{\lra} 
p^*_{1,m}E$ satisfying 
$\Delta_m^*(\epsilon_m) = \id$, 
$p_{12,m}^*(\epsilon_m) \circ 
p_{23,m}^*(\epsilon_m) = p_{13,m}^*(\epsilon_m).$ \par 
On the other hand, from the closed immersion $i(n)$, we can form the 
tubular neighbouhood $]X[^{\log}_{\cP(n)}$ and 
we have the projections 
$$ p_i: \, ]X[^{\log}_{\cP(1)} \lra ]X[^{\log}_{\cP} \,\,\,\, 
(i=1,2), \quad   
p_{ij}: \, ]X[^{\log}_{\cP(2)} \lra ]X[^{\log}_{\cP(1)} \,\,\,\, 
(1 \leq i < j \leq 3) $$
and the diagonal morphism $\Delta: 
\, ]X[^{\log}_{\cP} \lra ]X[^{\log}_{\cP(1)}.$ 
Let $\Strat''((X \hra \cP/\cY)^{\log})$ be the 
category of pairs $(E,\epsilon)$, where $E$ is a coherent 
${\cal O}_{]X[^{\log}_{\cP}}$-module and $\epsilon$ is an 
${\cal O}_{]X[^{\log}_{\cP(1)}}$-linear isomorphism 
$p_2^*E \overset{\sim}{\lra} p_1^*E$ satisfying 
$\Delta^*(\epsilon) = \id$, 
$p_{12}^*(\epsilon) \circ p_{23}^*(\epsilon) = 
p_{13}^*(\epsilon)$. Then we have the following proposition as in the 
absolute case (\cite[2.2.7]{shiho2}): 

\begin{prop}\label{strat}
With the above notation, we have the functorial equivalence of categories 
\begin{equation*}
I_{\conv}((X/\cY)^{\log}) \simeq 
\Strat'((X \hra \cP/\cY)^{\log}) \simeq 
\Strat''((X \hra \cP/\cY)^{\log}). 
\end{equation*}
\end{prop}

Let $\cE$ be an isocrystal on $(X,M_{X})/(\cY,M_{\cY})$ and let 
$(E,\epsilon)$ be the associated object in 
$\Strat''((X \hra \cP/\cY)^{\log})$. 
If we denote the first log infinitesimal neighborhood of 
$(\cP,M_{\cP})$ into $(\cP(1),M_{\cP(1)})$ by $(\cP^1, M_{\cP^1})$, 
$\epsilon$ induces the $\cO_{]X[^{\log}_{\cP^1}}$-linear isomorphism 
$$ \epsilon_1: 
\cO_{]X[^{\log}_{\cP^1}} \otimes_{\cO_{]X[^{\log}_{\cP}}} E \lra E 
\otimes_{\cO_{]X[^{\log}_{\cP}}} \cO_{]X[^{\log}_{\cP^1}} $$ 
and it induces the log connection 
$$ \nabla: E \lra E \otimes_{{\cal O}_{]X[^{\log}_{\cP}}} 
\omega^1_{]X[^{\log}_{\cP}/\cY_K} $$
by $\nabla(e) := \epsilon_1(1 \otimes e) - e \otimes 1$. 
(Here $\omega^1_{]X[^{\log}_{\cP}/\cY_K} := 
\omega^1_{\cP_K/\cY_K}|_{]X[^{\log}_{\cP}}$ and 
$\omega^1_{\cP_K/\cY_K}$ is the coherent $\cO_{\cP_K}$-module associated to 
$\Q \otimes_{\Z} \omega^1_{\cP/\cY} \in \Coh(\Q \otimes \cO_{\cP})$.) 
Then we have the following lemma. 

\begin{lem}
The log connection $\nabla$ above is integrable. 
\end{lem}

Note that the proof of this lemma is harder than the analogous 
one in the previous paper (\cite[2.2.8]{shiho2}) due to the lack of 
the `of Zarisiki type' hypothesis. 

\begin{pf} 
When there exists a chart of the log formal scheme $(\cP,M_{\cP})$, 
we can prove the lemma in the same way as \cite[1.2.7, 2.2.8]{shiho2}. 
So it suffices to prove that we may work etale locally on $\cP$. \par 
Let $\{(E_m,\epsilon_m)\}_m$ be the object in 
$\Strat'((X \hra \cP/\cY)^{\log})$ associated to 
$\cE$, let $({T_m(\cP)^n}', \allowbreak M\allowbreak {}_{{T_m(\cP)^n}'})$ be the 
$n$-th log infinitesimal neighborhood 
of $(T_m(\cP),M_{T_m(\cP)})$ 
in $(T_m(\cP), \allowbreak M_{T_m(\cP)}) \allowbreak 
\times_{(\cY,M_{\cY})} (T_m(\cP),M_{T_m(\cP)})$ and 
let $({T_m(\cP)^n},M_{{T_m(\cP)^n}}) \hra 
({T_m(\cP)^n}',M_{{T_m(\cP)^n}'})$ be the exact closed immersion 
defined by the ideal $\{x \in \cO_{T_m(\cP)^n}\,\vert\, \exists n, 
p^nx=0\}$. Then the canonical 
morphism 
$$(T_m(\cP)^n,M_{T_m(\cP)^n}) \lra 
(T_m(\cP),M_{T_m(\cP)}) \times_{(\cY,M_{\cY})} 
(T_m(\cP),M_{T_m(\cP)}) \lra (\cP(1),M_{\cP(1)})$$ 
factors through $(T_{m+l}(\cP(1)),M_{T_{m+l}(\cP(1))})$ for some 
$l \in \N$ 
by the universality of exactification and blow-up. So, by 
pulling-back $\epsilon_{m+l}$ to $T_m(\cP)^n$, we obtain the isomorphism 
$$ \epsilon'_{m,n}: \cO_{T_m(\cP)^n} \otimes_{\cO_{T_m(\cP)}} E_m 
\overset{\sim}{\lra} E_m \otimes_{\cO_{T_m(\cP)}} \cO_{T_m(\cP)^n} $$
and $\epsilon'_{m,1}$ induces the connection 
$$ \nabla_m: E_m \lra E_m \otimes_{\cO_{T_m(\cP)}} 
\omega^1_{T_m(\cP)/\cY}. $$
Then it is easy to see that the compatible family 
$\{(E_m,\nabla_m)\}$ induces $(E,\nabla)$ via the equivalence 
\begin{equation*}
\left( 
\begin{aligned}
& \text{compatible family of} \\
& \text{isocoherent sheaves on} \\
& \{ T_{m}(\cP) \}_m 
\end{aligned}
\right)
\overset{\sim}{\lra} 
\left( 
\begin{aligned}
& \text{coherent} \\
& \text{${\cal O}_{]X[^{\log}_{\cP}}$-module} \\
\end{aligned}
\right).
\end{equation*}
So it suffices to prove the integrability of $(E_m,\nabla_m)$ to prove 
the lemma, and we may work etale locally on $\cP$ to check it. 
So we are done. 
\end{pf}

Thanks to the above lemma, 
we can define the log de Rham complex 
{\small{
$$ \DR(]X[^{\log}_{\cP}/\cY_K,{\cal E}) := [0 @>>> E @>{\nabla}>> 
E \otimes_{{\cal O}_{]X[^{\log}_{\cP}}} \omega^1_{]X[^{\log}_{\cP}/\cY_K} 
@>{\nabla}>> \cdots @>{\nabla}>> 
E \otimes_{{\cal O}_{]X[^{\log}_{\cP}}} \omega^q_{]X[^{\log}_{\cP}/\cY_K} 
@>{\nabla}>> \cdots] $$ }}
on $]X[^{\log}_{\cP}/\cY_K$ 
associated to the isocrystal ${\cal E}$ 
in standard way. 

\begin{rem}\label{generalized-dr}
We can give the analogue of the categories $\Strat'((X \hra \cP/\cY)^{\log}), 
\Strat''((X \allowbreak \hra \cP/\cY)^{\log})$ and the log de Rham complex 
$\DR(]X[^{\log}_{\cP}/\cY_K,\cE)$ in a slightly more generalized 
situation. Assume given a commutative diagram 
\begin{equation*}
\begin{CD}
(X,M_X) @>i>> (\fP,M_{\fP}) \\ 
@VfVV @VgVV \\
(\cY,M_{\cY}) @= (\cY,M_{\cY}), 
\end{CD}
\end{equation*}
where $(X,M_X)$ is an object in $\L$, $f$ is a morphism in 
$\pLF$, $i$ is a closed immersion in $\LF$ and $g$ is a morphism 
in $\LF$ satisfying the following condition $(*)$: \\
\quad \\
$(*)$ \,\,\, Zariski locally on $\fP$, there exists a diagram 
\begin{equation*}
\begin{CD}
(X,M_X) @>i>> (\fP,M_{\fP}) @>{i'}>> (\ti{\fP},M_{\ti{\fP}})\\ 
@VfVV @VgVV @V{g'}VV \\
(\cY,M_{\cY}) @= (\cY,M_{\cY}) @= (\cY,M_{\cY}), 
\end{CD}
\end{equation*}
where $g'$ is a formally log smooth morphism in $\pLF$, 
$i'$ is a morphism in $\LF$ such that $i' \circ i$ is again a closed 
immersion and that $i'$ induces the isomorphism 
$(\fP^{\ex},M_{\fP^{\ex}}) \os{=}{\lra} (\ti{\fP}^{\ex},M_{\ti{\fP}^{\ex}})$, 
where $^{\ex}$ denotes the exactification of the closed immersion from 
$(X,M_X)$. \\
\quad \\
For $n \in \N$, let $(\fP(n),M_{\fP(n)})$ be the 
$(n+1)$-fold fiber product of $(\fP,M_{\fP})$ over $(\cY,M_{\cY})$ and let 
$i(n):(X,M_X) \hra (\fP(n),M_{\fP(n)})$ be the closed immersion induced by 
$i$. Then, if we denote the exactification of $i(n)$ by 
$(\fP(n)^{\ex},M_{\fP(n)^{\ex}})$, 
$((\fP(n)^{\ex},M_{\fP(n)^{\ex}}),(X,M_X))$ forms a widening of 
$(X,M_X)/(\cY,M_{\cY})$. \par 
Using $T_m(\fP(n)^{\ex})$ instead of $T_m(\cP(n))$, we can define the 
category $\Strat'((X \hra \fP/\cY)^{\log})$, and using 
$]X[^{\log}_{\fP(n)} = \fP(n)^{\ex}_K$ 
instead of $]X[^{\log}_{\cP(n)}$, we can define the category 
$\Strat''((X \hra \fP/\cY)^{\log})$. Then the analogue of 
Proposition \ref{strat} is true also in this situation. 
(We may work Zariski locally on $\fP$ and then we can reduce to Proposition 
\ref{strat} by the condition $(*)$.) Moreover, for an isocrystal $\cE$ on 
$(X/\cY)^{\log}_{\conv}$, we can define the log de Rham complex 
$\DR(]X[^{\log}_{\fP}/\cY_K,{\cal E})$ on 
$]X[^{\log}_{\fP} = \fP^{\ex}_K$ in the same way as above. 
\end{rem} 

Next we give a sketch of 
the proof of log convergent Poincar\'{e} lemma in relative 
situation. Let $(X,M_X)$ be an object in $\L$, 
$(X,M_{X}) \lra (\cY,M_{\cY})$ be a morphism in 
$\pLF$ and let $\cZ:=((\cZ,M_{\cZ}),(Z,M_Z),i,z)$ 
be a widening on $(X,M_{X})/(\cY,M_{\cY})$. 
Then we have morphisms of topoi 
$$ u: (X/\cY)^{\log,\sim}_{\conv,\Zar} \lra X^{\sim}_{\Zar}, \quad 
j_{\cZ}: (X/\cY)^{\log,\sim}_{\conv,\Zar} |_{\cZ} \lra 
(X/\cY)^{\log,\sim}_{\conv,\Zar} $$
and a morphism of ringed topoi 
$$ \phi_{\cZ}: 
(X/\cY)^{\log,\sim}_{\conv,\Zar} |_{\cZ} \lra \cZ^{\sim}_{\Zar}, $$
as in the absolute case \cite[pp.90--91]{shiho2}. On the other hand, 
we can define the direct limit site $\vec{\cZ}$ as follows, 
as in \cite[2.1.28]{shiho2} (note that we do not need 
the exactness of $\cZ$, because now we have the notion of system 
of universal enlargements in general): 

\begin{defn}
Let $\cZ$ be as above. Then we define the direct limit
 site $\vec{\cZ}$ as follows$:$ Objects are the open sets of some $T_n(\cZ)$. 
 For open sets $U \subseteq T_n(\cZ)$ and $V \subseteq T_m(\cZ)$,
 $\Hom_{\vec{\cZ}}(U,V)$ is empty unless $n \leq m$ and in the case 
 $n \leq m$, 
 $\Hom_{\vec{\cZ}}(U,V)$ 
 is defined to be the set of morphisms $f:U \lra V$ which
 commutes with the transition morphism $T_n(\cZ) \lra T_m(\cZ)$. 
 The coverings are defined by Zarisi topology for each object.  \par
 We define the structure sheaf ${\cal O}_{\vec{\cZ}}$ by
 ${\cal O}_{\vec{\cZ}}(U) := \Gamma(U,{\cal O}_U)$. 
A sheaf of $\Q \otimes_{\Z} {\cal O}_{\vec{\cZ}}$-modules $E$ is called
 crystalline, if, for any transition map $\psi:T_n(\cZ) \lra T_m(\cZ)$, 
the induced map of sheaves 
  $$ {\cal O}_{T_n(\cZ)} \otimes_{\psi^{-1}{\cal O}_{T_m(\cZ)}} \psi^{-1}E_m
   {\lra} E_n $$ 
is an isomorphism, 
 where $E_n$, $E_m$ denote the sheaves on $T_{n}(\cZ)$, $T_{m}(\cZ)$ 
 induced by $E$.
\end{defn}

We define the morphism of topoi $\gamma:
\vec{\cZ}^{\sim} \lra \cZ^{\sim}_{\Zar}$ in order that 
$\gamma^*$ is the pull-back functor and that $\gamma_*$ is the functor of 
taking inverse limit of the direct images. Then we have the following 
lemma, as \cite[2.1.31]{shiho2} and \cite[0.3.7]{ogus2}. 
(The proof is the same. Note that 
we use the affinity of $T_n(\cZ) \lra \cZ$ in the proof, 
which is true by Corollary \ref{affine} and Definition \ref{defsysenl}.) 

\begin{lem}\label{van}
Let $\cZ$ be an affine widening and let $E$ be a crystalline sheaf of 
$\Q \otimes_{\Z} {\cal O}_{\vec{\cZ}}$-modules.
Then we have $H^q(\vec{\cZ},E)=0$ for $q>0$. 
\end{lem}

We can define the functor 
$$\phi_{\vec{\cZ},*}: (\cX/\cY)^{\log,\sim}_{\conv,\Zar}\vert_{\cZ} 
  \lra \vec{\cZ}^{\sim} $$
as in \cite[p.91]{shiho2}. Then we have the equality 
$\phi_{\cZ,*}=\gamma_* \circ \phi_{\vec{\cZ},*}$ and 
the following commutative diagram: 

\begin{equation*}
\begin{CD}
(X/\cY)^{\log,\sim}_{\conv,\Zar}\vert_{\cZ} @>{\phi_{\vec{\cZ},*}}>>
\vec{\cZ}^{\sim} @>{\gamma_*}>> \cZ_{\Zar}^{\sim} \\
@V{j_{\cZ,*}}VV @. @\vert \\
(X/V)^{\log,\sim}_{\conv,\Zar} @>>{u_*}> X_{\Zar}^{\sim} @<<{z_*}< 
Z_{\Zar}^{\sim}. 
\end{CD}
\end{equation*}

By the same argument as in the absolute case, 
we have the following (\cite[2.3.1, 2.3.2]{shiho2}, 
\cite[0.4.1, 0.4.2]{ogus2}):

\begin{lem}\label{2.3.2}
Let the notations be as above and assume that $\cZ$ is an affine widening. 
Then$:$ 
\begin{enumerate}
\item 
The functor $\phi_{\vec{\cZ},*}$ sends an injective sheaf to a
flasque sheaf and exact. 
\item 
For an isocrystal $\cE$ on $(X/\cY)^{\log}_{\conv,\Zar}\vert_{\cZ}$, 
we have $H^q((X/\cY)^{\log}_{\conv,\Zar}\vert_{\cZ}, \cE) \isom
 H^q(\vec{\cZ},\allowbreak \phi_{\vec{\cZ},*}\cE)$ for all $q\geq 0$ and
 these groups vanish for $q>0$. 
\end{enumerate} 
\end{lem} 

\begin{pf} 
We omit the proof of (1), since the proof is the same as that of 
\cite[0.4.1, 0.4.2]{ogus2}. (2) is the immediate consequence of 
(1) and Lemma \ref{van}. 
\end{pf}

By using these lemmas, we can deduce the following, as 
\cite[2.3.3, 2.3.4]{shiho2}: 

\begin{prop}\label{vavvav}
Let the notations be as above and let 
$\cZ$ be an affine widening. Then, for an isocrystal $\cE$
in $(X/V)^{\log,\sim}_{\conv,\Zar}\vert_{\cZ}$, we have
$R^q(u \circ j_{\cZ})_*\cE=0$ and $R^qj_{\cZ,*}\cE=0$ for $q>0$.
\end{prop}

\begin{pf} 
We may assume $X$ is affine. $R^q(u \circ j_{\cZ})_*\cE$ is 
the sheaf associated to the presheaf 
$$U \mapsto H^q((U/\cY)^{\log}_{\conv,\Zar}\vert_{\cZ}, \cE)$$
 and it vanishes for $q>0$ by Lemma \ref{2.3.2}. \par 
To prove the vanishing $R^qj_{\cZ,*}\cE=0 \, (q>0)$,  
it suffices to show that the sheaf $(R^qj_{\cZ,*}\cE)_{\cZ'}$ on 
$\cZ'_{\Zar}$ induced by $R^qj_{\cZ,*}\cE$ vanishes  
 for any affine enlargement $\cZ'$ and $q>0$. 
Let $\cZ \times \cZ'$ be the direct product as widening (which is 
again affine), and consider the following diagram: 
{\scriptsize
\begin{equation*}
\begin{CD}
 (\cZ' \times \cZ)^{\ex,\sim}_{\Zar} @<{\gamma}<< 
 \vec{(\cZ' \times \cZ)^{\ex,\sim}}
@<{\phi_{\vec{(\cZ' \times \cZ)^{\ex}},*}}<<
 (X/\cY)^{\log,\sim}_{\conv,\Zar}\vert_{(\cZ' \times \cZ)^{\ex}} 
@>{j_{\cZ'}\vert_{\cZ}}>> 
 (X/V)^{\log,\sim}_{\conv,\Zar}\vert_{\cZ} \\
@V{\pr}VV @. @V{j_{\cZ}\vert_{\cZ'}}VV @V{j_{\cZ}}VV \\ 
{\cZ'}_{\Zar}^{\sim} @=  {\cZ'}_{\Zar}^{\sim} @<<{\phi_{\cZ'}}<
 (X/\cY)^{\log,\sim}_{\conv,\Zar}\vert_{\cZ'} 
@>>{j_{\cZ'}}>
 (X/\cY)^{\log,\sim}_{\conv,\Zar}.
\end{CD}
\end{equation*}}
We can prove that the functors $j_{\cZ'}^*, (j_{\cZ'}|_{\cZ})^*$ 
sends injectives and injectives (which can be shown as in \cite{ogus2}), 
and the functors 
$(j_{\cZ'}\vert_{\cZ})^*$, 
$\phi_{\cZ',*}$, $\phi_{\vec{(\cZ' \times \cZ)^{\ex}},*}$, $j_{\cZ'}^*$
 are exact. So we have 
\allowdisplaybreaks{
\begin{align*}
(R^qj_{\cZ,*}\cE)_{\cZ'} & = \phi_{\cZ',*}j^*_{\cZ'}R^qj_{\cZ,*}\cE \\ 
                   & = \phi_{\cZ',*}R^q(j_{\cZ}\vert_{\cZ'})_*
                           (j_{\cZ'}\vert_{\cZ})^*\cE \\
                   & = R^q(\pr \circ \gamma)_*
                       \phi_{\vec{(\cZ' \times \cZ)^{\ex}},*}
                       (j_{\cZ'}\vert_{\cZ})^*\cE
\end{align*}}
and the last term is equal to zero for $q>0$ by Lemma \ref{2.3.2}. 
(Note that, for any affine open $\cU' \subseteq \cZ'$ (which we can 
regard naturally as an enlargement), $\cU \times \cZ$ is an affine 
widening.) So we have the vanishing $R^qj_{\cZ,*}\cE = 0$ for $q > 0$. 
\end{pf} 

\begin{cor}\label{acyclic}
Let $\cZ$ be an affine widening and $\cE$ be an isocrystal on
$(\cX/\cY)^{\log,\sim}_{\conv,\Zar}\vert_{\cZ}$. Then $j_{\cZ,*}\cE$ is 
$u_*$-acyclic.
\end{cor}

\begin{pf} 
$$ Ru_*(j_{\cZ,*}\cE) = R(u \circ j_{\cZ})_* \cE = u_*j_{\cZ,*}\cE $$ 
by Proposition \ref{vavvav}. 
\end{pf}

Now let us consider the situation
\begin{equation}\label{embeddable}
\begin{CD}
(X,M_{X}) @>i>> (\cP,M_{\cP}) \\
@VfVV @VgVV \\
(\cY,M_{\cY}) @= (\cY,M_{\cY}), 
\end{CD}
\end{equation}
where $(X,M_X)$ is an object in $\L$, 
$f$ is a morphism in $\pLF$, $i$ is a closed immersion in 
$\pLF$ and $g$ is a formally log smooth morphism 
in $\pLF$. 
Then $\cP:=((\cP,M_{\cP}),(X,M_X),i,z)$ is a pre-widening 
of $(X,M_{X})/(\cY,M_{\cY})$. Denote the associated widening of 
$\cP$ by $\hat{\cP}$ and, 
for an isocrystal $\cE$ on 
$(X/\cY)^{\log}_{\conv}$, put 
$$ \omega^i_{\hat{\cP}}(\cE) := 
j_{\hat{\cP},*}(j^*_{\hat{\cP}}\cE \otimes_{\cO_{\cX/\cY}} 
\phi^*_{\hat{\cP}}(\omega^i_{\cP/\cY} \vert_{\hat{P}})). $$
Then we have the following theorem (relative version of 
\cite[2.3.5,2.3.6]{shiho2}), which should be called as the relative log 
convergent Poincar\'{e} lemma: 

\begin{thm}\label{lcpl1}
Let the notations be as above. Then, 
there exists a canonical structure of complex on 
$\omega^{\b}_{\hat{\cP}}(\cE)$ and the adjoint homomorphism 
$\cE \lra j_{\hat{\cP},*}j^*_{\hat{\cP}}\cE 
\allowbreak = \omega^0_{\hat{\cP}}(\cE)$ induces the quasi-isomorphism 
$\cE \overset{\simeq}{\lra} \omega^{\b}_{\hat{\cP}}(\cE).$ 
\end{thm} 

\begin{rem}
Note that, in the above theorem, $\cE$ is not assumed to be locally free: 
So it is slightly more general than \cite[2.3.5, 2.3.6]{shiho2} even 
in the absolute case. 
\end{rem}

\begin{pf} 
By definition, one can check the equality 
$$ \omega^i_{\hat{\cP}}(\cE)(\cZ) = 
\vpl_n (\cE_{T_n(\cZ \times \cP)} \otimes_{\cO_{T_n(\cZ \times \cP)}} 
\pi_{\cP,n}^*\omega^i_{\cP/\cY})(T_n(\cZ \times \cP)) $$ 
for an enlargement $\cZ$, where $\cZ \times \cP$ is the product of 
$\cZ$ and $\cP$ taken in the category of pre-widenings and $\pi_{\cP,n}$ 
is the canonical 
map $T_n(\cZ \times \cP) \lra \cP$. Varying $\cZ$, we see that 
the sheaf $\omega^i_{\hat{\cP}}(\cE)_{\cZ}$ on $\cZ_{\Zar}$ induced by 
$\omega^{\b}_{\hat{\cP}}$ is given by 
$$ \omega^i_{\hat{\cP}}(\cE)_{\cZ} = 
\vpl_n \pi_{\cZ,n,*}(\cE_{T_n(\cZ \times \cP)} 
\otimes_{\cO_{T_n(\cZ \times \cP)}} \pi_{\cP,n}^*\omega^i_{\cP/\cY}), $$
where $\pi_{\cZ,n}$ is the canonical map $T_n(\cZ \times \cP) \lra \cZ$. 
To define a canonical structure of a complex on $\omega^i_{\hat{\cP}}(\cE)$, 
it suffices to construct 
a canonical, functorial structure of a complex on 
$\omega^{\b}_{\hat{\cP}}(\cE)_{\cZ}$ for affine enlargements 
$\cZ := ((\cZ,M_{\cZ}),(Z,M_Z),i,z)$. 
Let $(\cP^m,M_{\cP^m})$ 
be the $m$-th log 
infinitesimal neighborhood of $(\cP,M_{\cP})$ in $(\cP,M_{\cP}) 
\times_{(\cY,M_{\cY})} (\cP,M_{\cP})$, let $p_{i,m}:(\cP^m, M_{\cP^m}) \lra 
(\cP,M_{\cP}) \, (i=1,2)$ be the morphism induced by the $i$-th 
projection and put $X_i^m := p_{i,m}^{-1}(X), 
M_{X_i^m} := M_X |_{X_i^m}$. Then we have pre-widenings 
$$ \cP^m := ((\cP^m,M_{\cP^m}), (X,M_X)), \quad 
\cP_i^m := ((\cP^m,M_{\cP^m}), (X_i^m,M_{X_i^m})), $$
and the following diagram of pre-widenings (for $i=1,2$): 
\begin{equation*}
\begin{CD}
\cP^m @>>> \cP_i^m \\
@VVV @V{p_{i,m}}VV \\
\cP @= \cP. 
\end{CD}
\end{equation*}
It induces the diagram of pre-widenings
\begin{equation}\label{zzzz}
\begin{CD}
\cP^m \times \cZ @>{r_i}>> \cP^m_i \times \cZ \\
@VVV @V{q_i}VV \\
\cP \times \cZ @= \cP \times \cZ 
\end{CD}
\end{equation}
for $i=1,2$. By taking the systems of universal enlargements of the 
above pre-widenings, we see that $r_i \, (i=1,2)$ induce 
the isomorphisms of inductive systems of $p$-adic fine log formal 
$\cB$-schemes 
\begin{equation}\label{hello}
\{T_n(\cP^m \times \cZ)\}_n \isom \{T_n(\cP_i^m \times \cZ) \}_n
\end{equation}
by Lemma \ref{zure} 
and that $q_i \, (i=1,2)$ induce the isomorphisms 
\begin{equation}\label{str-eq}
T_n(\cP_1^m \times \cZ) 
\isom T_n(\cP \times \cZ) \times_{\cP} \cP^m, \quad 
T_n(\cP_2^m \times \cZ) 
\isom \cP^m \times_{\cP} T_n(\cP \times \cZ) 
\end{equation}
by Lemma \ref{2.1.26}. 
By evaluating $\cE$ on $\{T_n(\cP \times \cZ)\}_n$, 
we see that $\cE$ naturally induces a coherent sheaf $E$ on 
$]Z[^{\log}_{\cP^m \times \cZ} = 
\bigcup_n T_n(\cP \times \cZ)_K$, and the isomorphisms 
\eqref{hello}, \eqref{str-eq} induce the isomorphisms 
$$ \theta_m: \pi_{\cP}^*\cO_{\cP^m} 
\otimes_{\cO_{]Z[^{\log}_{\cP \times \cZ}}} E 
\os{\simeq}{\lra} 
E \otimes_{\cO_{]Z[^{\log}_{\cP \times \cZ}}} \pi_{\cP}^*\cO_{\cP^m} 
\quad (m \in \N). $$
(Here $\pi_{\cP}^*$ denotes the functor 
$\Coh(\cO_{\cP}) \lra \Coh(\cO_{]Z[^{\log}_{\cP \times \cZ}})$ 
defined as composite 
$$\Coh(\cO_{\cP}) \lra \Coh(\Q \otimes \cO_{\cP}) \simeq 
\Coh(\cO_{\cP_K}) \lra \Coh(\cO_{]Z[^{\log}_{\cP \times \cZ}}),$$ 
where the last arrow is the functor induced by the map 
$]Z[^{\log}_{\cP \times \cZ} \lra \cP_K$.) Then we can define the map 
$\ti{d}: E \lra E \otimes_{\cO_{]Z[^{\log}_{\cP \times \cZ}}} 
\pi_{\cP}^*\omega^1_{\cP/\cY}$ by 
$\wt{d}(e) := \theta_1(1\otimes e) - e \otimes 1$, and extend it 
to the diagram 
$$ E @>{\wt{d}}>> E \otimes_{\cO_{]Z[^{\log}_{\cP \times \cZ}}} 
\pi_{\cP}^*\omega^1_{\cP/\cY} @>{\wt{d}}>> \cdots @>{\wt{d}}>> 
E \otimes_{\cO_{]Z[^{\log}_{\cP \times \cZ}}} 
\pi_{\cP}^*\omega^q_{\cP/\cY} @>{\wt{d}}>> \cdots $$ 
in standard way. By applying the direct image by the map 
$\pi_{\cZ}: ]Z[^{\log}_{\cP \times \cZ} \lra \cZ_K \lra \cZ$ 
to the above diagram, we obtain the diagram 
$$\omega^{\b}_{\hat{\cP}}(\cE)_{\cZ} := 
[\omega^0_{\hat{\cP}}(\cE)_{\cZ} @>d>> 
\omega^1_{\hat{\cP}}(\cE)_{\cZ} 
@>d>> \cdots @>d>> \omega^q_{\hat{\cP}}(\cE)_{\cZ} @>d>> \cdots ], 
$$
since we have 
$$ 
\pi_{\cZ,*} (E \otimes_{\cO_{]Z[^{\log}_{\cP \times \cZ}}} 
\pi_{\cP}^*\omega^q_{\cP/\cY}) = 
\vpl_n \pi_{\cZ,n,*}(\cE_{T_n(\cZ \times \cP)} 
\otimes_{\cO_{T_n(\cZ \times \cP)}} \pi_{\cP,n}^*\omega^q_{\cP/\cY})
$$ 
by definition. 
This construction is functorial with respect to affine enlargement 
$\cZ$ and so it induces the diagram 
$$ \omega^{\b}_{\hat{\cP}}(\cE) := 
[\omega^0_{\hat{\cP}}(\cE) @>d>> \omega^1_{\hat{\cP}}(\cE) 
@>d>> \cdots @>d>> \omega^q_{\hat{\cP}}(\cE) @>d>> \cdots ]. $$
We prove that the diagram $\omega^{\b}_{\hat{\cP}}(\cE)$ 
forms a complex and that the adjoint map 
$\cE \allowbreak \lra \allowbreak 
j_{\hat{\cP},*}j^*_{\hat{\cP}}\cE = \omega^0_{\hat{\cP}}(\cE)$ 
induces the quasi-isomorphism 
$ \cE \overset{\sim}{\lra} \omega^{\b}_{\hat{\cP}}(\cE)$. 
To prove it, it suffices to check it on an affine enlargement 
$\cZ := ((\cZ,M_{\cZ}),(Z,M_Z),i,z)$. Note first that we have the 
isomorphism 
\begin{align*}
\omega^q_{\hat{\cP}}(\cE)_{\cZ} & = 
\vpl_n\pi_{\cZ,n,*}(\cE_{T_n(\cZ \times \cP)} \otimes 
\pi_{\cP,n}^*\omega^q_{\cP/\cY}) \\ & = 
\vpl_n\pi_{\cZ,n,*}(\pi_{\cZ,n}^*\cE_{\cZ} \otimes 
\pi_{\cP,n}^*\omega^q_{\cP/\cY}) \\ & \cong 
\cE_{\cZ} \otimes_{\cO_{\cZ}} \vpl_n \pi_{\cZ,n,*}\pi_{\cP,n}^*
\omega^q_{\cP/\cY} = 
\cE_{\cZ} \otimes_{\cO_{\cZ}} \omega^q_{\hat{\cP}}(\cK_{X/\cY})_{\cZ}: 
\end{align*}
Indeed, we see that the functor 
$\vpl_n\pi_{\cZ,n,*}(\pi_{\cZ,n}^*(-) \otimes 
\pi_{\cP,n}^*\omega^q_{\cP/\cY})$ is exact on $\Coh(\Q \otimes \cO_{\cZ})$ 
by using Lemma \ref{add1} below and \cite[5.10.1]{bc}, and the functor 
$(-)\otimes_{\cO_{\cZ}} \vpl_n \pi_{\cZ,n,*}\pi_{\cP,n}^*
\omega^q_{\cP/\cY}$ is right exact. So it suffices to check the 
isomorphism in the case $\cE_{\cZ} \cong \Q \otimes_{\cZ} 
\cO_{\cZ}^{\oplus r}$, which is obvious. 
Then, 
by noting the fact that the diagram (\ref{zzzz}) is defined over $\cZ$ and 
by definition of $d$, 
we see that the map $d$ on $\omega^{\b}_{\hat{\cP}}(\cE)_{\cZ}$ is 
compatible with the map $\id \otimes d$ on 
$\cE_{\cZ} \otimes_{\cO_{\cZ}} \omega^{\b}_{\hat{\cP}}(\cK_{X/\cY})_{\cZ}$ 
via the above isomorphism. So it suffices to show that the diagram 
$\omega^{\b}_{\hat{\cP}}(\cK_{X/\cY})_{\cZ}$ forms a complex and that 
the diagram $\Q \otimes_{\Z} \cO_{\cZ} \lra \omega^{\b}_{\hat{\cP}}(\cK_{X/\cY})_{\cZ}$ 
is a complex which is locally homotopic to zero. \par 
By construction, the map 
$$d: \omega^{0}_{\hat{\cP}}(\cK_{X/\cY})_{\cZ} = 
\vpl_n \pi_{\cZ,n,*}(\pi_{\cP,n}^*\cO_{\cP})
\lra \vpl_n \pi_{\cZ,n,*}(\pi_{\cP,n}^*\omega^1_{\cP/\cY}) = 
\omega^{1}_{\hat{\cP}}(\cK_{X/\cY})_{\cZ}$$ is obtained by applying 
$\pi_{\cZ,*}$ to the relative differential 
$$ \ti{d}: \cO_{]Z[^{\log}_{\cP \times_{\cY} \cZ}} \lra 
\omega^1_{]Z[^{\log}_{\cP \times_{\cY} \cZ}/\cZ_K}. $$
Since $\ti{d} \circ \ti{d} = 0$ holds, we have $d \circ d = 0$. 
So the diagram $\omega^{\b}_{\hat{\cP}}(\cK_{X/\cY})$ forms a complex. 
Now let $\DR$ be the relative de 
Rham complex 
$$ \cO_{]Z[^{\log}_{\cP \times_{\cY} \cZ}} @>{\ti{d}}>> 
\omega^1_{]Z[^{\log}_{\cP \times_{\cY} \cZ}/\cZ_K} 
@>{\ti{d}}>> 
\omega^2_{]Z[^{\log}_{\cP \times_{\cY} \cZ}/\cZ_K} 
@>{\ti{d}}>> \cdots. $$
Then the diagram $\Q \otimes_{\Z} \cO_{\cZ} 
\lra \omega^{\b}_{\hat{\cP}}(\cK_{X/\cY})_{\cZ}$ is the same as the diagram  
$$ \sp_*\cO_{]Z[_\cZ} \lra \pi_{\cZ,*}\DR. $$ 
So it suffices to prove that this is a complex locally homotopic to zero. 
To prove this, we may assume that 
$]Z[^{\log}_{\cP \times_{\cY} \cZ}$ is isomorphic to $]Z[_{\cZ} \times D^r 
= \cZ_K \times D^r$, 
where $D^r$ is the $r$-dimensional open polydisc of radius $1$ (where 
$r$ is the rank of $\omega^1_{\cP/\cY}$), by Lemma \ref{add1} below. 
In this case, we can construct 
the desired homotopy by Lemma \ref{add2} below, 
as in the case of Poincar\'{e} lemma for a complex open 
polydisc. So the proof is finished 
(modulo the proof of Lemmas \ref{add1}, \ref{add2} 
below). 
\end{pf}

We give a proof of the following two lemmas which are used in the 
above proof: 

\begin{lem}\label{add1} 
Let $(X,M_X)$ be a fine log $\cB$-scheme and assume we are given the 
commutative diagram 
\begin{equation}\label{add11}
\begin{CD}
(X,M_X) @>>> (\cP,M_{\cP}) \\ 
@\vert @VfVV \\ 
(X,M_X) @>>> (\cQ,M_{\cQ}), 
\end{CD}
\end{equation}
where $(\cP,M_{\cP}), (\cQ,M_{\cQ})$ are $p$-adic fine log formal 
$\cB$-schemes, horizontal arrows are closed immersions and 
$f$ is formally log smooth. Then, Zariski locally on $\cQ$, 
we have the isomorphism 
$]X[^{\log}_{\cP} \cong ]X[^{\log}_{\cQ} \times D^r$ for some $r$ 
$($where $D^r$ is the $r$-dimensional open polydisc of radius $1)$ 
such that the morphism $]X[^{\log}_{\cP} \lra ]X[^{\log}_{\cQ}$ 
induced by $f$ is 
identified with the first projection 
$]X[^{\log}_{\cQ} \times D^r \lra ]X[^{\log}_{\cP}$ via the isomorphism. 
\end{lem}

\begin{pf} 
Let $\cP^{\ex}, \cQ^{\ex}$ be the exactification of the horizontal arrows 
of the diagram \eqref{add11}. Then we have the diagram 
\begin{equation*}
\begin{CD}
X @>>> \cP^{\ex} \\ 
@\vert @V{f^{\ex}}VV \\ 
X @>>> \cQ^{\ex}. 
\end{CD}
\end{equation*}
We can easily check that $f^{\ex}$ is formally smooth (that is, satisfies 
infinitesimal lifting property for nilpotent closed immersions of 
affine schemes). To prove the lemma, we may assume 
$\cQ^{\ex} =: \Spf A, \cP^{\ex} =: \Spf B$ are affine. Let $f^*: A \lra B$ be 
the ring homomorphism associated to $f^{\ex}$. By infinitesimal 
lifting property, there exists a ring homomorphism $s:B \lra A$ 
satisfying $s \circ f^* = \id$. Put $I := \Ker\, s$. Then, by 
\cite[19.5.3]{gd4}, $I/I^2$ is finitely generated formally projective 
$A$-module and we have the canonical isomorphism 
$\varphi_n: {\rm Sym}^n_A(I/I^2) \cong {\rm gr}^n_I B \, (n \in \N)$. 
By localizing $A$, we may assume that $I/I^2$ is a free $A$-module of 
finite rank. 
The canonical inclusion $I/I^2 \hra B/I^2$ lifts by formal projectivity to 
a homomorphism $I/I^2 \lra B/I^{n+1}$ and it induces the ring homomorphism 
$\psi_n: {\rm Sym}^{\b \leq n}_A(I/I^2) \lra B/I^{n+1}$. This is isomorphism 
since $\varphi_n$'s are isomorphisms. By taking inverse limit, we obtain 
the ring isomorphism 
$\varprojlim_n({\rm Sym}^{\b \leq n}_A(I/I^2)) \os{\cong}{\lra} B$, 
and the left hand side is isomorphic to 
$A[[x_1, \cdots, x_r]] = 
\varprojlim_n A[x_1, \cdots, x_r]/(x_1,\cdots,x_r)^{n+1}$. 
(Here $r$ is the rank of $I/I^2$.) So we have 
$$ ]X[^{\log}_{\cP} = (\Spf B)_K = (\Spf A[[x_1, \cdots, x_r]])_K = 
(\Spf A)_K \times D^r = ]X[^{\log}_{\cQ} \times D^r. $$
\end{pf} 

\begin{lem}\label{add2}
Let $\cZ := \Spf A_0$ be a $p$-adic affine formal $\cB$-scheme and 
$\pi:\cZ_K \times D^r \lra \cZ_K$ the projection. 
$($Here $D^r$ denotes the $r$-dimensional open unit polydisc of radius $1$.$)$
Then the complex 
\begin{equation}\label{add21}
\sp_*\cO_{\cZ_K} \lra 
(\sp \circ \pi)_*\Omega^{\b}_{\cZ_K \times D^r/\cZ_K} 
\end{equation}
is homotopic to zero. 
\end{lem} 

\begin{pf} 
Let us put $A:=\Q \otimes_{\Z} A_0$, 
$B := \Gamma(\cZ_K \times D^r, \cO_{\cZ_K \times D^r})$. 
Let $x_1, \cdots, x_r$ be the coordinate of $D^r$ and for $q \geq 0$, 
let us define $B^q$ by 
$$ B^q := \bigoplus_{1 \leq i_1 < \cdots < i_q \leq r} 
B dx_{i_1} \wedge dx_{i_2} \wedge \cdots \wedge dx_{i_q}. $$
Then, if we apply $\Gamma(\cZ,-)$ to the complex \eqref{add21}, 
we obtain the complex $\cC := [A \ra B^0 \ra B^1 \ra \cdots]$. 
To prove the lemma, it suffices to show that the complex $\cC$ is 
homotopic to zero via a homotopy which is functorial on $A$. \par 
For $1 \leq j \leq r$ and $q \geq 1$, let 
$\alpha^q_j: B^q \lra B^{q-1}$ be the map defined by 
\begin{align*}
& \phantom{=} 
\alpha^q_j(f(x_1,x_2,\cdots,x_r)dx_{i_1} \wedge \cdots \wedge dx_{i_q}) \\ 
& := 
\left\{ 
\begin{aligned}
& 0, \qquad \text{if $j \notin \{i_1, \cdots, i_q\}$}, \\
& (-1)^{k-1}
\left( \int_0^{x_j}f(x_1, \cdots, \underset{\widehat{j}}{t}, \cdots x_r)dt 
\right) dx_{i_1} \wedge \cdots \wedge 
\widehat{dx_{i_k}} \wedge \cdots \wedge dx_{i_q}, \qquad 
\text{if $j = i_k$}, 
\end{aligned}
\right. 
\end{align*}
and let $\alpha^0_j:B^0 \lra A$ be the zero map. 
Let $\beta_j:\cC \lra \cC$ be the map defined by 
$\beta := \id - (d \circ \alpha_j^{\b}+\alpha_j^{\b+1} \circ d)$, where 
$d$ denotes the differential of $\cC$. Then, by definition, $\beta_j$ is 
homotopic to $\id$. Moreover, 
$
\beta_j(B dx_{i_1} \wedge dx_{i_2} \wedge \cdots \wedge dx_{i_q}) 
$ is zero if $j=i_k$ for some $1 \leq k \leq q$ and we have 
$\beta_j(B dx_{i_1} \wedge dx_{i_2} \wedge \cdots \wedge dx_{i_q}) 
\subseteq  B dx_{i_1} \wedge dx_{i_2} \wedge \cdots \wedge dx_{i_q}
$ in general. We have $\beta_j(a) = a$ for $a \in A$ and 
$\beta_j(f(x_1, \cdots, x_r)) = f(x_1, \cdots, \underset{\widehat{j}}{0}, 
\cdots, x_r)$ for $f \in B^0$. \par 
Next let us define $\beta:\cC \lra \cC$ by $\beta := \beta_{r} \circ 
\beta_{r-1} \circ \cdots \circ \beta_1$. Then $\beta$ is homotopic to 
$\id$ and we have $\beta(B^q)=0 \, (q \geq 1)$, $\beta(a)=a\,(a \in A)$, 
$\beta(f(x_1,\cdots, x_r))=f(0,\cdots,0)\,(f \in B^0)$. 
Now we define $\gamma^q:B^q \lra B^{q-1}\,(q \geq 1)$, $\gamma^0:B^0 \lra A$ 
by $\gamma^q=0\,(q \geq 1), \gamma^0(f(x_1,\cdots, x_r)) := f(0,\cdots,0)$. 
Then $\{\gamma^q\}_q$ gives a homotopy between $\beta$ and zero map on 
$\cC$. So the identity map on $\cC$, being homotopic to $\beta$, is 
homotopic to zero map. So the complex $\cC$ is homotopic to zero. 
Moreover, it is easy to see that the homotopies we constructed are 
functorial on $A$. So we are done. 
\end{pf} 

Let $\ti{u}: (X/\cY)^{\log,\sim}_{\conv,\et} \lra X^{\sim}_{\Zar}$ be the 
composite 
$(X/\cY)^{\log,\sim}_{\conv,\et} \os{\epsilon}{\lra} 
(X/\cY)^{\log,\sim}_{\conv,\Zar} \os{u}{\lra}
 X^{\sim}_{\Zar}$. Then we can prove the following in 
the same way as \cite[2.3.8]{shiho2}: 

\begin{cor}\label{lcpl3}
Assume the diagram \eqref{embeddable} is given and 
let $\cE$ be an 
isocrystal on $(X/\cY)^{\log}_{\conv}$.  
Then we have the canonical quasi-isomorphism 
$$ R\ti{u}_*\cE = \sp_* \DR(]X[^{\log}_{\cP}/\cY_K, \cE). $$
\end{cor} 

\begin{pf}
Here we only sketch the proof. 
Since we may work Zariski locally, we may assume that 
$((\cP,M_{\cP}),(X,M_X))$ is an affine widening. 
By Proposition \ref{zid} and Theorem \ref{lcpl1}, 
we have the quasi-isomorphism 
$$ R\ti{u}_*\cE \os{\sim}{\lra} 
Ru_* \omega^{\b}_{\hat{\cP}}(\epsilon_*\cE). $$ 
Next, by Proposition \ref{acyclic}, we have the quasi-isomorphism 
\begin{align*}
Ru_*\omega^{\b}_{\hat{\cP}}(\epsilon_*\cE) & = 
Ru_* j_{\hat{\cP},*}(\phi^*_{\hat{\cP}}(\omega^{\b}_{\cP/\cY} 
\vert_{\hat{\cP}}) \otimes j^*_{\hat{\cP}} \epsilon_*\cE) \\
& \os{=}{\lla} 
u_* j_{\hat{\cP},*}(\phi^*_{\hat{\cP}}(\omega^{\b}_{\cP/\cY} 
\vert_{\hat{\cP}}) \otimes j^*_{\hat{\cP}} \epsilon_*\cE) 
= u_*\omega^{\b}_{\hat{\cP}}(\epsilon_*\cE). 
\end{align*}
Finally, we can prove the isomorphism (not only the quasi-isomorphism)
$$ \ti{u}_* \omega^{\b}_{\hat{\cP}}(\epsilon_*\cE) 
\isom \sp_* \DR(]X[^{\log}_{\cP}, \cE) $$
in the same way as \cite[2.3.8]{shiho2}. Combining these, we obtain the 
assertion. 
\end{pf}

We can prove the following in the same way as 
\cite[2.3.9]{shiho2} (we omit the proof). 

\begin{cor}\label{lcpl4}
Let $(X,M_X)$ be an object in $\L$, let 
$f:(X,M_{X}) \lra (\cY,M_{\cY})$ be a morphism in $\pLF$ and 
let $\cE$ be an isocrystal on 
$(X/\cY)^{\log}_{\conv}$. If we take an embedding system 
$$ (X,M_{X}) \overset{g}{\lla} (X^{(\bullet)}, M_{X^{(\bullet)}}) 
\overset{i}{\hra} (\cP^{(\bullet)},M_{\cP^{(\bullet)}}), $$
we have the isomorphism 
$$ R^qf_{X/\cY,\conv *}\cE \isom 
R^q(f\circ g)_* \sp_*\DR(]X^{(\b)}[^{\log}_{\cP^{(\b)}}/\cY_K,\cE). $$
\end{cor} 


Now we would like to compare relative log convergent cohomology and 
relative log crystalline cohomology when we are given a diagram as in 
\eqref{diag-p} with $f$ log smooth. 
To do this, first we construct 
a functor from the category of isocrystals on relative 
log convergent site to that on log crystalline site: 

\begin{prop}
Assume we are given the diagram 
$(X,M_X) \os{f}{\lra} (Y,M_Y) \os{\iota}{\hra} (\cY,M_{\cY})$ as in 
\eqref{diag-p} such that $f$ is log smooth. 
Then we have the canonical functor 
$$ \Phi: I_{\conv}((X/\cY)^{\log}) \lra I_{\crys}((X/\cY)^{\log}) $$
such that $\Phi$ sends locally free isocrystals on $(X/\cY)^{\log}_{\conv}$ 
to locally free isocrystals on $(X/\cY)^{\log}_{\crys}$. 
\end{prop}

\begin{pf}
We only sketch the outline because the proof is similar to 
\cite[5.3.1]{shiho1} (although the proof here is much simpler.) \par 
First we assume that $(X,M_X)$ admits a closed immersion 
$(X,M_X) \hra (\cP,M_{\cP})$ into a $p$-adic fine log formal $\cB$-scheme 
which is formally log smooth over $(\cY,M_{\cY})$ and we construct 
a functor $\Psi: I_{\conv}((X/\cY)^{\log}) \lra 
\HPDI((X \hra \cP/\cY)^{\log})$. 
Denote the $(i+1)$-fold fiber product of $(\cP,M_{\cP})$ over 
$(\cY,M_{\cY})$ by $(\cP(i),M_{\cP(i)})$ and let 
$(\cP(i)^{\ex},M_{\cP(i)^{\ex}})$ be the exactification of the closed 
immersion $(X,M_X) \hra (\cP(i),M_{\cP(i)})$. 
Put $\cI(i):= \Ker (\cO_{\cP(i)^{\ex}} \lra \cO_X)$, 
let $B_n(i)$ be the formal blow-up of $\cP(i)^{\ex}$ with respect to 
the ideal $p\cO_{\cP(i)^{\ex}} + \cI(i)^n$, let $T'_n(i)$ be the open 
sub formal scheme of $B_n(i)$ defined as the set of points 
$x \in B_n(i)$ satisfying 
$(p\cO_{\cP(i)^{\ex}}+\cI(i))\cO_{B_n(i),x} = p\cO_{B_n(i),x}$ and let 
$T_n(i)$ be the closed sub formal scheme of $T_n(i)$ defined by the ideal 
$\{x \in \cO_{T_n(i)} \,\vert\, p^nx=0 \,\text{for some $n>0$}\}$. 
(Then $\{T_n(i)\}_n$ is a system of universal enlargements of 
the pre-widening 
$((\cP(i),M_{\cP(i)}),(X,M_X))$. See Definition \ref{defsysenl}.) 
On the other hand, let $D(i)$ be 
the $p$-adically completed log PD-envelope of $(X,M_X)$ in 
$(\cP(i),M_{\cP(i)})$. 
Then, if $\Ker(\cO_{\cP^{\ex}(i)} \lra \cO_X)$ 
is generated by $m$ elements $(i=0,1,2)$, 
there exist canonical diagrams 
$D(i) \lra T'_n(i) \hookleftarrow T_n(i) \, (i=0,1,2)$ for 
$n\allowbreak = \allowbreak (p-1)m+1$, 
where the second map is the canonical closed immersion. Since we have 
the canonical equivalences of categories 
$\Coh(\Q \otimes \cO_{T_n(i)}) \cong \Coh(\Q \otimes \cO_{T'_n(i)})$, 
we can define, 
by the `pull-back by $D(i) \lra T'_n(i)$', 
the functor 
$$
\Psi: I_{\conv}((X/\cY)^{\log}) \os{\sim}{\lra}
\Strat''((X \hra \cP/\cY)^{\log}) \lra 
\HPDI((X \hra \cP/\cY)^{\log}). $$

Next we construct the desired functor $\Phi$ 
in general case. Take an embedding system 
$$ (X,M_{X}) \overset{g^{(\b)}}{\lla} (X^{(\bullet)}, M_{X^{(\bullet)}}) 
\overset{i}{\hra} (\cP^{(\bullet)},M_{\cP^{(\bullet)}}) $$
over $(X,M_X)/(\cY,M_{\cY})$ such that $X^{(0)}$ is affine, 
$(\cP^{(0)},M_{\cP^{(0)}}) \times_{(\cY,M_{\cY})} (Y,M_Y) = (X^{(0)}, 
M_{X^{(0)}})$ holds 
 and that 
$g$ is a Zariski \v{C}ech hypercovering. 
(The existence of such an embedding system 
follows from \cite[3.14]{kkato}.) 
Let 
$\Delta_{\conv}$ be the category of descent data on 
$I_{\conv}((X^{(i)}/\cY)^{\log})\,(i=0,1,2)$ (that is, the category of 
objects in $I_{\conv}((X^{(0)}/\cY)^{\log})$ endowed with 
isomorphism of glueing in $I_{\conv}((X^{(1)}/\cY)^{\log})$ satisfying 
the cocycle condition in $I_{\conv}(X^{(2)}/\cY)^{\log}).$)
Similarly, let $\Delta_{\crys}$ be the 
category of descent data on $I_{\crys}((X^{(i)}/\cY)^{\log}) 
\allowbreak \, \allowbreak (i=0,1,2)$ and 
let $\Delta_{\HPDI}$ be the 
category of descent data on 
$\HPDI((X^{(i)} \hra \cP^{(i)}/\cY)^{\log})\allowbreak \, 
\allowbreak (i=0,1,2)$. 
Then we have the diagram 
$$ \Delta_{\conv} \os{\Psi}{\lra} \Delta_{\HPDI} 
\os{\Lambda}{\lla} \Delta_{\crys} $$
and $\Lambda$ is an equivalence of categories because the functor 
$$ I_{\crys}((X^{(n)}/\cY)^{\log}) \lra 
\HPDI((X^{(n)} \hra \cP^{(n)}/\cY)^{\log}) $$ 
is an equivalence of catogories if $n=0$ and fully faithful in general. 
(See Proposition \ref{crys-hpdi} and the paragraph before it.) 
So we can define the desired functor $\Phi$ by 
$$ \Phi: I_{\conv}((X/\cY)^{\log}) \os{\sim}{\lra} \Delta_{\conv} 
\os{\Psi}{\lra} \Delta_{\HPDI} \os{\Lambda^{-1}}{\lra} \Delta_{\crys} 
\os{\sim}{\lra} I_{\crys}((X/\cY)^{\log}). $$

We can also check that the functor $\Phi$ defined above is independent 
of the choice of the embedding system chosen above, by a standard argument 
(cf. \cite[5.3.1 Step 3]{shiho1}). Hence we have constructed the 
functor $\Phi$. It is easy to see that $\Phi$ preserves the local freeness. 
So we are done. 
\end{pf} 

Now we prove the comparison theorem between relative log convergent 
cohomology and relative log crystalline cohomology. 
Assume we are given the diagram \eqref{diag-p} such that $f$ is log smooth 
and assume for the moment that 
$(X,M_X)$ admits a closed immersion 
$(X,M_X) \hra (\cP,M_{\cP})$ into a $p$-adic fine log formal $\cB$-scheme 
which is formally log smooth over $(\cY,M_{\cY})$. 
Let $\cE$ be an 
isocrystal on $(X/\cY)^{\log}_{\conv}$ and let $\Phi(\cE)$ be the 
associated isocrystal on $(X/\cY)^{\log}_{\crys}$. Then we have 
the log de Rham complex 
$\DR(]X[^{\log}_{\cP}/\cY_K, \cE)$ on $]X[^{\log}_{\cP}$ 
associated to $\cE$. On the other hand, if we denote 
the $p$-adically completed log PD-envelope of $(X,M_X) \hra (\cP,M_{\cP})$ 
by $(D,M_D)$, we have 
the log de Rham complex $\DR(\D/\cY,\Phi(\cE))$ associated to $\Phi(\cE)$. 
Moreover, by the same method as 
\cite[3.1.3]{shiho2}, we have the canonical morphism of complexes 
$$ \sp_*\DR(]X[^{\log}_{\cP}/\cY_K, \cE) \lra \DR(\cD/\cY,\Phi(\cE)). $$
Since we have the quasi-isomorphisms 
$Rf_*\sp_*\DR(]X[^{\log}_{\cP}/\cY_K, \cE) = Rf_{X/\cY,\conv *}\cE$ and 
$Rf_* \allowbreak \DR(\cD/\cY,\Phi(\cE)) \allowbreak = \allowbreak 
Rf_{X/\cY,\crys *}\Phi(\cE)$, 
we have the map 
$Rf_{X/\cY,\conv *}\cE \lra Rf_{X/\cY,\crys *}\Phi(\cE)$ 
in this case. Even if $(X,M_X)$ does not admit the closed immersion 
$(X,M_X) \hra (\cP,M_{\cP})$ as above, we can define the map 
$Rf_{X/\cY,\conv *}\cE \lra Rf_{X/\cY,\crys *}\Phi(\cE)$ 
by taking an embedding system. Then we have the following: 

\begin{thm}\label{conv-crys} 
Assume we are given the diagram \eqref{diag-p} with $(X,M_X)$ log smooth 
over $(Y,M_Y)$ and let $\cE$ be a locally free 
isocrystal on $(X/\cY)^{\log}_{\conv}$. 
Then the homomorphism 
$$ Rf_{X/\cY,\conv *}\cE \lra Rf_{X/\cY,\crys *}\Phi(\cE) $$ 
defined above is a quasi-isomorphism. 
\end{thm}

\begin{pf}
It suffices to prove that, under the condition of the existence 
of the closed immersion 
$(X,M_X) \hra (\cP,M_{\cP})$ as above, the homomorphism 
$Rf_*\sp_*\DR(]X[^{\log}_{\cP}/\cY_K, \cE) \lra Rf_*
\DR(\cD/\cY,\Phi(\cE))$ is a 
quasi-isomorphism. We can reduce it to the case where 
$(\cP,M_{\cP}) \times\allowbreak {}_{(\cY,M_{\cY})} (Y,M_Y) = 
(X,M_X)$ holds, because 
both hand sides 
are known to be independent of the choice 
of the closed immersion $(X,M_X) \hra (\cP,M_{\cP})$ up to 
quasi-isomorphism. In this case, the both hand sides are the same. 
\end{pf}


As a corollary, we have the following, which are the starting points of 
the argument in the next section: 

\begin{cor}\label{perf1}
Assume we are given the diagram 
\begin{equation*}
(X,M_X) \os{f}{\lra} (Y,M_Y) \os{\iota}{\hra} (\cY,M_{\cY}), 
\end{equation*}
where $f$ is a proper log smooth integral morphism in $\LB$, 
 $(\cY,M_{\cY})$ is an object in $\pLF$ and $\iota$ is the exact closed 
immersion defined by the ideal sheaf $p{\cal O}_{\cY}$. Then, for 
a locally free isocrystal $\cE$ on $(X/\cY)^{\log}_{\conv}$, 
$Rf_{X/\cY,\conv *}\cE$ is a perfect complex of 
$\Q \otimes_{\Z} \cO_{\cY}$-modules on $\cY_{\Zar}$. 
\end{cor}

\begin{pf}
It is immediate from Theorems \ref{perfcrys} and 
\ref{conv-crys}. 
\end{pf}

\begin{cor}\label{bccrys}
Assume we are given a diagram 
\begin{equation}
\begin{CD}
(X',M_{X'}) @>>> (Y',M_{Y'}) @>>> (\cY',M_{\cY'}) \\
@VVV @VVV @V{\varphi}VV \\
(X,M_X) @>f>> (Y,M_Y) @>{\iota}>> (\cY,M_{\cY}), 
\end{CD}
\end{equation}
where 
$f$ is a proper log smooth integral morphism in $\LB$, 
$\iota$ is the exact closed immersion defined by 
the ideal sheaf $p\cO_{\cY}$ and the squares are Cartesian. 
Then, for a 
locally free isocrystal $\cE$ on $(X/\cY)^{\log}_{\conv}$, 
we have the quasi-isomorphism 
$$ 
L\varphi^*
Rf_{X/\cY,\conv *}\cE \os{\sim}{\lra} 
Rf_{X'/\cY',\conv *} \varphi^*\cE. $$
\end{cor}

\begin{pf}
It is immediate from Theorems \ref{crysbc} and \ref{conv-crys}. 
\end{pf} 


\section{Relative log convergent cohomology (II)}

Assume we are given a diagram 
\begin{equation}\label{diag-s2}
(X,M_X) \os{f}{\lra} (Y,M_Y) \os{\iota}{\hra} (\cY,M_{\cY}), 
\end{equation}
where $f$ is a proper log smooth integral morphism in $\LB$, 
 $(\cY,M_{\cY})$ is an object in $\pLF$ and $\iota$ is a homeomorphic 
exact closed immersion in $\pLF$. 
In this section, we prove the coherence and the base change property 
of relative log convergent cohomology of $(X,M_X)$ over 
$(\cY,M_{\cY})$ when `$f$ has log smooth parameter' 
(for definition, see Definition \ref{parameter}). 
Note that, in the above situation, $\iota$ is not assumed 
to be defined by $p\cO_{\cY}$. So $\iota$ does not necessarily 
admit a PD-structure. Nevertheless, we can prove 
them by reducing to Corollaries \ref{perf1}, \ref{bccrys} (that is, the case 
where $\iota$ admits canonical PD-structure). \par 
The first proposition we need is 
the topological invariance of the category of isocrystals on 
log convergent site and that of 
relative log convergent cohomology. 

\begin{prop}\label{inv}
Assume we are given an object $(\cY,M_{\cY})$ 
in $\pLF$ and a homeomorphic exact closed immersion 
$i: (X_1,M_{X_1}) \hra (X_2,M_{X_2})$ in $\L$ over $(\cY,M_{\cY})$. 
Then$:$ 
\begin{enumerate}
\item 
The restriction functor 
$$ I_{\conv}((X_2/\cY)^{\log}) \lra I_{\conv}((X_1/\cY)^{\log}) $$ 
is an equivalence of categories. 
\item 
For $\cE \in I_{\conv}((X_2/\cY)^{\log})$, the restriction 
$$ Rf_{X_2/\cY,\conv *}\cE {\lra} 
Rf_{X_1/\cY,\conv *}i^*\cE $$
is a quasi-isomorphism. 
\end{enumerate}
\end{prop}

\begin{pf}
To prove (1), we may assume the existence of the closed immersion 
$(X_2, M_{X_2}) \hra (\cP,M_{\cP})$ into a $p$-adic fine log 
formal scheme $(\cP,M_{\cP})$ which is formally log smooth over 
$(\cY,M_{\cY})$. Let $(\cP(i),M_{\cP(i)})$ be the $(i+1)$-fold 
fiber product of $(\cP,M_{\cP})$ over $(\cY,M_{\cY})$. 
Then, since the natural morphism 
$]X_1[^{\log}_{\cP(i)} \lra ]X_2[^{\log}_{\cP(i)}$ is an isomorphism 
by definition of log tubular neighborhood, the restriction functor 
$$ \Strat''((X_2 \hra \cP/\cY)^{\log}) \lra 
\Strat''((X_1 \hra \cP/\cY)^{\log}) $$ 
is an equivalence of categories. The assertion (1) follows from this. \par 
Let us prove the assertion (2). 
Take an embedding system 
$$ (X_2,M_{\cX_2}) \lla (X_2^{(\bullet)},
 M_{X_2^{(\bullet)}}) 
\hra (\cP^{(\bullet)},M_{\cP^{(\bullet)}}) $$
of $(X_2,M_{\cX_2})/(\cY,M_{\cY})$ and put 
$(X_1^{(\bullet)},
 M_{X_1^{(\bullet)}}) := (X_2^{(\bullet)},
 M_{X_2^{(\bullet)}}) \times_{(X_2,M_{X_2})} (X_1,M_{X_1})$. Then 
the diagram 
$$ (X_1,M_{\cX_1}) \lla (X_1^{(\bullet)},
 M_{X_1^{(\bullet)}}) 
\hra (\cP^{(\bullet)},M_{\cP^{(\bullet)}}) $$
is also an embedding system. Then, since 
the natural morphism $]X_1^{(\b)}[^{\log}_{\cP^{(\b)}} \lra 
]X_2^{(\b)}[^{\log}_{\cP^{(\b)}}$ is an isomorphism, 
we have the isomorphism of the associated log de Rham complexes 
$$ 
\DR(]X_2^{(\b)}[^{\log}_{\cP^{(\b)}}/\cY_K, \cE) \os{=}{\lra} 
\DR(]X_1^{(\b)}[^{\log}_{\cP^{(\b)}}/\cY_K, i^*\cE). 
$$ 
The assertion of the lemma follows from this and Corollary \ref{lcpl4}. 
\end{pf}

By combining Proposition \ref{inv} and Corollaries \ref{perf1}, \ref{bccrys}, 
we immediately obtain the following: 

\begin{cor}\label{perf2}
Assume we are given the diagram \eqref{diag-s2} and a locally free 
isocrystal $\cE$ on 
$(X/\cY)^{\log}_{\conv}$. Put 
$(Y_1, M_{Y_1}) := (\cY,M_{\cY}) \otimes_{\Z_p} \Z/p\Z$ and 
assume that there exists a proper 
log smooth integral morphism $(X_1,M_{X_1}) \lra (Y_1, M_{Y_1})$ 
satisfying $(X_1,M_{X_1}) \times_{(Y_1,M_{Y_1})} (Y,M_Y) = \allowbreak 
(X, \allowbreak M_X)$. Then 
$Rf_{X/\cY,\conv *}\cE$ is a perfect complex 
of $\Q \otimes_{\Z} \cO_{\cY}$-modules. 
\end{cor}

\begin{cor}\label{bc1}
Assume we are given a diagram 
\begin{equation}
\begin{CD}
(X',M_{X'}) @>>> (Y',M_{Y'}) @>>> (\cY',M_{\cY'}) \\
@VVV @VVV @V{\varphi}VV \\
(X,M_X) @>>> (Y,M_Y) @>>> (\cY,M_{\cY}), 
\end{CD}
\end{equation}
where 
the horizontal lines are as in the diagram \eqref{diag-s2} and 
the left square is Cartesian. Let 
$(Y_1,M_{Y_1}) := (\cY,M_{\cY}) \otimes_{\Z_p} \Z/p\Z$ 
and assume that there exists a proper 
log smooth integral morphism $(X_1,M_{X_1}) \lra (Y_1, M_{Y_1})$ 
satisfying $(X_1,M_{X_1}) \times_{(Y_1,M_{Y_1})} (Y,M_Y) = (X,M_X)$. 
Then, for a locally free isocrystal $\cE$ on 
$(X/\cY)^{\log}_{\conv}$, 
we have the canonical quasi-isomorphism 
$$ L\varphi^*Rf_{X/\cY,\conv *}\cE \os{\sim}{\lra} 
Rf_{X'/\cY',\conv *} \varphi^*\cE. $$
\end{cor}

Now we introduce the notion of `having log smooth parameter': 

\begin{defn}\label{parameter}
We say that a proper log smooth integral 
morphism $f: (X,M_X) \lra \allowbreak (Y, \allowbreak 
M_Y)$ of fine log $B$-schemes has 
log smooth parameter in strong sense $($over $(B,M_B))$, 
if there exists a diagram of 
fine log formal $B$-schemes
\begin{equation}\label{star}
\begin{CD}
(X,M_X) @<<< (X',M_{X'}) @>>> (X_0,M_{X_0}) \\
@VfVV @V{f'}VV @V{f_0}VV \\
(Y,M_Y) @<g<< (Y',M_{Y'}) @>{g'}>> (Y_0,M_{Y_0}), 
\end{CD}
\end{equation}
where two squares are Cartesian, $g$ is strict etale and surjective, 
$f_0$ is proper log smooth integral and $(Y_0,M_{Y_0})$ is log smooth 
over $(B,M_B)$. We say a proper log smooth 
morphism $f: (X,M_X) \lra (Y,M_Y)$ of fine log $B$-schemes 
has log smooth parameter 
if we have a decomposition $X := \coprod_{i}X_i$ into open and closed 
subschemes such that the composite 
$$ (X_i,M_X|_{X_i}) \hra (X,M_X) \os{f}{\lra} (Y,M_Y)$$
has log smooth parameter in strong sense for each $i$.
\end{defn}

\begin{rem}
\begin{enumerate}
\item 
Let $f: (X,M_X) \lra (Y,M_Y)$ be a morphism of 
fine log $B$-schemes having log smooth parameter. Then any base change 
of $f$ in the category $\LB$ also has log smooth parameter. 
\item 
If Let $f: (X,M_X) \lra (Y,M_Y)$ is a morphism of 
fine log $B$-schemes having log smooth parameter in strong sense, 
we can take a diagram \eqref{star} with $Y_0$ affine. 
\end{enumerate}
\end{rem}

Then we have the following theorem, which is one of 
the main results in this section: 

\begin{thm}\label{perfmain}
Assume we are given a diagram 
\begin{equation}\label{perfdiag}
(X,M_X) \os{f}{\lra} (Y,M_Y) \os{\iota}{\hra} (\cY,M_{\cY}), 
\end{equation}
where $f$ is a proper log smooth integral morphism having log smooth parameter 
in $\LB$ and $\iota$ is a homeomorphic 
exact closed immersion in $\pLF$. 
Then, for a locally free isocrystal $\cE$ on 
$(X/S)^{\log}_{\conv}$, $Rf_{X/\cY,\conv *}\cE$ is a perfect complex 
of $\Q \otimes_{\Z} \cO_{\cY}$-modules. 
\end{thm}

\begin{lem}\label{perf3}
Let the notation be as in Theorem \ref{perfmain} and assume moreover that 
$f$ has log smooth parameter in strong sense. Then 
there exists a strict etale surjective morphism 
$\varphi: (\cY'',M_{\cY''}) \lra (\cY,M_{\cY})$ such that, 
if we denote the base change of the diagram \eqref{perfdiag} by $\varphi$ by 
$$ 
(X'',M_{X''}) \lra (Y'',M_{Y''}) \os{\iota''}{\hra} (\cY'',M_{\cY''}) $$
and if we put 
$(Y''_1, M_{Y''_1}) := (\cY'',M_{\cY''}) \otimes_{\Z_p} \Z/p\Z$, 
there exists a proper log smooth integral morphism 
$(X''_1, M_{X''_1}) \lra (Y''_1,M_{Y''_1})$ 
satisfying $(X''_1,M_{X''_1}) \times_{(Y''_1,M_{Y''_1})} 
(Y'',M_{Y''}) = (X'',M_{X''})$. 
\end{lem}

\begin{pf}
Let us take a diagram \eqref{star} with $Y_0$ affine 
and let $(\cY',M_{\cY'})$ be the 
fine log formal $\cB$-scheme which is strict formally 
etale over $(\cY,M_{\cY})$ 
satisfying $(\cY',M_{\cY'}) \times_{(\cY,M_{\cY})} (Y,M_Y) = (Y',M_{Y'})$. 
Then, since $Y_0$ is affine and $(Y_0,M_{Y_0})$ is log smooth over 
$(B,M_B)$, we can lift $(Y_0,M_{Y_0})$ to a fine log formal $\cB$-scheme 
$(\cY_0, M_{\cY_0})$ which is formally log smooth over $(\cB,M_{\cB})$. 
Then, by the infinitesimal lifting property of log smooth morphism, 
there exists 
a strict etale surjective morphism $\psi: (\cY'',M_{\cY''}) \lra 
(\cY',M_{\cY'})$ and a morphism $h:(\cY'',M_{\cY''}) \lra (\cY_0, M_{\cY_0})$ 
which fits into the commutative diagram 
\begin{equation*}
\begin{CD}
(Y'',M_{Y''}) @>{\iota''}>> (\cY'',M_{\cY''}) \\
@VVV @VhVV \\
(Y_0,M_{Y_0}) @>{\subset}>> (\cY_0, M_{\cY_0}), 
\end{CD}
\end{equation*}
where $\iota''$ denotes the base change of $\iota$ by the composite 
$$ (\cY'',M_{\cY''}) \os{\psi}{\lra} (\cY',M_{\cY'}) \lra (\cY,M_{\cY})$$
(which we denote by $\varphi$) and 
the left vertical arrow is the composite 
$$ 
(Y'',M_{Y''}) = (\cY'',M_{\cY''}) \times_{(\cY,M_{\cY})} (Y,M_Y) 
\os{\psi \times \id}{\lra} (\cY',M_{\cY'}) = (Y',M_{Y'}) \os{g'}{\lra} 
(Y_0,M_{Y_0}). $$
Then, if we denote the base change of 
$$ (X_0,M_{X_0}) \lra (Y_0,M_{Y_0}) \lra (\cY_0, M_{\cY_0}) $$
by $h$ by 
$$ (X''_1,M_{X''_1}) \lra (Y''_1,M_{Y''_1}) \lra (\cY'', M_{\cY''}), $$
we have $(Y''_1,M_{Y''_1}) = (\cY'',M_{\cY''}) \otimes_{\Z_p} \Z/p\Z$ and 
$(X''_1,M_{X''_1}) \times_{(Y''_1,M_{Y''_1})} 
(Y'',M_{Y''}) = (X'',M_{X''})$. So we are done. 
\end{pf}

\begin{pf*}{Proof of Theorem \ref{perfmain}}
We may assume that $X$ has log smooth parameter in strong sense. 
By Lemma \ref{perf3}, we have a strict formally etale \v{C}ech 
hypercovering $\epsilon: 
(\cY^{(\b)},M_{\cY^{(\b)}}) \lra (\cY,M_{\cY})$, such that, if we denote 
the base change of the diagram \eqref{perfdiag} by 
$(\cY^{(\b)},M_{\cY^{(\b)}}) \lra (\cY,M_{\cY})$ by 
$$ 
(X^{(\b)},M_{X^{(\b)}}) \lra 
(Y^{(\b)},M_{Y^{(\b)}}) \hra (\cY^{(\b)},M_{\cY^{(\b)}}) $$
and if we put $(Y^{(\b)}_1,M_{Y^{(\b)}_1}) := 
(\cY^{(\b)},M_{\cY^{(\b)}}) \otimes_{\Z_p} \Z/p\Z$, 
there exists a proper log smooth morphism 
$(X_1^{(n)},M_{X_1^{(n)}}) \lra (Y_1^{(n)},M_{Y_1^{(n)}})$ satisfying 
$(X_1^{(n)},M_{X_1^{(n)}}) \times_{(Y_1^{(n)},M_{Y_1^{(n)}})}
(Y^{(n)},\allowbreak M_{Y^{(n)}}) \allowbreak 
= (X^{(n)},M_{X^{(n)}})$ for each $n \in \N$. 
So, if we denote 
the restriction of $\cE$ to 
 $I_{\conv}((X^{(\b)}/\cY^{(\b)})^{\log})$ by $\cE^{(\b)}$, 
$Rf_{X^{(n)}/\cY^{(n)},\conv *}\cE^{(n)}$ is a perfect complex 
for each $n$, by 
Corollary \ref{perf2}. \par 
Now we claim that the canonical morphism 
$$Rf_{X/\cY,\conv *}\cE \lra 
R\epsilon_*Rf_{X^{(\b)}/\cY^{(\b)}, \conv *}\cE^{(\b)} $$
is a quasi-isomorphism. 
We can reduce the claim to 
the case where $X$ is affine. So we may assume that 
$(X,M_X)$ admits a closed immersion 
$(X,M_X) \hra (\cP,M_{\cP})$ to a fine log formal $\cB$-scheme 
$(\cP,M_{\cP})$ 
which is formally log smooth over $(\cY,M_{\cY})$ and satisfies 
$(\cP,M_{\cP}) \times_{(\cY,M_{\cY})} (Y,M_Y) 
= (X,M_X)$. In this case, we have the strict formally 
etale \v{C}ech hypercovering 
$\delta: (\cP^{(\b)},M_{\cP^{(\b)}}) \lra (\cP,M_{\cP})$ 
satisfying $(\cP^{(\b)},M_{\cP^{(\b)}}) \times_{(\cY,M_{\cY})} (Y,M_Y) 
= (X^{(\b)},M_{X^{(\b)}})$ since we have $\cP_{\et} \simeq X_{\et}$. 
If we denote the morphism 
$\cP \lra \cY$ by $h$, we have 
$$ Rf_{X/\cY,\conv *}\cE = Rh_* \sp_* \DR(\cP_K/\cY_K, \cE), $$
$$ R\epsilon_*Rf_{X^{(\b)}/\cY^{(\b)}, \conv *}\cE^{(\b)}
=  Rh_* R\delta_* \sp^{(\b)}_* \DR(\cP^{(\b)}_K/\cY^{(\b)}_K, 
\cE^{(\b)}), $$
where $\sp, \sp^{(\b)}$ is the spacialization map 
$\cP_K \lra \cP, \cP^{(\b)}_K \lra \cP^{(\b)}$, respectively. 
So, to prove the claim, it suffices to show 
that the map 
$$ 
\sp_* \DR(\cP_K/\cY_K, \cE) \lra
R\delta_* \sp^{(\b)}_* \DR(\cP^{(\b)}_K/\cY^{(\b)}_K 
\cE^{(\b)}) = R\sp_* R\delta_{K,*}\DR(\cP^{(\b)}_K/\cY^{(\b)}_K,\cE^{(\b)}) 
$$ 
is a quasi-isomorphism 
and it is reduced to 
the quasi-isomorphism 
$\cF \os{=}{\lra} R\delta_{K,*}\delta^*_K\cF$ for a coherent sheaf $\cF$ on 
$\cP_K$. This is true by \cite[7.1.2]{chts}. 
So we proved the claim. \par 
By the above claim, we have the spectral sequence 
$$ 
E_2^{s,t} = R^s\epsilon_*R^tf_{X^{(\b)}/\cY^{(\b)},\conv *}\cE^{(\b)} 
\,\Longrightarrow\, R^{s+t}f_{X/\cY,\conv *}\cE. $$
Note that 
$R^tf_{X^{(n)}/\cY^{(n)},\conv *}\cE^{(n)}$ is known to be isocoherent and 
that they are compatible with respect to $n$ by Corollary \ref{bc1}. 
So there exists (by etale descent of isocoherent sheaves) 
an isocoherent sheaf $\cF^t$ on $\cY$ such that 
$\cF^t \otimes_{\cO_{\cY}} \cO_{\cY^{(\b)}} = 
R^tf_{X^{(\b)}/\cY^{(\b)},\conv *}\cE^{(\b)}$ holds. 
Then we have 
$$ 
R^s\epsilon_*R^tf_{X^{(\b)}/\cY^{(\b)},\conv *}\cE^{(\b)} = 
\cF^t \,\, (s=0), \,\,\, 0 \,\, (s>0). $$
So we have $R^tf_{X/\cY,\conv *}\cE = \cF^t$ and it is isocoherent. 
Moreover, $Rf_{X/\cY,\conv *}\cE$ is bounded and it 
has finite tor-dimension Zariski locally 
because so does $Rf_{X^{(0)}/\cY^{(0)},\conv *}\cE^{(0)}
= Rf_{X^{(0)}/\cY^{(0)},\crys *}\allowbreak \Phi(\cE^{(0)})$. 
Therefore $Rf_{X/\cY,\conv *}\cE$ is a perfect complex
of $\Q \otimes_{\Z} \cO_{\cY}$-modules. 
\end{pf*}

\begin{rem}\label{remrem}
In the above proof, we have shown the isomorphism 
$$ R^qf_{X/\cY,\conv *}\cE \otimes_{\cO_{\cY}} \cO_{\cY^{(n)}} = 
R^qf_{X^{(n)}/\cY^{(n)},\conv *}\cE^{(n)}. $$
\end{rem}

Next we prove the base change property. 

\begin{thm}\label{bc2}
Assume we are given a diagram 
\begin{equation}\label{eq1eq1}
\begin{CD}
(X',M_{X'}) @>>> (Y',M_{Y'}) @>>> (\cY',M_{\cY'}) \\
@VVV @VVV @V{\varphi}VV \\
(X,M_X) @>f>> (Y,M_Y) @>>> (\cY,M_{\cY}), 
\end{CD}
\end{equation}
where 
the horizontal lines are as in Theorem \ref{perfmain} and the 
left square is Cartesian. 
Then, for a locally free isocrystal $\cE$ on $(X/\cY)^{\log}_{\conv}$, 
we have the quasi-isomorphism 
$$ 
L\varphi^*
Rf_{X/\cY,\conv *}\cE \os{\sim}{\lra} 
Rf_{X'/\cY',\conv *} \varphi^*\cE. $$
\end{thm}

\begin{pf}
We may assume that $f$ has log smooth parameter in strong sense and 
it suffices to prove the induced morphisms on cohomologies are isomorphisms 
etale locally. Let us take a strict formally etale surjective morphism 
$(\cY_0,M_{\cY_0}) \lra (\cY,M_{\cY})$ such that, 
if we denote the base change of the bottom horizontal line in \eqref{eq1eq1} 
by 
$(\cY_0,M_{\cY_0}) \lra (\cY,M_{\cY})$ by 
$$ 
(X_0,M_{X_0}) \lra (Y_0,M_{Y_0}) \hra (\cY_0,M_{\cY_0}), $$
there exists a proper 
log smooth morphism $(X_1,M_{X_1}) \lra (Y_1, M_{Y_1})$ 
saisfying $(X_1,M_{X_1}) \allowbreak 
\times_{(Y_1,M_{Y_1})} (Y,M_Y) = (X_0,M_{X_0})$, 
where $(Y_1, M_{Y_1}) := (\cY_0,M_{\cY_0}) \otimes_{\Z_p} \Z/p\Z$. 
Put $(\cY'_0, M_{\cY'_0}) := 
(\cY_0,M_{\cY_0}) \times_{(\cY,M_{\cY})} (\cY', M_{\cY'})$, 
denote the projection $(\cY'_0, M_{\cY'_0}) \lra (\cY', M_{\cY'})$ by 
$\varphi'$, denote the base change of the top horizontal line by  
$\varphi'$ by 
$$(X'_0,M_{X'_0}) \lra (Y'_0,M_{Y'_0}) \hra (\cY'_0,M_{\cY'_0}) $$
and denote the restriction of $\cE$ to $I_{\conv}((X_0/\cY_0)^{\log})$ by 
$\cE_0$. Then we have, by Remark \ref{remrem}, the isomorphisms 
$$ 
\cH^q(L\varphi^*Rf_{X/\cY,\conv *}\cE)|_{\cY'_0} = 
\cH^q(L{\varphi'}^*Rf_{X_0/\cY_0,\conv}\cE_0), $$
$$ 
\cH^q(Rf_{X'/\cY',\conv *} \varphi^*\cE)|_{\cY'_0} = 
\cH^q(Rf_{X'_0/\cY'_0, \conv *}{\varphi'}^*\cE_0).
$$ 
So it suffices to prove the quasi-isomorphism 
$$ L{\varphi'}^*Rf_{X_0/\cY_0,\conv}\cE_0 \os{\sim}{\lra} 
Rf_{X'_0/\cY'_0, \conv *}{\varphi'}^*\cE_0, 
$$ 
and it follows from Corollary \ref{bc1}. 
\end{pf}

Let us denote the category of $p$-adic fine log formal schemes over 
$(\cB,M_{\cB})$ of the form $\coprod_{\lam \in \Lam}(\cX_{\lam},M_{\lam})$ 
\,$((\cX_{\lam},M_{\lam}) \in \pLF)$ by $\pLF'$. For 
$(\cX,M_{\cX}) = \coprod_{\lam \in \Lam}(\cX_{\lam},M_{\lam}) \in \pLF'$ 
with $(\cX_{\lam},M_{\lam}) \in \pLF\, (\lam \in \Lam)$, we put 
$\Coh(\Q \otimes \cO_{\cX}) := \prod_{\lam \in \Lam} 
\Coh(\Q \otimes \cO_{\cX_{\lam}})$. 
A morphism $(\cY',M_{\cY'}) \lra (\cY,M_{\cY})$ 
in $\pLF'$ is said to be analytically flat if the induced functor 
$\Coh(\Q \otimes \cO_{\cY}) \lra \Coh(\Q \otimes \cO_{\cY'})$ is 
exact. Then, as a corollary of Corollary \ref{bc1} and Theorem 
\ref{bc2}, we have the 
following result. 

\begin{cor}\label{bc3}
Assume we are given a diagram in $\pLF$ 
\begin{equation*}
\begin{CD}
(X',M_{X'}) @>>> (Y',M_{Y'}) @>>> (\cY',M_{\cY'}) \\
@VVV @VVV @V{\varphi}VV \\
(X,M_X) @>f>> (Y,M_Y) @>{\iota}>> (\cY,M_{\cY}), 
\end{CD}
\end{equation*}
where $f$ is a proper log smooth integral morphism in 
$\LB$, $\iota$ is a 
homeomorphic exact closed immersion and the left square is 
Cartesian. 
Let $\cE$ be a locally free isocrystal on $(X/\cY)^{\log}_{\conv}$ and 
assume one of the conditions $(1)$, $(1)'$ below and 
one of the conditions $(2)$, $(2)'$ below are true$:$ \\
$(1)$ \,\, $\varphi$ is analytically flat. \\
$(1)'$ \,\, $R^qf_{X/\cY, \conv *}\cE$ is a locally free $\Q \otimes_{\Z}
\cO_{\cY}$-module for any $q$. \\
$(2)$ \,\, $f$ has log smooth parameter. \\
$(2)'$ \,\, If we put $(Y_1,M_{Y_1}) := (\cY,M_{\cY}) \otimes_{\Z_p} \Z/p\Z$, 
there exists a proper log smooth integral morphism 
$(X_1,M_{X_1}) \lra (Y_1,M_{Y_1})$ satisfying 
$(X_1,M_{X_1}) \times_{(Y_1,M_{Y_1})} (Y,M_Y) = (X,M_X)$. \\ 
Then we have the canonical isomorphism 
$$ \varphi^*R^qf_{X/\cY,\conv *}\cE \os{\cong}{\lra} 
R^qf_{X'/\cY',\conv *}\varphi^*\cE \quad (q \in \N). $$
\end{cor}

In the case where the condition (1) is satisfied, we call the result 
of Corollary \ref{bc3} `the analytically flat base change theorem'. \par 
In the rest of this section, we give some examples of 
analytically flat morphisms. Let us call a morphism 
$(\cY',M_{\cY'}) \lra (\cY,M_{\cY})$ in $\pLF'$ analytically faithfully flat 
if the induced functor $\Coh(\Q \otimes \cO_{\cY}) \lra \Coh(\Q \otimes 
\cO_{\cY'})$ is exact and faithful. First we prove a lemma, which 
provides the first examples of analytically (faithfully) flat morphisms:

\begin{lem}\label{anflex0}
\begin{enumerate}
\item 
A morphism $\varphi:(\cY',M_{\cY'}) \lra (\cY,M_{\cY})$ in $\pLF$ 
is analytically flat $($resp. analytically faithfully flat$)$ 
if $\cY' \lra \cY$ is flat $($resp. faithfully flat$)$. 
\item 
A morphism $\varphi: (\cY',M_{\cY'}) = \coprod_{\lam \in \Lam} 
(\cY'_{\lam},M_{\cY'_{\lam}}) \lra (\cY,M_{\cY})$ in 
$\pLF'$ is analytically flat if and only if each 
$\varphi|_{(\cY'_{\lam},M_{\cY'_{\lam}})}$ is analytically flat. 
\item 
Let us assume given the commutative diagram in $\pLF'$
\begin{equation*}
\begin{CD}
(\cZ',M_{\cZ'}) @>{\varphi'}>> (\cZ,M_{\cZ}) \\ 
@V{f'}VV @VfVV \\
(\cY',M_{\cY'}) @>{\varphi}>> (\cY,M_{\cY}) 
\end{CD}
\end{equation*}
and assume that $f$ is analytically flat and that $f'$ is 
analytically faithfully flat. Then, if $\varphi'$ is analytically flat, 
so is $\varphi$. 
\item 
A morphism in $\pLF'$ of the form 
$\varphi: (\cY',M_{\cY'}) = \coprod_{\lam \in \Lam} 
(\cY'_{\lam},M_{\cY'_{\lam}})\lra (\cY,M_{\cY}) \, 
((\cY'_{\lam},M_{\cY'_{\lam}}), (\cY,M_{\cY}) \in \pLF)$ 
is analytically flat 
$($resp. analytically faithfully flat$)$ 
if the morphisms 
$\varphi_K: \cY'_{\lam,K} \lra \cY_K \,(\lam \in \Lam)$ 
induced by $\varphi$ are 
admissible open immersions $($resp. if the morphisms 
$\varphi_K: \cY'_{\lam,K} \lra \cY_K \,(\lam \in \Lam)$ 
induced by $\varphi$ are 
admissible open immersions and if the morphism 
$\coprod_{\lam \in \Lam} \cY'_{\lam,K} \lra \cY_K$ is surjective$)$. 
\end{enumerate}
\end{lem}

\begin{pf} 
First we prove (1). 
If $\cY' \lra \cY$ is flat, $\Q \otimes_{\Z} \varphi^{-1}\cO_{\cY} \lra 
\Q \otimes_{\cZ} \cO_{\cY'}$ is flat and so $\varphi^* : \Coh(\Q 
\otimes \cO_{\cY}) \lra \Coh(\Q \otimes \cO_{\cY'})$ is exact. 
We prove $\varphi^*$ is faithful if $\varphi$ is faithfully flat. 
 To prove this, we may assume that $\cY, \cY'$ are in $\pLF$. 
Let us take a morphism $f:\Q \otimes M \lra \Q \otimes N$ such that 
$\varphi^*f: \Q \otimes \varphi^*M \lra \Q \otimes \varphi^*N$ is zero. 
It suffices to prove that $f$ is zero. To prove this, 
we may assume that $f$ comes from a morphism $M \lra N$. 
Then there exists a non-zero integer 
$n$ such that $n \varphi^*f:\varphi^*M \lra \varphi^*N$ is zero. Then 
we have $nf=0$ by the fully-faithfulness of $\varphi$ and so 
$f=0$ as the morphism $\Q \otimes M \lra \Q \otimes N$. \par 
(2) and (3) are immediate. (4) follows from the equivalence of categories 
$\Coh(\Q \otimes \cO_{\cY}) = \Coh(\cO_{\cY_K}), 
\Coh(\Q \otimes \cO_{\cY'_{\lam}}) = \Coh(\cO_{\cY'_{\lam,K}}).$ 
\end{pf} 

The following proposition gives an important example of 
analytically flat morphisms, which is useful in the next section: 

\begin{prop}\label{anflex}
Assume we are given a commutative diagram 
\begin{equation}\label{anfldiag}
\begin{CD} 
(X',M_{X'}) @>{\iota'}>> (\cP',M_{\cP'}) \\
@V{f}VV @V{g}VV \\ 
(X,M_{X}) @>{\iota}>> (\cP,M_{\cP}), 
\end{CD}
\end{equation}
where $f$ is a strict morphism in $\LB$, $\iota, \iota'$ are closed 
immersions and $g$ is a formally log smooth morphism in $\pLF$. 
Let 
$\{(T_{X,m},M_{T_{X,m}})\}_{m \in \N}, \{(T_{X',m},M_{T_{X',m}})\}_{m \in \N}$ 
be the system of universal enlargements 
of $((\cP,M_{\cP}),(X,M_{X})), 
((\cP',M_{\cP'}),(X',M_{X'}))$ respectively 
and let $g_m: (T_{X',m},M_{T_{X',m}}) \lra 
(T_{X,m},M_{T_{X,m}}) \, (m \in \N)$ be the morphisms induced by $g$. Then 
$g_m$'s are analytically flat. 
\end{prop}

\begin{pf} 
For any closed point $x'$ in $X'$, let us put $x:=f(x')$ and take 
a finite morphism of spectra of fields $\wt{x} \lra x$ such that, 
if we put $\wt{x}' := \wt{x} \times_x x'$, $\wt{x}_{\red}$ is 
isomorphic (over $\wt{x}$) to the disjoint union of finite number of 
$\wt{x}$'s. (This is possible because 
the morphism $x' \lra x$ corresponds to a finite extension of fields.) 
Also, take a strict morphism 
$(\wt{\cP},M_{\wt{\cP}}) \lra (\cP,M_{\cP})$ with $\wt{\cP} \lra \cP$ flat 
which fits into the following Cartesian diagram 
\begin{equation}\label{anfl-diag1}
\begin{CD}
(\wt{x},M_X|_{\wt{x}}) @>>> (\wt{\cP},M_{\wt{\cP}}) \\ 
@VVV @VVV \\ 
(x,M_X|_x) @>>> (\cP,M_{\cP}), 
\end{CD}
\end{equation}
where the lower horizontal arrow is the composite 
$(x,M_X|_x) \hra (X,M_X) \os{\iota}{\lra} (\cP,M_{\cP})$. 
(This is possible: Indeed, we may reduce to the case where $\wt{x} \lra x$ 
corresponds to a monogenic extension and it is easy in this case.) 
Then let us put $(\wt{\cP}',M_{\wt{\cP}}') := 
(\wt{\cP},M_{\wt{\cP}}) \times_{(\cP,M_{\cP})} (\cP',M_{\cP'})$. \par 
Let 
$\{(T_{x,m},M_{T_{x,m}})\}_{m \in \N}$, 
$\{(T_{x',m},M_{T_{x',m}})\}_{m \in \N}$, 
$\{(T_{\wt{x},m},M_{T_{\wt{x},m}})\}_{m \in \N}$, 
$\{(T_{\wt{x}',m},M_{T_{\wt{x}',m}})\}_{m \in \N}$, 
$\{(T_{\wt{x}'_{\red},m},M_{T_{\wt{x}'_{\red},m}})\}_{m \in \N}$ be 
the system of universal enlargements of $(x,M_X|_x)$, 
$(x',M_{X'}|_{x'})$, $(\wt{x},M_X|_{\wt{x}})$, $(\wt{x}',M_{X'}|_{\wt{x}'})$, 
$(\wt{x}'_{\red}, M_{X'}|_{\wt{x}'_{\red}})$ in 
$(\cP,M_{\cP})$, $(\cP',M_{\cP'})$, $(\wt{\cP},M_{\wt{\cP}})$, 
$(\wt{\cP}',M_{\wt{\cP}'})$, $(\wt{\cP}',M_{\wt{\cP}'})$, respectively. \par 
First let us consider the canonical morphism 
$$ \coprod_{n=m}^{\infty} \coprod_{x' \in X'} T_{x',n,K} \lra 
]X'[^{\log}_{\cP'}. $$
It is surjective and each component of this map is an admissible open 
immersion. (In this proof, 
we call such a morphism a surjective map by admissible open sets.) 
By pulling back this morphism by $T_{X',m,K} \os{\subset}{\lra} 
]X'[^{\log}_{\cP'}$, we obtain a surjective map by admissible open 
sets 
$$ \alpha'_K: \coprod_{n=m}^{\infty} \coprod_{x' \in X'} 
(T_{x',n} \times_{T_{X',n}} T_{X',m})_K \lra T_{X',m,K} $$  
which comes from the morphism in $\pLF'$ 
$$ \alpha': \coprod_{n=m}^{\infty} \coprod_{x' \in X'} 
T_{x',n} \times_{T_{X',n}} T_{X',m} \lra T_{X',m}. $$  
By the same argument, we obtain another surjective map by 
admissible open sets 
$$ \alpha_K : 
\coprod_{n=m}^{\infty} \coprod_{x' \in X'} (T_{x,n} \times_{T_{X,n}} 
T_{X,m})_K \lra T_{X,m,K} $$ 
which again comes from the morphism in $\pLF'$ 
$$ \alpha : 
\coprod_{n=m}^{\infty} \coprod_{x' \in X'} T_{x,n} \times_{T_{X,n}} 
T_{X,m} \lra T_{X,m}. $$ 
It is easy to see that the morphisms $\alpha, \alpha'$ 
fit into the following commutative diagram: 
\begin{equation*}
\begin{CD}
\coprod_{n=m}^{\infty} \coprod_{x' \in X'} 
T_{x',n} \times_{T_{X',n}} T_{X',m} @>>> 
\coprod_{n=m}^{\infty} \coprod_{x' \in X'} T_{x,n} \times_{T_{X,n}} 
T_{X,m} \\ 
@V{\alpha'}VV @V{\alpha}VV \\ 
T_{X',m} @>{g_{m}}>> T_{X,m}. 
\end{CD}
\end{equation*}
By Lemma \ref{anflex0} (4), $\alpha, \alpha'$ are analytically 
faithfully flat, and by Lemma \ref{anflex0} (1), (2), (3) and the above 
diagram, $g_m$ is analytically flat if the map 
$T_{x',n} \times_{T_{X',n}} T_{X',m} \lra T_{x,n} \times_{T_{X,n}} 
T_{X,m}$ (which we denote by $\beta$)
is analytically flat for any $n \geq m$ and $x'$. So we are reduced to 
showing that the morphism $\beta$ is analytically flat. \par 
Next, let us note that the morphisms 
$$ T_{\wt{x}',n} \lra T_{x',n}, \quad T_{\ti{x},n} \lra T_{x,n} $$ 
induced by \eqref{anfl-diag1} is faithfully flat: Indeed, since 
$T_{\wt{x}',n}, T_{x',n}, T_{\ti{x},n}, T_{x,n}$ are unchanged if 
we shrink $\cP$ and $\wt{\cP}$ as long as $\cP$ (resp. $\wt{\cP}$) 
contains $x$ (resp. $\wt{x}$), we may assume that $\wt{\cP} \lra 
\cP$ is faithfully flat to prove this claim. 
In this case, the claim follows from the isomorphisms 
$T_{\wt{x}',n} = T_{x',n} \times_{\cP} \ti{\cP}, 
T_{\ti{x},n} = T_{x,n} \times_{\cP} \ti{\cP}$, which is true by 
Lemma \ref{2.1.26}. Hence the morphisms 
$$ \gamma': T_{\wt{x}',n} \times_{T_{X',n}} T_{X',m} \lra 
T_{x',n} \times_{T_{X',n}} T_{X',m}, \,\,\,\, 
\gamma: T_{\wt{x},n} \times_{T_{X,n}} T_{X,m} \lra T_{x,n} 
\times_{T_{X,n}} T_{X,m}$$
are faithfully flat and they fit into the following commutative diagram: 
\begin{equation*}
\begin{CD}
T_{\wt{x}',n} \times_{T_{X',n}} T_{X',m} @>>> 
T_{\wt{x},n} \times_{T_{X,n}} T_{X,m} \\ 
@V{\gamma'}VV @V{\gamma}VV \\ 
T_{x',n} \times_{T_{X',n}} T_{X',m} @>{\beta}>> 
T_{x,n} \times_{T_{X,n}} T_{X,m}. 
\end{CD}
\end{equation*}
So, by Lemma \ref{anflex0} (3), we are reduced to prove that the morphism 
$T_{\wt{x}',n} \times_{T_{X',n}} T_{X',m} \lra 
T_{\wt{x},n} \times_{T_{X,n}} T_{X,m}$ (which we denote by $\delta$) 
is analytically flat. \par 
Next let us consider the morphism 
$$ \coprod_{l=n}^{\infty} (T_{\wt{x}'_{\red},l})_K \lra 
]\wt{x}'[^{\log}_{\wt{\cP}'}. $$ 
It is a surjective map by admissible open sets since we have 
$]\wt{x}'_{\red}[^{\log}_{\wt{\cP}'}= ]\wt{x}'[^{\log}_{\wt{\cP}'}$. 
By pulling back this morphism by $(T_{\wt{x},n} \times_{T_{X',n}} 
T_{X',m})_K \os{\subset}{\lra} ]\wt{x}'[^{\log}_{\wt{\cP}'}$, 
we obtain a surjective map by admissible open sets
$$ \delta'_K: \coprod_{l=n}^{\infty} 
(T_{\wt{x}'_{\red},l} \times_{T_{\wt{x}',l}} T_{\wt{x}',n} 
\times_{T_{X',n}} T_{X',m})_K \lra (T_{\wt{x},n} \times_{T_{X',n}} 
T_{X',m})_K$$ 
which comes from the morphism in $\pLF'$ 
$$ \delta': \coprod_{l=n}^{\infty} 
T_{\wt{x}'_{\red},l} \times_{T_{\wt{x}',l}} T_{\wt{x}',n} 
\times_{T_{X',n}} T_{X',m} \lra T_{\wt{x},n} \times_{T_{X',n}} 
T_{X',m}. $$ 
So, again by Lemma \ref{anflex0}, we are reduced to showing that, 
for any $l \geq n \geq m$, 
the composite morphism 
$$ 
T_{\wt{x}'_{\red},l} \times_{T_{\wt{x}',l}} T_{\wt{x}',n} 
\times_{T_{X',n}} T_{X',m} 
\os{\delta'}{\lra} T_{\wt{x},n} \times_{T_{X',n}} T_{X',m} 
\os{\delta}{\lra} T_{\wt{x},n} \times_{T_{X,n}} T_{X,m}, 
$$ 
which we denote by $\epsilon$, is analytically flat. \par 
Now let us put $S := T_{\wt{x},n} \times_{T_{X,n}} T_{X,m}$. 
Then the morphism $\epsilon: 
T_{\wt{x}'_{\red},l} \times_{T_{\wt{x}',l}} T_{\wt{x}',n} 
\times_{T_{X',n}} T_{X',m} \lra S$ is factorized as follows: 
\begin{align*}
& T_{\wt{x}'_{\red},l} \times_{T_{\wt{x}',l}} T_{\wt{x}',n} 
\times_{T_{X',n}} T_{X',m} \\ 
& \,\lra\, 
T_{\wt{x}'_{\red},l} \times_{T_{\wt{x}',l}} 
(T_{\wt{x}',l} \times_{T_{\wt{x},l}} T_{\wt{x},n}) 
\times_{T_{X',n}} (T_{X',n} \times_{T_{X,n}} T_{X,m}) \\ 
& \,=\,  
T_{\wt{x}'_{\red},l} \times_{T_{\wt{x},l}} S \lra S. 
\end{align*}
Note that the first arrow is analytically flat, since the associated morphism 
between rigid analytic spaces is an admissible open immersion. Hence we are 
reduced to showing that the projection 
$T_{\wt{x}'_{\red},l} \times_{T_{\wt{x},l}} S \lra S$ (which we denote by 
$\pi$) is analytically flat. Now let us consider 
the following commutative diagram 
\begin{equation*}
\begin{CD}
(T_{\wt{x}'_{\red},l} \times_{T_{\wt{x},l}} S)_K @>{\subset}>> 
]\wt{x}'_{\red}[^{\log}_{\wt{\cP}'} \\ 
@V{\pi_K}VV @V{\varphi}VV \\
S_K @>{\subset}>> ]\wt{x}[^{\log}_{\wt{\cP}}, 
\end{CD}
\end{equation*}
where $\varphi$ is the morphism induced by $\wt{x}'_{\red} \lra \wt{x}' \lra 
\wt{x}$ and the horizontal arrows are natural admissible open immersions. 
Since $\wt{x}'_{\red}$ is isomorphic to the disjoint union 
of finite number of $\wt{x}$'s 
and $(\wt{\cP}',M_{\wt{\cP}'}) \lra (\wt{\cP},M_{\wt{\cP}})$ 
is formally log smooth, there exists an isomorphism 
$]\wt{x}'_{\red}[^{\log}_{\wt{\cP}'} \cong 
\coprod_{i=1}^n ]\wt{x}[^{\log}_{\wt{\cP}} \times D^{r_i}$ 
for some $n, r_i \in \N$ such that $\varphi$ is identified with 
the projection 
$$ \coprod_{i=1}^n ]\wt{x}[^{\log}_{\wt{\cP}} \times D^{r_i} \lra 
\coprod_{i=1}^n ]\wt{x}[^{\log}_{\wt{\cP}} \lra 
]\wt{x}[^{\log}_{\wt{\cP}} $$ 
via this isomorphism. Then the morphism $\pi_K$ factors as 
$$ 
(T_{\wt{x}'_{\red},l} \times_{T_{\wt{x},l}} S)_K
\os{\subset}{\lra} \varphi^{-1}(S_K) = 
\coprod_{i=1}^n S_K \times D^{r_i} \lra S_K, 
$$ 
where the first map is an admissible open immersion and the second 
map is the projection. So it is easy to see that the functor 
$$ \Coh(\Q \otimes \cO_S) = \Coh(\cO_{S_K}) 
\os{\epsilon_K^*}{\lra} 
\Coh(\cO_{(T_{\wt{x}'_{\red},l} \times_{T_{\wt{x},l}} S)_K})
= 
\Coh(\Q \otimes \cO_{T_{\wt{x}'_{\red},l} \times_{T_{\wt{x},l}} S}), 
$$ 
which is the same as $\epsilon^*$, is exact. So we have proved that 
$\epsilon$ is analytically flat and so we are done. 
\end{pf}

\section{Relative log analytic cohomology}

In this section, first we introduce a rigid analytic variant of 
relative log convergent cohomology, which we call relative log analytic 
cohomology. It is regarded as a relative version of analytic cohomology 
introduced in \cite{shiho2}. 
Then we will prove 
a relation between relative log convergent cohomology and 
relative log analytic cohomology for proper log smooth integral morphisms 
having log smooth parameter. This implies the coherence of 
relative log analytic cohomology. After that, we 
prove the existence of a canonical structure of an isocrystal on 
relative log analytic cohomology. \par 
First we give a definition of relative log analytic 
cohomology. 

\begin{defn}\label{defrellogancoh}
Assume we are given a diagram 
\begin{equation}\label{diag3-1}
(X,M_X) \os{f}{\lra} (Y,M_Y) \os{\iota}{\hra} (\cY,M_{\cY}), 
\end{equation}
where $f$ is a morphism in $\LB$ and $\iota$ is a closed immersion 
in $\pLF$, 
and let $\cE$ be an isocrystal on $(X/\cY)^{\log}_{\conv}$. 
Take an embedding system over $\cY$ 
\begin{equation}\label{diag3-2}
 (X,M_{X}) \overset{g}{\lla} (X^{(\bullet)}, M_{X^{(\bullet)}}) 
\overset{i}{\hra} (\cP^{(\bullet)},M_{\cP^{(\bullet)}}), 
\end{equation}
let $\cE^{(\b)}$ be the restriction of $\cE$ to 
$(X^{(\b)},M_{X^{(\b)}})$, 
let $\DR(]X^{(\b)}[^{\log}_{\cP^{(\b)}}/\cY_K,\cE^{(\b)})$ be the 
log de Rham complex associated to $\cE^{(\b)}$ and let $h$ be the morphism 
$]X^{(\b)}[^{\log}_{\cP^{(\b)}} \lra ]Y[^{\log}_{\cY}$ induced by the 
embedding system \eqref{diag3-2}. 
Then we define $Rf_{X/\cY, \an *}\cE$, 
$R^qf_{X/\cY,\an *}\cE$ by 
$$ 
Rf_{X/\cY,\an *}\cE := 
Rh_*\DR(]X^{(\b)}[^{\log}_{\cP^{(\b)}}/\cY_K,\cE^{(\b)}), \,\,\,\, 
R^qf_{X/\cY,\an *}\cE := 
R^qh_*\DR(]X^{(\b)}[^{\log}_{\cP^{(\b)}}/\cY_K,\cE^{(\b)}) $$
and we call $R^qf_{X/\cY,\an *}\cE$ the $q$-th relative log analytic 
cohomology of $(X,M_X)/(\cY,M_{\cY})$ with coefficient $\cE$. 
It is a sheaf of $\cO_{]Y[^{\log}_{\cY}}$-modules. 
\end{defn}

\begin{rem}\label{rem3-1}
With the above notation, let 
$\{(\cY_m,M_{\cY_m})\}_m$ be the system of universal enlargements of 
$((\cY,M_{\cY}), (Y,M_Y), \iota,\id)$ and denote the base change of 
the diagrams \eqref{diag3-1}, \eqref{diag3-2} by 
$(\cY_m,M_{\cY_m}) \lra (\cY,M_{\cY})$ by 
\begin{equation}
(X_m,M_{X_m}) \lra (Y_m,M_{Y_m}) \hra (\cY_m,M_{\cY_m}), 
\end{equation}
\begin{equation}
 (X_m,M_{X_m}) \lla (X_m^{(\bullet)}, M_{X_m^{(\bullet)}}) 
\hra (\cP_m^{(\bullet)},M_{\cP_m^{(\bullet)}}), 
\end{equation}
respectively. Let $\cE_m$, $\cE_m^{(\b)}$ be the restriction of $\cE$ to 
$(X_m,M_{X_m})$, 
$(X_m^{(\b)},\allowbreak M_{X_m^{(\b)}})$ respectively, 
let $\DR(]X_m^{(\b)}[^{\log}_{\cP_m^{(\b)}}/\cY_{m,K},\cE_m^{(\b)})$ be the 
log de Rham complex associated to $\cE_m^{(\b)}$ and let $h_m$ be the morphism 
$]X_m^{(\b)}[^{\log}_{\cP_m^{(\b)}} \lra ]Y[^{\log}_{\cY_m} = \cY_{m,K}$. 
Then we have $]Y[^{\log}_{\cY}= \bigcup_m  \cY_{m,K}, 
h_m = h \vert_{]X^{(\b)}_m[^{\log}_{\cP_m^{(\b)}}}, \,
]X^{(\b)}_m[^{\log}_{\cP^{(\b)}_m} = h_m^{-1}(\cY_{m,K})$ and 
$$ R^qf_{X_m/\cY_m,\an *}\cE_m = 
R^qh_{m *}\DR(]X_m^{(\b)}[^{\log}_{\cP_m^{(\b)}}/\cY_{m,K},\cE_m^{(\b)}) $$
is nothing but the restriction of $R^qf_{X/\cY,\an *}\cE$ to 
$\cY_{m,K}$. 
\end{rem}

In order to assure that Definition \ref{defrellogancoh} is well-defined, 
we should prove the following proposition: 

\begin{prop}\label{wdprop}
Let the notations be as in Definition \ref{defrellogancoh}. Then 
the definition of the relative log analytic cohomology 
$R^qf_{X/\cY,\an *}\cE$ of $(X,M_X)/(\cY,M_{\cY})$ with coefficient $\cE$ 
is independent of the choice of the embedding system. 
\end{prop}

The method of proof is similar to that of \cite[2.2.14]{shiho2}. 
First we prove the following descent property: 

\begin{lem}\label{wdlem1}
Assume we are given the Cartesian diagram 
\begin{equation}\label{wd1}
\begin{CD}
(X^{(\b)},M_{X^{(\b)}}) @>>> (\cP^{(\b)},M_{\cP^{(\b)}}) \\ 
@VVV @VgVV \\
(X,M_X) @>{\iota}>> (\cP,M_{\cP}), 
\end{CD}
\end{equation}
where $(X,M_X)$ is an object in $\L$, $(\cP,M_{\cP})$ is an object in 
$\LF$, $\iota$ is a closed immersion 
and $g$ is a strict formally etale hypercovering. Let 
$g_K:]X^{(\b)}[^{\log}_{\cP^{(\b)}} \lra ]X[^{\log}_{\cP}$ be the 
morphism induced by $g$ and let $\cE$ be a coherent sheaf on 
$]X[^{\log}_{\cP}$. Then we have the quasi-isomorphism 
$$ \cE \overset{=}{\lra} Rg_{K,*}g_K^*\cE. $$
\end{lem} 

\begin{pf} 
Let $\{(\cP_m,M_{\cP_m})\}_m$ be the system of universal enlargements 
of the pre-widening $((\cP,M_{\cP}),(X,M_X),\iota,i)$ and let 
\begin{equation*}
\begin{CD}
(X_m^{(\b)},M_{X_m^{(\b)}}) @>>> (\cP_m^{(\b)},M_{\cP_m^{(\b)}}) \\ 
@VVV @V{g_m}VV \\
(X_m,M_{X_m}) @>{\iota}>> (\cP_m,M_{\cP_m}) 
\end{CD}
\end{equation*}
be the diagram obtained by applying $\times_{(\cP,M_{\cP})} (\cP_m,M_{\cP_m})$ 
to the diagram \eqref{wd1}. Let $g_{m,K}:P^{(\b)}_{m,K} \lra P_{m,K}$ be the 
map induced by $g_m$. Then we have 
$$ ]X[^{\log}_{\cP} = \bigcup_m \cP_{m,K}, \,\, 
]X^{(\b)}[^{\log}_{\cP^{(\b)}} = \bigcup_m \cP^{(\b)}_{m,K}, \,\, 
g_K^{-1}(\cP_{m,K}) = \cP^{(\b)}_{m,K}, \,\, 
g_{m,K} = g_K \vert_{\cP^{(\b)}_{m,K}}.$$ 
So we have 
$(Rg_{K,*}g_K^*\cE)\vert_{\cP_{m,K}} = 
Rg_{m,K,*}g^*_{m,K}(\cE \vert_{\cP_{m,K}})$ and so it suffices to prove the 
quasi-isomorphism $\cE \vert_{\cP_{m,K}} \os{=}{\lra} 
Rg_{m,K,*}g^*_{m,K}(\cE \vert_{\cP_{m,K}})$. 
This is already proven in \cite[7.3.3]{chts}. So we are done. 
\end{pf} 

Next we prove a lemma which corresponds to \cite[2.2.15]{shiho2}: 

\begin{lem}\label{logsminv}
Let $(X,M_X) \lra (Y,M_Y) \hra (\cY,M_{\cY})$ as in Definition 
\ref{defrellogancoh} and assume we are given a commutative diagram 
over $(\cY,M_{\cY})$ 
\begin{equation}\label{wd5}
\begin{CD}
(X,M_X) @>{\iota_1}>> (\cP_1,M_{\cP_1}) \\ 
@\vert @V{\varphi}VV \\ 
(X,M_X) @>{\iota_2}>> (\cP_2,M_{\cP_2}), 
\end{CD}
\end{equation}
where $(\cP_j,M_{\cP_j})$ are objects in $\pLF$ which are formally log 
smooth over $(\cY,M_{\cY})$, 
$\iota_j$ are closed immersions $(j=1,2)$ and $\varphi$ is a formally 
log smooth morphism. Let 
$\varphi_K: ]X[^{\log}_{\cP_1} \lra ]X[^{\log}_{\cP_2}$ be the morphism 
induced by $\varphi$ and let $\cE$ be 
an isocrystal on $(X/\cY)^{\log}_{\conv}$. Then we have a quasi-isomorphism 
$$ \DR(]X[^{\log}_{\cP_2}/\cY_K,\cE) \os{=}{\lra} 
R\varphi_{K,*}\DR(]X[^{\log}_{\cP_1}/\cY_K,\cE). $$
\end{lem}

\begin{pf} 
For $j=1,2$, let $\cP^{\ex}_j$ be the exactification of the closed immersion 
$\iota_j$. Then it suffices to prove the lemma Zariski locally on 
$\cP^{\ex}_2$. So we may assume that $\cP^{\ex}_{1,K}$ is isomorphic to 
$\cP^{\ex}_{2,K} \times D^r$ and the morphism $\varphi_K$ is equal to the 
projection $\cP^{\ex}_{2,K} \times D^r \lra \cP^{\ex}_{2,K}$. 
Let $(t_1, \cdots, t_r)$ be the coordinate of $D^r$ and let 
$\Omega^{\b}$ be the complex 
$$ 
[\cO_{\cP^{\ex}_{1,K}} \lra 
\bigoplus_{i=1}^r\cO_{\cP^{\ex}_{1,K}} dt_i \lra 
\bigoplus_{1 \leq i_1 < i_2 \leq r} \cO_{\cP^{\ex}_{1,K}} dt_{i_1} \wedge 
dt_{i_2} \lra \cdots \lra 
\cO_{\cP^{\ex}_{1,K}} dt_1 \wedge \cdots \wedge dt_r]. $$ 
Then $\DR(]X[^{\log}_{\cP_1}/\cY,\cE)$ is equal to the total complex 
associated to the double complex 
$\varphi_K^*\DR(]X[^{\log}_{\cP_2}/\cY,\cE) 
\otimes_{\cO_{\cP^{\ex}_{1,K}}} \Omega^{\b}$. 
So, to prove the lemma, it suffices to prove the quasi-isomorphism 
$$ E \os{\simeq}{\lra} R\varphi_{K,*}(\varphi_K^*E 
\otimes_{\cO_{\cP^{\ex}_{1,K}}} \Omega^{\b})
= \varphi_{K,*}(\varphi_K^*E \otimes_{\cO_{\cP^{\ex}_{1,K}}} \Omega^{\b}) = 
E \otimes_{\cO_{\cP^{\ex}_{2,K}}} \varphi_{K,*}\Omega^{\b} $$ 
for a coherent $\cO_{\cP^{\ex}_{2,K}}$-module $E$. This is reduced to 
showing that the complex $\cO_{\cP^{\ex}_{2,K}} \lra \varphi_{K,*}\Omega^{\b}$
 is homotopic to zero. We can prove it in the same way as the proof of 
 Lemma \ref{add2}. 
\end{pf} 

\begin{pf*}{Proof of Proposition \ref{wdprop}}
The proof is similar to \cite[2.2.14]{shiho2}. \par 
Let 
$$ (X,M_X) \overset{g_j}{\lla} (X^{(\bullet)}_j,M_{X^{(\bullet)}_j}) 
\overset{\iota^{(\bullet)}_j}{\hra} (\cP^{(\bullet)}_j,M_{\cP^{(\bullet)}_j})
\quad 
(j=1,2) $$ 
be embedding systems and denote the pull-back of 
$\cE$ to $(X_j^{(\bullet)}/\cY)^{\log}_{\conv}$ by 
$\cE^{(\bullet)}_j \, (j=1,2)$. 
Put $(X^{(m,n)}, M_{X^{(m,n)}}) := (X^{(m)}_1, M_{X^{(m)}_1}) 
\times_{(X,M_X)} (X^{(n)}_2, M_{X^{(n)}_2})$ and 
$(\cP^{(m,n)}, M_{\cP^{(m,n)}}) \allowbreak := \allowbreak 
(\cP^{(m)}_1, M_{\cP^{(m)}_1}) 
\times_{(\cY,M_{\cY})} (\cP^{(n)}_2, M_{\cP^{(n)}_2})$. 
(
$(X^{(\bullet,\bullet)}, M_{X^{(\bullet,\bullet)}})$ 
forms a bisimplicial fine log $B$-scheme and 
$(\cP^{(\bullet,\bullet)}, M_{\cP^{(\bullet,\bullet)}})$ 
 forms a bisimplicial $p$-adic fine log formal $\cB$-scheme.) 
Let $\varphi_{j,K}:]X^{(\b,\b)}[^{\log}_{\cP^{(\b,\b)}} \lra 
]X_j^{(\b)}[^{\log}_{\cP_j^{(\b)}}$ be the morphism induced by 
$(\cP^{(\b,\b)},M_{\cP^{(\b,\b)}}) \lra (\cP_j^{(\b)},M_{\cP_j^{(\b)}})$ and 
denote the pull-back of $\cE$ to 
$(X^{(\bullet,\bullet)}/\cY)^{\log}_{\conv}$ by 
$\cE^{(\bullet,\bullet)}$. Then, 
to prove the proposition, it suffices to show the quasi-isomorphisms 
$$ 
   \DR(]X^{(\b)}_j[^{\log}_{\cP^{(\b)}_j}/\cY_K, \cE^{(\b)}_j) \isom 
   R\varphi_{j,K,*} 
   \DR(]X^{(\b,\b)}[^{\log}_{\cP^{(\b,\b)}}/\cY_K, 
   \cE^{(\b,\b)}) \quad (j=1,2), $$ 
and it suffice to treat the case $j=1$. 
Let us fix $n \in \N$ and let 
$\varphi^{(n)}_K: ]X^{(n,\b)}[^{\log}_{\cP^{(n,\b)}} \lra 
]X^{(n)}[^{\log}_{\cP_1^{(n)}}$ be the morphism 
induced by the first projection 
$(\cP^{(n,\b)},M_{\cP^{(n,\b)}}) \lra (\cP_1^{(n)},M_{\cP^{(n)}_1})$. 
Then, to prove the above quasi-isomorphism (for $j=1$), 
it suffices to prove the quasi-isomorphism 
\begin{equation}\label{adada}
\DR(]X^{(n)}_1[^{\log}_{\cP^{(n)}_1}/\cY_K, \cE^{(n)}_1) \isom 
R\varphi^{(n)}_{K,*}
\DR(]X^{(n,\b)}[^{\log}_{\cP^{(n,\b)}}/\cY_K, \cE^{(n,\b)}). 
\end{equation}

In the following, we give a proof of the quasi-isomorphism 
(\ref{adada}). 
Let $\hat{\cP}^{(n)}_1$ 
be the formal completion of $\cP^{(n)}_1$ 
along $X^{(n)}_1$, put $M_{\hat{\cP}^{(n)}_1} := 
M_{\cP^{(n)}_1}|_{\hat{\cP}^{(n)}_1}$ and 
let 
$h_n: (\hat{\cP}^{(n,\b)}_1,M_{\hat{\cP}^{(n,\b)}})
\lra (\hat{\cP}^{(n)}_1,M_{\hat{\cP}^{(n)}_1})$ be the unique 
strict formally 
etale hypercovering satisfying 
$$(\hat{\cP}^{(n,\b)},M_{\hat{\cP}^{(n,\b)}}) 
\times_{(\hat{\cP}^{(n)}_1,\allowbreak M_{\hat{\cP}^{(n)}_1})} (X^{(n)}_1,M_{X^{(n)}_1}) 
\isom (X^{(n,\b)},M_{X^{(n,\b)}}).$$ 
Let $(\ti{\cP}^{(n,m)},M_{\ti{\cP}^{(n,m)}})$ 
be $(\cP^{(n,\b)},M_{\cP^{(n,\b)}}) \times_{(\cY,M_{\cY})} 
(\hat{\cP}^{(n,\b)}_1,M_{\hat{\cP}^{(n,\b)}_1})$. Then we have the 
following diagram: 
\begin{equation*}
\begin{CD}
(X_1^{(n)},M_{X_1^{(n)}}) @<{g_n}<<
(X^{(n,\b)},M_{X^{(n,\b)}}) @= (X^{(n,\b)},M_{X^{(n,\b)}}) @= 
(X^{(n,\b)},M_{X^{(n,\b)}}) \\ 
@VVV @VVV @VVV @VVV \\ 
(\hat{\cP}^{(n)}_1,M_{\hat{\cP}^{(n)}_1}) @<{h_n}<< 
(\hat{\cP}^{(n,\b)}_1, M_{\hat{\cP}^{(n,\b)}_1}) @<{\pr_1^{(\b)}}<< 
(\ti{\cP}^{(n,\b)},M_{\ti{\cP}^{(n,\b)}})  
@>{\pr_2^{(\b)}}>> ({\cP}^{(n,\b)},M_{{\cP}^{(n,\b)}}), 
\end{CD}
\end{equation*}
where the vertical arrows are the canonical 
closed immersions. Then, by the claim shown in \cite[p.81]{shiho2}, 
Zariski locally on $X^{(n,m)}$ 
there exists an exact closed immersion 
$$ (X^{(n,m)},M_{X^{(n,m)}}) \hra (\cP^{(n,m)}_1, M_{\cP^{(n,m)}_1})$$ 
of $(X^{(n,m)},M_{X^{(n,m)}})$ into an object 
in $\pLF$ which is strict formally etale 
over $(\cP^{(n)}_1, \allowbreak M_{\cP^{(n)}_1})$ such that 
$(\hat{\cP}^{(n,m)}_1, M_{\hat{\cP}^{(n,m)}_1})$ is nothing but 
the completion of 
$(\cP^{(n,m)}_1, M_{\cP^{(n,m)}_1})$ along $(X^{(n,m)}, \allowbreak 
M_{X^{(n,m)}})$. 
Then we have the canonical isomorphism $\cP^{(n,m),\ex}_1 \cong 
\hat{\cP}^{(n,m),\ex}_1$ of exactifications. On the other hand, 
if we put $(\ol{\cP}^{(n,m)},M_{\ol{\cP}^{(n,m)}}) := 
(\cP^{(n,m)}_1,M_{\cP^{(n,m)}_1}) 
\times_{(\cY,M_{\cY})} (\cP^{(n,m)}, \allowbreak M_{\cP^{(n,m)}})$, 
we have the 
canonical isomorphism 
$\ol{\cP}^{(n,m),\ex} \cong \ti{\cP}^{(n,m),\ex}$. So, by 
Remark \ref{generalized-dr}, we can define log de Rham complexes 
$$ 
\DR(]X^{(n,\b)}[^{\log}_{\hat{\cP}^{(n,\b)}_1}/\cY_K, \cE^{(n,\b)}) =: 
\wh{\DR}^{(n,\b)}_1, 
\qquad 
\DR(]X^{(n,\b)}[^{\log}_{\ti{\cP}^{(n,\b)}}/\cY_K, \cE^{(n,\b)}) =: 
\ti{\DR}^{(n,\b)}. $$ 
Put $\DR^{(n,\b)}:= 
\DR(]X^{(n,\b)}[^{\log}_{\cP^{(n,\b)}}/\cY_K,\cE^{(n,\b)}), 
\DR^{(n)}_1 := \DR(]X_1^{(n)}[^{\log}_{\cP^{(n)}_1}/\cY_K, \cE^{(n)})$. 
Then the proof of the 
isomorphism (\ref{adada}) is reduced to the three quasi-isomorphisms 
$$ R \pr^{(\b)}_{1,K,*} \wt{\DR}^{(n,\b)} = \wh{\DR}_1^{(n,\b)}, 
\,\,\,\, R \pr^{(\b)}_{2,K,*} \wt{\DR}^{(n,\b)} = \wh{\DR}^{(n,\b)}, 
\,\,\,\, R h_{n,K,*} \wh{\DR}_1^{(n,\b)} = \DR^{(n)}_1. $$
The first and the second quasi-isomorphisms follow from 
Lemma \ref{logsminv} and the fact that the morphisms 
$$ \pr^{(m)}_1: (\ti{\cP}^{(n,m)}, M_{\ti{\cP}^{(n,m)}}) 
\lra (\hat{\cP}^{(n,m)}_1, M_{\hat{\cP}^{(n,m)}_1}), $$
$$ \pr^{(m)}_2: (\ti{\cP}^{(n,m)}, M_{\ti{\cP}^{(n,m)}}) 
\lra (\hat{\cP}^{(n,m)}, M_{\hat{\cP}^{(n,m)}}) $$
are locally the completions of the formally log smooth morphisms 
$$ (\ol{\cP}^{(n,m)}, M_{\ol{\cP}^{(n,m)}}) \lra  
(\cP^{(n,m)}_1, M_{\cP^{(n,m)}_1}), $$  
$$(\ol{\cP}^{(n,m)}, M_{\ol{\cP}^{(n,m)}}) \lra  
(\cP^{(n,m)}, M_{\cP^{(n,m)}}), $$
respectively. 
The third quasi-isomorphism follows from Lemma \ref{wdlem1}. So we are done. 
\end{pf*}

The next theorem establishes the relation of relative log convergent 
cohomology and relative log analytic cohomology: 

\begin{thm}\label{coherence0}
Assume we are given a diagram 
\begin{equation}\label{diag-coherence0}
(X,M_X) \os{f}{\lra} (Y,M_Y) \os{\iota}{\hra} (\cY,M_{\cY}), 
\end{equation}
where $f$ is a proper log smooth integral morphism having log smooth parameter 
in $\LB$ and $\iota$ is a homeomorphic exact closed immersion 
in $\pLF$. Then, for a locally free isocrystal $\cE$ on 
$(X/\cY)^{\log}_{\conv}$ and $q \geq 0$, 
$R^qf_{X/\cY,\an *}\cE$ is a coherent 
sheaf on $]Y[^{\log}_{\cY} = \cY_K$ and we have the isomorphism 
$\sp_*R^qf_{X/\cY,\an *}\cE = R^qf_{X/\cY,\conv *}\cE$. 
\end{thm}

\begin{pf}
By definition of relative log analytic cohomology and relative log 
convergent cohomology, we have the quasi-isomorphism 
$$ 
Rf_{X/\cY,\conv *}\cE = R\,\sp_* Rf_{X/\cY,\an *}\cE. 
$$ 
So we have the spectral sequence 
$$ 
E_2^{s,t} = R^s\sp_*R^tf_{X/\cY,\an *}\cE \,\Longrightarrow\, 
R^{s+t}f_{X/\cY,\conv *}\cE. 
$$
Now we prove the theorem by induction on $q$. Assume the theorem is 
true up to $q-1$. Then we have 
$R^s\sp_*R^tf_{X/\cY,\an *}\cE = 0$ for $s>0, t<q$. So 
we have 
\begin{equation}\label{spstar}
\sp_*R^qf_{X/\cY,\an *}\cE = R^qf_{X/\cY,\conv *}\cE. 
\end{equation}
Now let us take a strict morphism 
$(\cY',M_{\cY'}) \lra (\cY,M_{\cY})$ such that the induced morphism 
$\cY'_K \lra \cY_K$ is an affinoid admissible open immersion and 
denote the base change of \eqref{diag-coherence0} by 
the morphism $(\cY',M_{\cY'}) \lra (\cY,M_{\cY})$ by 
\begin{equation}
(X',M_{X'}) {\lra} (Y',M_{Y'}) \os{\iota}{\hra} (\cY',M_{\cY'}). 
\end{equation}
Let us denote the restriction of $\cE$ to $I_{\conv}((X'/\cY')^{\log})$ by 
$\cE'$. Then we have, by the above argument and analytically flat 
base change theorem of 
relative log convergent cohomology, the isomorphism 
$$ 
\sp_*((R^qf_{X/\cY,\an *}\cE) |_{\cY'_K}) = 
R^qf_{X'/\cY',\conv *}\cE' = R^qf_{X/\cY,\conv *}\cE \otimes_{\cO_{\cY}} 
\cO_{\cY'}. 
$$
So $R^qf_{X/\cY,\an *}\cE$ is the sheaf of $\cO_{\cY_K}$-modules whose 
value at $\cY'_K$ is equal to 
$\Gamma(\cY', 
R^qf_{X/\cY,\conv *}\allowbreak 
\cE \allowbreak \otimes_{\cO_{\cY}} \cO_{\cY'})$ 
for any affinoid admissible open set of the form $\cY'_K \subseteq \cY_K$. 
This shows that $R^qf_{X/\cY,\an *}\cE$ is the coherent sheaf on 
$\cY_K$ induced by $R^qf_{X/\cY,\conv *}\cE$. So the assertion is 
proved. 
\end{pf}

\begin{cor}\label{coherence}
Assume we are given a diagram 
\begin{equation}
(X,M_X) \os{f}{\lra} (Y,M_Y) \os{\iota}{\hra} (\cY,M_{\cY}), 
\end{equation}
where $f$ is a proper log smooth integral morphism having log smooth parameter 
in $\LB$ and $\iota$ is a closed immersion in $\pLF$. 
Then, for a locally free isocrystal $\cE$ on 
$(X/\cY)^{\log}_{\conv}$, $R^qf_{X/\cY,\an *}\cE$ is a coherent 
sheaf on $]Y[^{\log}_{\cY}$. 
\end{cor}

\begin{pf}
By the argument in Remark \ref{rem3-1}, we can reduce to the case where 
$\iota$ is a homeomorphic 
exact closed immersion, and in 
this case, the claim is already proven in Theorem \ref{coherence0}. 
\end{pf}

Next we prove the existence of 
a structure of isocrystal on relative log analytic 
cohomology. To prove this, first we discuss a structure of an isocrystal on 
relative log convergent cohomology. 

\begin{thm}\label{iso1}
Assume we are given a diagram 
$$ (X,M_X) \os{f}{\lra} (Y,M_Y) \os{g}{\lra} (\cS,M_{\cS}), $$
where $f$ is a proper log smooth integral morphism having log smooth parameter 
in $\LB$ and $g$ is a morphism in $\pLF$. Then, for a locally free isocrystal 
$\cE$ in $I_{\conv}((X/\cS)^{\log})$ and a non-negative integer $q$, 
there exists a unique isocrystal $\cF$ on $(Y/\cS)^{\log}_{\conv}$ 
satisfying the following condition$:$ 
For any pre-widening $\cZ := ((\cZ,M_{\cZ}),(Z,M_Z),i,z)$ such that 
$z$ is a strict morphism and that $(\cZ,M_{\cZ})$ is formally log smooth 
over $(\cS,M_{\cS})$, the restriction of $\cF$ to 
$I_{\conv}((Z/\cS)^{\log}) \cong \Strat'((Z \hra \cZ/\cS)^{\log})$ is 
functorially given by 
$\{(R^qf_{X \times_Y Z_n/T_n(\cZ),\conv *}\cE, \allowbreak 
\epsilon_n)\}_n$, 
where 
$\{T_n(\cZ)\}_n := \{((T_n(\cZ),M_{T_n(\cZ)}),(Z_n,M_{Z_n}))\}_n$ is 
the system of universal enlargements of $\cZ$ and $\epsilon_n$ is the 
isomorphism 
{\small{ 
$$ 
p^*_{2,n}R^qf_{X \times_Y Z_n/T_n(\cZ), \conv *}\cE 
\os{\simeq}{\rightarrow} 
R^qf_{X \times_Y Z(1)_n/T_n(\cZ(1)), \conv *}\cE \os{\simeq}{\leftarrow} 
p^*_{1,n}R^qf_{X \times_Y Z_n/T_n(\cZ), \conv *}\cE. $$ }}
$($Here $(\cZ(1),M_{\cZ(1)}):=(\cZ,M_{\cZ}) \times_{(\cS,M_{\cS})} 
(\cZ,M_{\cZ})$, 
$$\{T_n(\cZ(1))\}_n := \{((T_n(\cZ(1)),M_{T_n(\cZ(1))}), 
(Z(1)_n, M_{Z(1)_n}))\}_n$$ 
is the system of universal enlargements of 
$((\cZ(1),M_{\cZ(1)}),(Z,M_Z))$ and $p_{i,n}$ is the morphism 
$(T_n(\cZ(1)),M_{T_n(\cZ(1))}) \lra (T_n(\cZ),M_{T_n(\cZ)})$ induced by 
the $i$-th projection.$)$ 
\end{thm} 

\begin{pf}
First let us define $\cF$ in the case where there exists a closed immersion 
$(Y,M_Y) \os{\iota}{\hra} (\cP,M_{\cP})$ in $\pLF$ over $(\cS,M_{\cS})$ 
such that $(\cP,M_{\cP})$ is formally log smooth over $(\cS,M_{\cS})$. 
Let $(\cP(i),M_{\cP(i)}) \,(i=0,1,2)$ be the $(i+1)$-fold fiber product 
of $(\cP,M_{\cP})$ over $(\cS,M_{\cS})$. 
Let $\{T_n(\cP(i))\} := 
\{((T_n(\cP(i)),M_{T_n(\cP(i))}),(Z_n(i),M_{Z_n(i)}))\}_n$ 
be the system of universal enlargements of $((\cP(i),M_{\cP(i)}),(Y,M_Y))$ 
and let us put 
$\cF_n(i) \allowbreak := R^qf_{X \times_Y Z_n(i)/T_n(\cP(i)),\conv *}\cE$. 
Let us note first that the transition morphisms 
$(T_n(\cP(0)),M_{T_n(\cP(0))}) \allowbreak 
\lra (T_{n'}(\cP(0)), \allowbreak M_{T_{n'}(\cP(0))}) \, (n<n')$ are 
analytically flat by Lemma \ref{anflex0} (2). So, by 
analytically flat base change theorem, the family 
$\{\cF_n(0)\}_n$ defines a compatible family of isocoherent sheaves 
on $\{T_n(\cP(0))\}_n$. Let us note next that 
the projections 
$(T_n(\cP(i+1)),M_{T_n(\cP(i+1))}) \allowbreak 
\lra (T_n(\cP(i)),M_{T_n(\cP(i))})$ 
are analytically flat by Proposition \ref{anflex}. So, again by 
analytically flat base change theorem, $\{\cF_n(1)\}_n$ and 
$\{\cF_n(2)\}_n$ 
induces a structure of compatible family of stratifications 
on $\{\cF_n(0)\}_n$. In this way, 
$\{\cF_n(i)\}_{n,i}$ induces an object $\cF$ in 
$\Strat'((Y \allowbreak 
\hra \cP/\cS)^{\log}) = I_{\conv}((Y/\cS)^{\log}).$ \par 
Now we define the isocrystal $\cF$ in general case. 
Take an embedding system 
\begin{equation}\label{emb3}
(Y,M_{Y}) \overset{g}{\lla} (Y^{(\bullet)}, M_{Y^{(\bullet)}}) 
{\hra} (\cP^{(\bullet)},M_{\cP^{(\bullet)}})
\end{equation}
such that $g$ is a strict formally etale \v{C}ech hypercovering and that 
$(\cP^{(i)},M_{\cP^{(i)}})$ is the $(i+1)$-fold fiber product 
of $(\cP^{(0)},M_{\cP^{(0)}})$ over $(\cS,M_{\cS})$. 
Then, by the construction in the previous paragraph for 
$(Y^{(i)},M_{Y^{(i)}}) \hra (\cP^{(i)},M_{\cP^{(i)}})$, we have 
$\cF^{(i)} \in I_{\conv}((Y^{(i)}/\cS)^{\log}) \, (i=0,1,2)$ and they are 
compatible thanks to Proposition \ref{anflex} and 
analytically flat base change theorem. 
So, by etale descent of isocrystals on relative log convergent site, 
$\{\cF^{(i)}\}_{i=0,1,2}$ descents to an isocrystal $\cF$ on 
$(Y/\cS)^{\log}_{\conv}$. \par 
Note that the isocrystal $\cF$ constructed in the previous paragraph 
satisfies the required property for $\cZ = ((\cP^{(i)},M_{\cP^{(i)}}), 
(Y^{(i)},M_{Y^{(i)}})) \,(i=0,1,2)$. By etale descent for isocrystals 
on relative log convergent site, this property characterizes $\cF$. 
So we have proved the uniqueness of $\cF$. \par 
Next, let us take a pre-widening $\cZ:=((\cZ,M_{\cZ}),(Z,M_Z),i,z)$ as 
in the statement of the theorem and we check the required property for $\cF$. 
Let us consider the following diagram of functors 
\begin{equation}\label{funct0}
\begin{CD}
I_{\conv}((Z/\cS)^{\log}) @>{{g'}^*}>> 
I_{\conv}((Z \times_{Y} Y^{(\b)}/\cS)^{\log}) \\ 
@A{z^*}AA @A{{z'}^*}AA \\ 
I_{\conv}((Y/\cS)^{\log}) @>{g^*}>> I_{\conv}((Y^{(\b)}/\cS)^{\log}), 
\end{CD}
\end{equation} 
where $I_{\conv}((Z \times_{Y} Y^{(\b)}/\cS)^{\log}), 
I_{\conv}((Y^{(\b)}/\cS)^{\log})$ denotes the category of descent data 
with respect to $I_{\conv}((Z \times_{Y} Y^{(n)}/\cS)^{\log}), 
I_{\conv}((Y^{(n)}/\cS)^{\log}) \, (n=0,1,2)$ respectively and 
$g^*, {g'}^*$ (resp. $z^*,{z'}^*$) are the functors induced by $g$ (resp. 
$z$). Then $g^*, {g'}^*$ are equivalence of categories. Let 
$$ \{T_n(\cZ)\}=\{((T_n(\cZ),M_{T_n(\cZ)}), (Z_n,M_{Z_n}))\}, $$
$$ \{T_n(\cP^{(\b)})\}=\{((T_n(\cP^{(\b)}),M_{T_n(\cP^{(\b)})}), 
(Y^{(\b)}_n,M_{Y^{(\b)}_n}))\}, $$
$$ \{T_n(\cZ \times \cP^{(\b)})\}=
\{((T_n(\cZ \times \cP^{(\b)}),M_{T_n(\cZ \times \cP^{(\b)})}), 
((Z \times_Y Y^{(\b)})_n,M_{(Z \times_Y Y^{(\b)})_n}))\}$$ be 
the system of universal enlargements of the widenings 
$\cZ, \cP^{(\b)}, \cZ \times \cP^{(\b)}$ (the product is taken in the 
category of widenings) respectively. Then, by definition, $\cF$ is sent 
by $g^*$ to the object 
$\{ (Rf_{X \times_Y Y^{(\b)}_n/T_n(\cP^{(\b)}), \conv *}\cE, 
\epsilon_{Y^{(\b)},n}) \}_n$ (where $\epsilon_{Y^{(\b)},n}$ is defined as in 
$\epsilon_n$ in the statement of the theorem), and by analytically flat 
base change, it is sent by ${z'}^*$ to the object 
$\{ (Rf_{X \times_Y (Z \times_Y Y^{(\b)})_n/
T_n(\cZ \times \cP^{(\b)}), \conv *}\cE, 
\epsilon_{Z \times_Y Y^{(\b)},n}) \}_n$ 
(where $\epsilon_{Z \times Y^{(\b)},n}$ 
is also defined as in $\epsilon_n$ in the statement of the theorem). 
On the other hand, we have an object 
$\{(R^qf_{X \times_Y Z_n/T_n(\cZ),\conv *}\cE, \epsilon_n)\}_n$ in 
$I_{\conv}((Z/\cS)^{\log})$, and by analytically flat base change theorem, 
it is sent by ${g'}^*$ also to 
 $\{ (Rf_{X \times_Y (Z \times_Y Y^{(\b)})_n/
T_n(\cZ \times \cP^{(\b)}), \conv *}\cE, 
\epsilon_{Z \times_Y Y^{(\b)},n}) \}_n$. 
So, by the above diagram, we conclude that the restriction of $\cF$ to 
$I_{\conv}((Z/\cS)^{\log})$ (by $z$) is given by 
$\{(R^qf_{X \times_Y Z_n/T_n(\cZ),\conv *}\cE, \epsilon_n)\}$, 
as required. \par 
Finally we prove the functoriality of the above expression. 
Let 
$$ \varphi: \cZ' := ((\cZ',M_{\cZ'}),(Z',M_{Z'}),i',z') \lra \cZ := 
((\cZ,M_{\cZ}),(Z,M_Z),i,z)$$ 
be a morphism of pre-widenings such that 
$z,z'$ are strict and that $(\cZ,M_{\cZ}), (\cZ',M_{\cZ'})$ are formally 
log smooth over $(\cS,M_{\cS})$, and we put 
$\cF_{\cZ} := \{(R^qf_{X \times_Y Z_n/T_n(\cZ),\conv *}\cE, \epsilon_n)\}_n$, 
$\cF_{\cZ'} \allowbreak := \allowbreak 
\{(R^qf_{X \times_Y Z'_n/T_n(\cZ'),\conv *}\cE, \epsilon'_n)\}_n$, 
where $Z'_n, T_n(\cZ'),\epsilon'_n$ are defined from $\cZ'$
in analogous way as $Z_n, T_n(\cZ), \epsilon_n$. 
Then, we have two isomorphisms of the form 
\begin{equation}\label{func1}
\varphi^*\cF_{\cZ} \lra \cF_{\cZ'}: 
\end{equation}
One is the isomorphism 
induced by the functoriality of relative log convergent 
cohomology, and the other is the isomorphism induced by the isocrystal 
structure of $\cF$. We should prove the coincidence of them. 
Let us 
put $(\cZ'',M_{\cZ''}) := (\cZ,M_{\cZ}) \times_{(\cS,M_{\cS})} 
(\cZ',M_{\cZ'})$ and denote the projections 
$(\cZ'',M_{\cZ''}) \lra (\cZ,M_{\cZ}), (\cZ'',M_{\cZ''}) \lra 
(\cZ',M_{\cZ'})$ by $\pr_1,\pr_2$, respectively. Let 
$\gamma: (\cZ',M_{\cZ'}) \lra (\cZ'',M_{\cZ''})$ be the graph of 
$\varphi$ and put 
$\cF_{\cZ''} := 
\{(R^qf_{X \times_Y Z''_n/T_n(\cZ''),\conv *}\cE, \epsilon''_n)\}$ 
(where $Z''_n, T_n(\cZ''),\epsilon''_n$ are defined from 
$((\cZ'',M_{\cZ''}),(Z',M_{Z'}))$
in analogous way to $Z_n, T_n(\cZ), \epsilon_n$). Then we have 
two diagrams of the form 
\begin{equation}
\varphi^*\cF_{\cZ} = \gamma^*\pr_1^*\cF_{\cZ} 
\os{\simeq}{\lra} \gamma^*\cF_{\cZ''} 
\os{\simeq}{\lla} \gamma^*\pr_2^*\cF_{\cZ'} = \cF_{\cZ'}: 
\end{equation} 
One is the diagram induced by the functoriality of relative log convergent 
cohomology, and the other is the isomorphism induced by the isocrystal 
structure of $\cF$. In both cases, the composite is equal to 
the isomorphism \eqref{func1}. So, to prove the coincidence of the maps 
\eqref{func1}, we may replace $\varphi$ by $\pr_i \,(i=1,2)$, 
that is, we may assume that $(\cZ',M_{\cZ'})$ is 
formally log smooth over $(\cZ,M_{\cZ})$. In this case, we have the 
following diagram of functors 
\begin{equation*}
\begin{CD}
I_{\conv}((Z'/\cS)^{\log}) @>{{g'}^*}>> 
I_{\conv}((Z' \times_Y Y^{(\b)}/\cS)^{\log}) \\ 
@A{\varphi^*}AA @A{{\varphi'}^*}AA \\ 
I_{\conv}((Z/\cS)^{\log}) @>{g^*}>>
I_{\conv}((Z \times_Y Y^{(\b)}/\cS)^{\log}), 
\end{CD}
\end{equation*}
where the horizontal lines are as the top horizontal line 
in the diagram \eqref{funct0} and 
 $\varphi^*,{\varphi'}^*$ are induced by $\varphi^*$. 
Let us define $\cF_{\cZ \times \cP^{(\b)}} \in 
I_{\conv}((Z \times_Y Y^{(\b)}/\cS)^{\log}), 
\cF_{\cZ' \times \cP^{(\b)}} \in 
I_{\conv}((Z' \times_Y Y^{(\b)}/\cS)^{\log})$ in the same way 
as $\cF_{\cZ}, \cF_{\cZ'}$. 
Then we have the two commutative diagrams of the form 
\begin{equation}\label{funct2} 
\begin{CD}
{g'}^*\cF_{\cZ'} @>>> \cF_{\cZ' \times \cP^{(\b)}} \\ 
@AAA @AAA \\ 
{g'}^*\varphi^*\cF_{\cZ} @>>> {\varphi'}^*\cF_{\cZ \times \cP^{(\b)}}, 
\end{CD}
\end{equation}
where the left vertical arrow is given by applying ${g'}^*$ to the maps 
\eqref{func1}, that is, the map in one diagram is given by
 the isocrystal structure and the map in another diagram is 
 given by the functoriality of relative log analytic cohomology. 
Other arrows in the diagrams \eqref{funct2} 
are given by either of them, which coincide by definition 
of $\cF$ given above. Since all the arrows in \eqref{funct2} are 
isomorphisms and ${g'}^*$ is an equivalence, we conclude from 
two diagrams \eqref{funct2} that the two isomorphisms \eqref{func1} 
are equal. So we have proved the functoriality and the proof of 
the theorem is now finished. 
\end{pf} 

\begin{rem}
In the above situation, 
we do not claim that the value $\cF_{\cZ}$ of $\cF$ on any 
enlargement $\cZ := ((\cZ,M_{\cZ}),(Z,M_Z),i,z)$ 
is given by $R^qf_{X \times_Y Z/\cZ,\conv *}\cE$. 
If we know that $R^qf_{X \times_Y Z/\cZ,\conv *}\cE$ is locally free 
for any enlargement $\cZ := ((\cZ,M_{\cZ}),(Z,M_Z),i,z)$, the equality 
$\cF_{\cZ} = R^qf_{X \times_Y Z/\cZ,\conv *}\cE$ holds for any 
enlargement $\cZ := ((\cZ,M_{\cZ}),(Z,M_Z),i,z)$. (The proof is 
similar to the above proof, using Corollary \ref{bc3} under the 
assumption (1)' and (2).) 
This is the case if $(\cS,M_{\cS}) = (\Spf V, \tls)$ 
(where $V$ is a complete discrete valuation ring 
of mixed characteristic $(0,p)$ with residue field $k$) 
and $M_Y$ is trivial (cf. \cite{ogus2}). 
\end{rem}

As a corollary, we have the existence of 
a structure of an isocrystal on relative log analytic 
cohomology. 

\begin{cor}\label{iso2}
Assume we are given a diagram 
$$ (X,M_X) \os{f}{\lra} (Y,M_Y) \os{g}{\lra} (\cS,M_{\cS}), $$
where $f$ is a proper log smooth integral 
morphism having log smooth parameter 
in $\LB$ and $g$ is a morphism in $\pLF$. Then, for a locally free isocrystal 
$\cE$ on $I_{\conv}((X/\cS)^{\log})$ and a non-negative integer $q$, 
there exists a unique isocrystal $\cF$ on $(Y/\cS)^{\log}_{\conv}$ 
such that, for any pre-widening 
$\cZ:=((\cZ,M_{\cZ}),(Z,M_Z),i,z)$ of $(Y,M_Y)/(\cS,M_{\cS})$
 such that $z$ is strict and $(\cZ,M_{\cZ})$ is formally 
 log smooth over 
$(\cS,M_{\cS})$, $\cF$ induces, via the functor 
$$ I_{\conv}((Y/\cS)^{\log}) \lra 
I_{\conv}((Z/\cS)^{\log}) \simeq \Strat''((Z \hra \cZ/\cS)^{\log}), $$
an object of the form 
$(R^qf_{X \times_Y Z/\cZ,\an *}\cE, \epsilon)$, where 
$\epsilon$ is the canonical isomorphism 
$$ p_2^*R^qf_{X \times_Y Z/\cZ,\an *}\cE  \os{\sim}{\ra} 
R^qf_{X \times_Y Z/\cZ(1),\an *}\cE \os{\sim}{\leftarrow}
p_1^*R^qf_{X \times_Y Z/\cZ,\an *}\cE. $$ 
$($Here $(\cZ(1),M_{\cZ(1)}) := (\cZ,M_{\cZ}) \times_{(\cS,M_{\cS})} 
(\cZ,M_{\cZ})$ and $p_i$ 
denotes the $i$-th projection  
$]Z[^{\log}_{\cZ(1)} \allowbreak 
\lra \allowbreak ]Z[^{\log}_{\cZ}.)$ 
\end{cor}

\begin{pf}
It is immediate from Theorem \ref{iso1} and 
the relation between relative log convegent cohomology and 
relative log analytic cohomology. 
\end{pf}


\section{Relative log analytic cohomology and relative rigid cohomology}

In this section, 
we prove a comparison 
theorem between relative log analytic cohomology and relative rigid 
cohomology (defined in \cite{berthelot1} and \cite{chts}) for 
certain proper log smooth integral morphism having 
log smooth parameter. 
By using it, we prove a version of a conjecture of Berthelot on 
the coherence and the overconvergence of 
relative rigid cohomology with coefficient in special case, 
that is, the case where the given situation admits a `nice' log structure. 
We also discuss on the Frobenius structure on relative rigid cohomology 
when a given coefficient admits a Frobenius structure. \par 
Throughout this section, we assume that the log structure $M_{\cB}$ on 
the base log formal scheme $\cB$ is trivial. \par 
First let us recall some terminologies of \cite{chts} and the definition 
of relative rigid cohomology given in \cite{chts}. (See also 
\cite{berthelot2}, \cite{shiho2}. Note that some notations here are 
different from those in \cite{chts}.) 
A pair $(X,\ol{X})$ is, by definition, a pair of schemes 
$X,\ol{X}$ of characteristic $p$ endowed with an open immersion 
$X \hra \ol{X}$. A map of pairs $f:(X,\ol{X}) \lra (Y,\ol{Y})$ is 
a morphism of schemes $f:\ol{X} \lra \ol{Y}$ satisfying $f(X) \subseteq Y$. 
$f$ is called separated of finite type 
if so is $f:\ol{X} \lra \ol{Y}$. 
A pair $(Y,\ol{Y})$ over a given pair $(X,\ol{X})$ is a pair endowed 
with structure 
morphism $(Y,\ol{Y}) \lra (X,\ol{X})$. In this paper, a pair 
is always assumed to be 
a pair over $(B,B)$ whose structure morphism is separated of finite 
type and all the morphisms of pairs are assumed to be separated 
morphisms of finite type over $(B,B)$. A morphism of pairs 
$f: (X,\ol{X}) \lra (Y,\ol{Y})$ is called strict if $f^{-1}(Y)=X$ holds. 
A triple $(X,\ol{X},\cX)$ is a pair $(X,\ol{X})$ endowed with 
a $p$-adic formal scheme $\cX$ over $\Z_p$ and a closed immersion 
$\ol{X} \hra \cX$ over $\Z_p$. A map of triples 
$f:(X,\ol{X},\cX) \lra (Y,\ol{Y},\cY)$ 
is defined in an obvious way. $f$ is called separated of finite type 
if so is $\cX \lra \cY$. A triple $(Y,\ol{Y},\cY)$ over a given triple 
$(X,\ol{X},\cX)$ is a triple endowed with structure morphism 
$(Y,\ol{Y},\cY) \lra (X,\ol{X},\cX)$. 
In this paper, all the triples 
are assumed to be 
treiples over $(B,B,\cB)$ whose structure morphism is separated of finite 
type and all the morphisms of triples are assumed to be separated 
morphisms of finite type over $(B,B,\cB)$. 
For a triple $(S,\ol{S},\cS)$ and a pair $(X,\ol{X})$ over 
$(S,\ol{S})$, an $(X,\ol{X})$-triple over $(S,\ol{S},\cS)$ is a triple 
$(Y,\ol{Y},\cY)$ over $(S,\ol{S},\cS)$ endowed with a morphism 
of pairs $(Y,\ol{Y}) \lra (X,\ol{X})$ over $(S,\ol{S})$. \par 
Now let us assume given morphisms 
\begin{equation}\label{chts1}
f: (X,\ol{X}) \lra (Y,\ol{Y}), \qquad \iota: \ol{Y} \hra \cY, 
\end{equation} 
where $f$ is a morphism of pairs and $\iota$ is a closed immersion 
of $\ol{Y}$ into a $p$-adic formal $\cB$-scheme $\cY$ over 
$\cB$. Assume for the moment that the above morphism fits into 
a diagram 
\begin{equation}\label{chts2}
\begin{CD}
\ol{X} @>i>> \cP \\ 
@VfVV @V{f'}VV \\ 
\ol{Y} @>{\iota}>> \cY, 
\end{CD}
\end{equation}
where $i$ is a closed immersion of $\ol{X}$ into a $p$-adic formal 
$\cB$-scheme $\cP$ and 
$f'$ is a morphism which is formally smooth on a neighborhood of $X$.
For $n \in \N$, we denote by $\cP(n)$ the $(n+1)$-fold fiber product of 
$\cP$ over $\cY$ and by $i(n): \ol{X} \hra \cP(n)$ the natural closed 
immersion induced by $i$. Then we can form admissible open immersions 
$j_X:\,]X[_{\cP(n)} \hra ]\ol{X}[_{\cP(n)}$ of tubular neighborhoods and 
we have the projections 
$$ p_i\,:\,]\ol{X}[_{\cP(1)} \lra ]\ol{X}[_{\cP} \,\, (i=1,2),\qquad 
p_{ij}\,:\,]\ol{X}[_{\cP(2)} \lra ]\ol{X}[_{\cP(1)} \,\, 
(1 \leq i < j \leq 3) $$ 
and the diagonal morphism 
$\Delta\,:\,]\ol{X}[_{\cP} \hra \, ]X[_{\cP(1)}.$ 
Then the category $I^{\dag}((X,\ol{X})/\cY_K, \cP)$ of 
realization of overconvergent isocrystals of $(X,\ol{X})/\cY_K$ 
over $\cP$ is defined to be the category of pairs $(\cE,\epsilon)$, 
where $\cE$ is a coherent $j^{\dag}\cO_{]\ol{X}[_{\cP}}$-module and 
$\epsilon$ is a $j^{\dag}\cO_{]\ol{X}[_{\cP(1)}}$-linear 
isomorphism $p_2^*\cE \os{\sim}{\lra} p_1^*\cE$ satisfying 
$\Delta^*(\epsilon)=\id, p_{12}^*(\epsilon) \circ p_{23}^*(\epsilon) = 
p_{13}^*(\epsilon)$. (This definition is different from that given in 
\cite[10.2]{chts}, but they are equivalent by \cite[Lemma 10.2.2]{chts}.) 
When we are given an object $\cE := (\cE,\epsilon) \in 
I^{\dag}((X,\ol{X})/\cY_K, \cP)$, we can form an integrable connection of 
the form $\cE \lra \cE \otimes j_X^{\dag}
\Omega^1_{]\ol{X}[_{\cP}/\cY_K}$ by using $\epsilon$, and it induces 
the de Rham complex of the form 
$\cE \otimes j_X^{\dag} \Omega^{\b}_{]\ol{X}[_{\cP}/\cY_K}$, which 
we denote by $\DR^{\dag}(]\ol{X}[_{\cP}/\cY_K,\cE)$. \par 
In the case where there does not necessarily have a diagram 
\eqref{chts2}, we can always find a diagram 
\begin{equation}\label{chts3} 
\ol{X} \os{g}{\lla} \ol{X}^{(\b)} \os{i^{(\b)}}{\hra} \cP^{(\b)}, 
\end{equation}
where $g$ is a Zariski hypercovering, $i^{(\b)}$ is a closed immersion 
over $\cY$ such that each $\cP^{(n)}$ is a $p$-adic formal $\cB$-scheme 
which is formally smooth on a neighborhood of $X^{(n)}:= X 
\times_{\ol{X}} \ol{X}^{(n)}$. Then, we can form the category 
$I^{\dag}((X^{(\b)},\ol{X}^{(\b)})/\cY_K, \cP^{(\b)})$ of descent data 
with respect to the categories 
$I^{\dag}((X^{(n)},\ol{X}^{(n)})/\cY_K,\cP^{(n)}) \, (n=0,1,2)$ and it is 
known that this category is independent of the choice of the diagram 
\eqref{chts3}. The category of overconvergent isocrystals on 
$(X,\ol{X})/\cY_K$ is defined to be this category and it is simply 
denoted by $I^{\dag}((X,\ol{X})/\cY_K)$. \par 
For an object 
$\cE := (\cE^{(\b)},\epsilon^{(\b)})$ in 
$I^{\dag}((X,\ol{X})/\cY_K) = 
I^{\dag}((X^{(\b)},\ol{X}^{(\b)})/\cY_K, \cP^{(\b)})$, 
we have the associated de Rham complex 
$\DR^{\dag}(]\ol{X}^{(\b)}[_{\cP^{(\b)}}/\cY_K,\cE^{(\b)})$ on 
$]\ol{X}^{(\b)}[_{\cP^{(\b)}}$. Let $h$ be the morphism 
$]\ol{X}^{(\b)}[_{\cP^{(\b)}} \lra ]\ol{Y}[_{\cY}$ induced by $f \circ g$. 
Then we define $Rf_{(X,\ol{X})/\cY,\rig *}\cE, \allowbreak 
R^qf_{(X,\ol{X})/\cY,\rig *}\cE$ by 
\begin{align*}
Rf_{(X,\ol{X})/\cY,\rig *}\cE &:= Rh_*
\DR^{\dag}(]\ol{X}^{(\b)}[_{\cP^{(\b)}}/\cY_K, \cE^{(\b)}), \\
R^qf_{(X,\ol{X})/\cY,\rig *}\cE &:= R^qh_*
\DR^{\dag}(]\ol{X}^{(\b)}[_{\cP^{(\b)}}/\cY_K, \cE^{(\b)}), 
\end{align*}
respectively and we call $R^qf_{(X,\ol{X})/\cY,\rig *}\cE$ the $q$-th 
relative rigid cohomology of $(X,\ol{X})/\cY$ with coefficient $\cE$. 
It is a sheaf of $j_Y^{\dag}\cO_{]\ol{Y}[_{\cY}}$-modules, where 
$j_Y$ denotes the admissible open immersion $]Y[_{\cY} \hra 
]\ol{Y}[_{\cY}$. \par 
The category of overconvergent isocrystals 
$I^{\dag}((X,\ol{X})/\cY)$ satisfies the descent property 
for Zariski covering of $\ol{X}$ 
(it is rather easy consequence from the definition given above). 
Here we give a proof of the descent property of $I^{\dag}((X,\ol{X})/\cY)$ 
for etale covering of $\ol{X}$: 

\begin{prop}\label{etaledescent}
Assume we are given morphisms \eqref{chts1} and let 
$\ol{X}^{(\b)} \lra \ol{X}$ be an etale hypercovering. 
Let us put $X^{(\b)} := X \times_{\ol{X}} \ol{X}^{(\b)}$ and 
let us denote the 
category of descent data with respect to the categories 
$I^{\dag}((X^{(n)},\ol{X}^{(n)})/\cY_K) \,(i=0,1,2)$ by 
$I^{\dag}((X^{(\b)}, \ol{X}^{(\b)})/\cY_K)$. Then the restriction functor 
\begin{equation}\label{rest}
I^{\dag}((X,\ol{X})/\cY_K) \lra I^{\dag}((X^{(\b)},\ol{X}^{(\b)})/\cY_K)
\end{equation}
is an equivalence of categories. 
\end{prop} 

\begin{rem} 
In the case where $(Y,\ol{Y},\cY) = (\Spec k, \Spec k, \Spf V)$ (where 
$V$ is a complete discrete valuation ring of mixed characteristic 
with residue field $k$) and 
$\ol{X}$ is proper over $k$, this is a special case of a result of 
Etesse \cite{etesse}. (Etesse treats the category of overconvergent 
$F$-isocrystals, 
but his argument works also in the case without Frobenius structure.) 
\end{rem} 

\begin{pf} 
By using \cite[3.3.4.2]{saintdonat}, 
one can reduce the proposition to the case where 
$\ol{X}^{(\b)}$ is the \v{C}ech hypercovering associated to the etale 
surjective morphism $\ol{X}^{(0)} \lra \ol{X}$. Since the assertion is 
Zariski local on $\ol{X}$, we may assume the existence of the diagram 
\eqref{chts2} such that $\ol{X}$ and $\cP$ are affine. 
Moreover, since one can replace $\ol{X}^{(0)}$ by 
a Zariski refinement of it, we may assume 
(by \cite[claim in p.81]{shiho2}) the existence of the Cartesian diagram 
\begin{equation}
\begin{CD} 
\ol{X}^{(0)} @>{i^{(0)}}>> \cP^{(0)} \\ 
@VVV @VVV \\
\ol{X} @>i>> \cP
\end{CD}
\end{equation} 
for some affine $p$-adic formal $\cB$-scheme $\cP^{(0)}$ 
formally etale over $\cP$. 
Let $\cP^{(\b)}$ be the \v{C}ech hypercovering associated to 
$\cP^{(0)} \lra \cP$ and define $A, B^{(\b)}$ by 
$\cP = \Spf A, \cP^{(\b)} = \Spf B^{(\b)}$. 
Let $\{f_i\}_{i \in I}$ be a set of generators of the ideal 
$\Ker (A \lra \Gamma(\ol{X},\cO_{\ol{X}}))$, let 
$\{\ol{g}_j\}_{j \in J}$ be a set of generators of the ideal 
$\Ker(\Gamma(\ol{X},\cO_{\ol{X}}) \lra \Gamma(\ol{X}-X,\cO_{\ol{X}-X}))$ 
and let $g_j \,(j \in J)$ be a lift of $\ol{g}_j$ to $A$. \par 
Under the above notation, first we prove that the restriction functor 
\begin{equation}\label{cohmod}
(\text{coherent $j_X^{\dag}\cO_{]\ol{X}[_{\cP}}$-modules}) 
\lra 
\left(
\begin{aligned}
& \text{compatible family of coherent} \\
& \text{$j_X^{\dag}\cO_{]\ol{X}^{(n)}[_{\cP}^{(n)}}$-modules \,$(n=0,1,2)$}
\end{aligned}
\right) 
\end{equation} 
is an equivalence of categories. Fix a strictly increasing 
sequence $\ul{\eta} := 
\{\eta_m\}$ in $\Q \otimes_{\Z} |K^{\times}|$ converging to $1$. 
For another strictly increasing sequence $\ul{\nu} := \{\nu_m\}$ in 
$\Q \otimes_{\Z} |K^{\times}|$ converging to $1$ with 
$\nu_m > \eta_m$ and $J' \subseteq J$, 
we define 
\begin{align*}
V_{\eta_m, \nu_m, J'} & := \{x \in ]X[_{\cP} \,\vert\, 
|f_i(x)|\leq \eta_m\,(i \in I), |g_j(x)| \geq \nu_m \,(j \in J')\}, \\ 
V_{\eta_m,\nu_m} & := \bigcup_{j \in J}V_{\eta_m,\nu_m,\{j\}}, 
\qquad V_{\ul{\eta},\ul{\nu}} := \bigcup_{m}V_{\eta_m,\nu_m} 
\end{align*}
and let $V^{(\b)}_{\eta_m, \nu_m, J'}, 
V^{(\b)}_{\eta_m,\nu_m}, V^{(\b)}_{\ul{\eta},\ul{\nu}}$ be the pull-back of 
$V_{\eta_m, \nu_m, J'}, 
V_{\eta_m,\nu_m}, V_{\ul{\eta},\ul{\nu}}$ by 
$]X^{(\b)}[_{\cP^{(\b)}} \lra ]X[_{\cP}$. Assume given an object 
$\{E^{(n)}\}_n$ in the right hand side \eqref{cohmod}. 
Then $\{E^{(n)}\}_n$ comes from a compatible family of coherent modules 
on $V^{(n)}_{\ul{\eta},\ul{\nu}}\,(n=0,1,2)$ 
for some $\ul{\nu}$. So we obtain the 
compatible family of coherent modules 
$\{\ti{E}^{(n)}_{m,J'}\}_{n,m,J'}$ on 
$V^{(\b)}_{\eta_m,\nu_m,J'}$. Now let us note that the morphism 
$V^{(\b)}_{\eta_m,\nu_m,J'} \lra V_{\eta_m,\nu_m,J'}$ comes from a 
formally etale \v{C}ech covering of certain $p$-adic formal $\cB$-schemes. 
So, by rigid analytic faifufully flat descent, the family 
of modules 
$\{\ti{E}^{(n)}_{m,J'}\}_n$ is obtained as the pull-back of 
a coherent module $\ti{E}_{m,J'}$ on $V_{\eta_m,\nu_m,J'}$ (for each 
$m,J'$). Since the family $\{\ti{E}_{m,J'}\}_{m,J'}$ 
are compatible with respect to $m,J'$, it is obtained as the pull-back 
of the module $\ti{E}$ on $V_{\ul{\eta},\ul{\nu}}$ and $E:= j^{\dag}\ti{E}$ 
is the object in the left hand side of \eqref{cohmod} which restricts to 
$E^{(\b)}$. So we proved the essential surjectivity of the functor 
\eqref{cohmod}. One can prove the full-faithfulness in the same way. 
So the functor \eqref{cohmod} is an equivalence of categories. \par 
Now let us prove the proposition. Let us given an object 
$\{E^{(n)}\}_n := \{(E^{(n)},\epsilon^{(n)})\}_n$ in 
$I^{\dag}((X^{(\b)},\ol{X}^{(\b)})/\cY) = 
I^{\dag}((X^{(\b)},\ol{X}^{(\b)})/\cY, \cP^{(\b)})$. Then 
$\{E^{(n)}\}_n$ defines a compatible family of 
coherent $j_X^{\dag}\cO_{]\ol{X}^{(n)}[_{\cP^{(n)}}}$-modules $(n=0,1,2)$. 
So it descents to a coherent $j_X^{\dag}\cO_{]\ol{X}[_{\cP}}$-module 
$E$. On the other hand, $\epsilon^{(n)}$ is 
an isomorphism $p_2^*E^{(n)} \os{\simeq}{\lra} p_1^*E^{(n)}$ of 
$j^{\dag}\cO_{]\ol{X}^{(n)}[_{\cP^{(n)} \times_{\cY} \cP^{(n)}}}$-modules 
and we know that $\{p_i^*E^{(n)}\}_n\,(i=1,2)$ descents to the 
$j^{\dag}_X\cO_{]\ol{X}[_{\cP \times_{\cY} \cP}}$-module $p_i^*E$. 
Now note that $]\ol{X}^{(\b)}[_{\cP^{(\b)} \times_{\cY} \cP^{(\b)}}$ is 
isomorphic to $]\ol{X}^{(\b)}[_{\cP \times_{\cY} \cP^{(\b)}}$ and that 
the diagram 
\begin{equation*}
\begin{CD}
\ol{X}^{(\b)} @>>> \cP \times_{\cY} \cP^{(\b)} \\ 
@VVV @VVV \\ 
\ol{X} @>>> \cP \times_{\cY} \cP
\end{CD}
\end{equation*} 
is Cartesian. So, by the equivalence of categories \eqref{cohmod} 
($\cP, \cP^{(\b)}$ replaced by $\cP \times_{\cY} \cP, 
\cP \times_{\cY} \cP^{(\b)}$ respectively), $\{\epsilon^{(n)}\}_n$ 
descents to 
an isomorphism $\epsilon: p_2^*E \lra p_1^*E$ of 
$j_X^{\dag}\cO_{]\ol{X}[_{\cP \times_{\cY} \cP}}$-modules. 
In this way we obtain an object $(E,\epsilon)$ in 
$I^{\dag}((X,\ol{X})/\cY_K)$. So we have shown the essential surjectivity 
of the functor \eqref{rest}. The full-faithfulness can be shown in the same 
way. So we have proved the assertion of the proposition. 
\end{pf} 

Now we recall the statement of the conjecture of Berthelot 
(\cite{berthelot1}). Note that we need some discussion here because 
we would like to introduce several versions of the conjecture. 
First we quote \cite[4.3]{berthelot1}, where the conjecture is stated: \\
\quad \\
Soient $S$ un $k$-sch\'{e}ma, $f:X \lra S$ un morphisme propre et lisse: 
supposons pour simplifier que $k$ soit parfait, et qu'il existe une 
compactification $S \hra \ol{S}$ poss\'{e}dant un plongement $\ol{S} \hra T$ 
dans un sch\'{e}ma formel lisse $T$ sur l'anneau $\cW$ de vecteurs de Witt 
\`{a} coefficients dans $k$ (par exemple si $S$ est quasi projectif sur 
$k$). On peut alors former, comme on l'a vu en (2.5) c), les faisceaux de 
cohomologie rigide relative $R^qf_{\rig *}(X/T)$; Si $j^{\d}$ est associ\'{e} 
comme (2.2) \`{a} l'inclusion $]S[_T \hra ]\ol{S}[_T$, ce sont de fa\c{c}on 
naturelle des $j^{\d}\cO_{]\ol{S}[}$-modules. Je conjecture que les 
$R^qf_{\rig *}(X/T)$ poss\'{e}dent une structure canonique de $F$-cristal 
surconvergent ($\Phi$ \'{e}tant induit par fonctorialit\'{e} par le 
Frobenius de $X$). ... \\
\quad \\
(Here, $R^qf_{\rig *}(X/T)$ is the relative rigid cohomology of 
the structure overconvergent isocrystal on $(X,\ol{X})/T$, 
where $\ol{X}$ denotes any compactification of $X$ over $\ol{S}$.) 
In the above form of the 
conjecture, $\ol{X}, \ol{S}$ is assumed to be complete, 
the coefficient is assumed to be trivial and the base formal scheme 
is assumed to be $\Spf \cW$. \par 
Next, take a triple of the form 
$(S,S,\cS)$ and assume that we are given a diagram 
of pairs 
$$ (X,\ol{X}) \os{f}{\lra} (Y,\ol{Y}) \os{g}{\lra} (S,S). $$ 
Then, based on \cite[4.3]{berthelot1}, Tsuzuki gives a generalized form 
of the conjecture (he calls it Berthelot's conjecture) in 
\cite[\S 4]{tsuzuki3}, whose statement 
is as follows: Assume that $f:\ol{X} \lra \ol{Y}$ is proper, $f^{-1}(Y)=X$ 
and $f|_X:X \lra Y$ is smooth. Then, for an overconvergent ($F$-)isocrystal 
$\cE$ on $(X,\ol{X})/\cS_K$ and $q \geq 0$, 
`the $q$-th rigid cohomology overconvergent ($F$-)isocrystal 
$Rf_{\rig *}\cE$' on $(Y,\ol{Y})/\cS_K$ 
(for definition, see \cite[3.3.1]{tsuzuki3} and Remark \ref{tsurem} below) 
exists. By examining the property which 
`the $q$-th rigid cohomology overconvergent ($F$-)isocrystal' should satisfy, 
we arrive at 
the following version of Berthelot's conjecture, which is stronger than 
the above-mentioned forms of Berthelot's conjecture: 

\begin{conj}\label{bconj1}
Let us take a triple of the form 
$(S,S,\cS)$ $($over $(B,B,\cB))$ 
and assume we are given a diagram of pairs 
$$ (X,\ol{X}) \os{f}{\lra} (Y,\ol{Y}) \os{g}{\lra} (S,S) $$ 
such that $f:\ol{X} \lra \ol{Y}$ is proper, $f^{-1}(Y)=X$ 
and that $f|_X:X \lra Y$ is smooth. Then, for an overconvergent 
$(F$-$)$isocrystal $\cE$ on $(X,\ol{X})/\cS_K$ and $q \geq 0$, 
there exists uniquely an overconvergent $(F$-$)$isocrystal $\cF$ 
$($which is called the $q$-th rigid cohomology overconvergent 
isocrystal$)$ on 
$(Y,\ol{Y})/\cS_K$ satisfying the following condition$:$ 
For any $(Y,\ol{Y})$-triple $(Z,\ol{Z},\cZ)$ over $(S,S,\cS)$ 
with $\cZ$ formally smooth over $\cS$ on a neighborhood of $Z$, 
the restriction of $\cF$ to $I^{\d}((Z,\ol{Z})/\cS_K, \cZ)$ is given 
functorially by 
$(R^qf_{(X \times_Y Z,\ol{X} \times_{\ol{Y}} \ol{Z})/\cZ, \rig *}\cE, 
\epsilon)$, where $\epsilon$ is given by 
$$ 
p_2^*R^qf_{(X \times_Y Z,\ol{X} \times_{\ol{Y}} \ol{Z})/\cZ, \rig *}\cE 
\os{\simeq}{\ra} 
R^qf_{(X \times_Y Z,\ol{X} \times_{\ol{Y}} \ol{Z})/\cZ \times_{\cS} \cZ, 
\rig *}\cE \os{\simeq}{\leftarrow} 
p_1^*R^qf_{(X \times_Y Z,\ol{X} \times_{\ol{Y}} \ol{Z})/\cZ, \rig *}\cE. 
$$ 
$($Here $p_i$ is the morphism 
$]\ol{Z}[_{\cZ \times_{\cS} \cZ} \lra \,]\ol{Z}[_{\cZ}$ induced by the 
$i$-th projection.$)$ 
\end{conj} 

\begin{rem} 
The conjecture contains 
the claim that the two morphisms used in the definition of $\epsilon$ are 
isomorphisms. 
\end{rem} 

The above conjecture is so strong that we do not know any non-trivial 
example for which Conjecture \ref{bconj1} is proved. We would like to 
consider the following version of Berthelot's conjecture, which is 
slightly weaker than Conjecture \ref{bconj1} but strong enough to assure 
the unique existence of `the $q$-th rigid cohomology overconvergent 
isocrystal': 

\begin{conj}\label{bconj2}
Take a triple of the form 
$(S,S,\cS)$ $($over $(B,B,\cB))$ 
and assume we are given a diagram of pairs 
$$ (X,\ol{X}) \os{f}{\lra} (Y,\ol{Y}) \os{g}{\lra} (S,S) $$ 
such that $f:\ol{X} \lra \ol{Y}$ is proper, $f^{-1}(Y)=X$ 
and that $f|_X:X \lra Y$ is smooth. 
Then 
there exists a subcategory $\cC$ of the category of 
$(Y,\ol{Y})$-triples over $(S,S,\cS)$ such that, 
for an overconvergent 
$(F$-$)$isocrystal $\cE$ on $(X,\ol{X})/\cS_K$ and $q \in \N$, 
there exists 
uniquely an overconvergent $(F$-$)$isocrystal $\cF$ 
$($which is called the $q$-th rigid cohomology overconvergent 
isocrystal$)$ on 
$(Y,\ol{Y})/\cS_K$ satisfying the following condition$:$ 
For any $(Z,\ol{Z},\cZ) \in \cC$ 
with $\cZ$ formally smooth over $\cS$ on a neighborhood of $Z$, 
the restriction of $\cF$ to $I^{\d}((Z,\ol{Z})/\cS_K, \cZ)$ is given 
functorially by 
$(R^qf_{(X \times_Y Z,\ol{X} \times_{\ol{Y}} \ol{Z})/\cZ, \rig *}\cE, 
\epsilon)$, where $\epsilon$ is given by 
$$ 
p_2^*R^qf_{(X \times_Y Z,\ol{X} \times_{\ol{Y}} \ol{Z})/\cZ, \rig *}\cE 
\os{\simeq}{\rightarrow} 
R^qf_{(X \times_Y Z,\ol{X} \times_{\ol{Y}} \ol{Z})/\cZ \times_{\cS} \cZ, 
\rig *}\cE \os{\simeq}{\leftarrow} 
p_1^*R^qf_{(X \times_Y Z,\ol{X} \times_{\ol{Y}} \ol{Z})/\cZ, \rig *}\cE. 
$$ 
$($Here $p_i$ is the morphism 
$]\ol{Z}[_{\cZ \times_{\cS} \cZ} \lra \,]\ol{Z}[_{\cZ}$ induced by the 
$i$-th projection.$)$ 
\end{conj} 

The difference of Conjecture \ref{bconj2} compared to Conjecture \ref{bconj1} 
is the introduction of the category $\cC$: Conjecture \ref{bconj2} is weaker 
than Conjecture \ref{bconj1} in the sense that $\cC$ need not be equal to 
the category of all $(Y,\ol{Y})$-triples $(Z,\ol{Z},\cZ)$ over 
$(S,S,\cS)$. However, 
Conjecture \ref{bconj2} is strong enough because it requires 
the unique existence 
of $\cF$: In other words, the category $\cC$ must be big enough to 
characterize $\cF$. 

\begin{rem}\label{tsurem}
In the situation of Conjecture \ref{bconj2}, let us take 
a $(Y,\ol{Y})$-triple $(Z,\ol{Z},\cZ)$ over $(S,S,\cS)$ such that 
$(Z,\ol{Z}) \lra (Y,\ol{Y})$ is a strict Zariski open covering and that 
$\cZ$ is formally smooth over $\cS$ on a neighborhood of 
$Z$, and let us define the category 
$\cC_0$ as the category of triples $(Z',\ol{Z}',\cZ')$ over $(Z,\ol{Z},\cZ)$ 
such that $(Z',\ol{Z}') \lra (Z,\ol{Z})$ is a strict open immersion and that 
$\cZ'$ is formally smooth over $\cZ$ on a neighborhood of $Z'$. 
By looking at the definition of `$q$-th rigid cohomology overconvergent 
isocrystal' in \cite[3.3.1]{tsuzuki3}, we see that, in the conjecture given 
in \cite{tsuzuki3}, Tsuzuki allows only the category of the form 
$\cC_0$ (for some $(Z,\ol{Z},\cZ)$) as the category $\cC$ in Conjecture 
\ref{bconj2}. In this sense, Conjecture \ref{bconj2} is slightly weaker than 
the conjecture given in \cite{tsuzuki3}. 
\end{rem} 

\begin{rem} 
Conjecture \ref{bconj2} is proved in the case where there exists a 
Cartesian diagram 
\begin{equation*}
\begin{CD}
\ol{X} @>f>> \ol{Y} @>g>> S \\ 
@VVV @VVV @VVV \\ 
\cP @>{f'}>> \cY @>{g'}>> \cS, 
\end{CD}
\end{equation*} 
where lower horizontal lines are morphisms in $\pLF$ such that 
$f'$ is formally smooth on a neighborhood of $X$ and 
$g'$ is formally smooth on a neighborhood of $Y$ (\cite[Thm 5]{berthelot1}, 
\cite[4.1.4]{tsuzuki3}). 
A weaker result (generic overconvergence) is proved in the case of 
relative dimension $1$ (\cite[4.2.7]{tsuzuki3}). 
In the case $\cE$ is trivial, a result closely related to 
Conjecture \ref{bconj2} is proved when 
$f:X \lra Y$ is an abelian scheme (\cite{etesse}) and when 
$Y$ is a smooth curve (\cite{matsudatrihan}, using \cite{kedlaya0}). 
\end{rem} 

The main purpose of this section is to prove 
Conjecture \ref{bconj2} in the case where a given morphism 
$f:\ol{X} \lra \ol{Y}$ admits a nice log structure such that 
the coefficient $\cE$ extends to $\ol{X}$ logarithmically. 
(For precise statement, see Theorem \ref{main2}, Corollary \ref{main3}.) 
To do this, we should relate overconvergent isocrystals to 
isocrystals on log convergent site. So, we 
give some preliminary results which are needed to 
construct a functor from 
the category of isocrystals on relative log convergent site 
to the category of overconvergent isocrystals. \par 
In the following, for a sheaf of rings $\cA$ on a site, we denote 
the category of $\cA$-modules by $\Mod(\cA)$ and the category of 
coherent $\cA$-modules by $\Coh(\cA)$. \par 
First let us assume given a diagram 
\begin{equation}\label{dag-diag1}
\begin{CD}
(\ol{X},M_{\ol{X}}) @>i>> (\cP,M_{\cP}) \\ 
@VVV @VVV \\ 
\ol{X} @>>> \cP_0, 
\end{CD} 
\end{equation}
where the notations are as follows: 
$(\ol{X},M_{\ol{X}})$ is a fine log $B$-scheme, 
$(\cP,M_{\cP})$ is a $p$-adic fine log formal $\cB$-scheme and 
$\cP_0$ is a $p$-adic formal $\cB$-scheme. The horizontal arrows are 
closed immersions and 
the left vertical arrow is the `forgetting log' morphism. 
(Here we regard (formal) schemes naturally as fine log (formal) schemes 
with trivial log structures.) Moreover, 
let us assume given an open immersion $j_X:X \hra \ol{X}$ satisfying 
$j_X(X) \subseteq (\cP,M_{\cP})_{\triv}$ 
(this implies $j_X(X) \subseteq (\ol{X},M_{\ol{X}})_{\triv}$) 
such that 
the morphism $\cP \lra \allowbreak \cP_0$ in the diagram \eqref{dag-diag1} is 
formally etale on a neighborhood of $X$. Then we have a commutative diagram 
\begin{equation}\label{dag-diag2} 
\begin{CD}
]X[^{\log}_{\cP} @>{j_X^{\log}}>> ]\ol{X}[^{\log}_{\cP} \\ 
@V{\psi_X}VV @V{\varphi_X}VV \\ 
]X[_{\cP_0} @>{j_X}>>]\ol{X}[_{\cP_0}
\end{CD}
\end{equation} 
induced by \eqref{dag-diag1} and $j_X$. \par 
For the moment, we assume moreover that the closed immersion $i$
admits a factorization 
$$ (*): \qquad (\ol{X},M_{\ol{X}}) \hra (\cP',M_{\cP'}) \lra (\cP,M_{\cP}), $$
where the first map is an exact closed immersion and 
the second map is formally log etale morphism of $p$-adic 
fine log formal $\cB$-schemes. 
Let us recall the following: 

\begin{lem}\label{daglem1}
Under the above assumption, 
$\psi_X$ is an isomorphism and 
$\varphi_X$ induces an isomorphism between some strict neighborhood of 
$]X[^{\log}_{\cP}$ in $]\ol{X}[^{\log}_{\cP}$ and some strict neighborhood 
of $]X[_{\cP_0}$ in $]\ol{X}[_{\cP_0}$. 
\end{lem} 

\begin{pf} 
Since $\cP \lra \cP_0$ is formally etale on a neighborhood of $X$, 
there exists an isomorphism of a strict neighborhood of $]X[_{\cP}$ in 
$]\ol{X}[_{\cP}$ and a strict neighborhood of $]X[_{\cP_0}$ in 
$]\ol{X}[_{\cP_0}$ by strong fibration theorem. So we may reduce to the 
case $\cP=\cP_0$ to prove the lemma. 
In this case, 
the upper horizontal line of \eqref{dag-diag2} 
is rewritten as $]X[_{\cP'} \lra ]\ol{X}[_{\cP'}$. 
Since the log structures on 
$(\cP',M_{\cP'})$ and $(\cP,M_{\cP})$ are trivial on a neighborhood of 
$X$, the map $\cP' \lra \cP$ is formally etale on a neighborhood of $X$. 
So $\psi_X$ is an isomorphism by weak fibration theorem and $\varphi_X$ 
induces 
an isomorphism of some strict neighborhood of 
$]X[_{\cP'}$ in $]\ol{X}[_{\cP'}$ and some strict neighborhood 
of $]X[_{\cP}$ in $]\ol{X}[_{\cP}$, by strong fibration theorem. 
\end{pf} 

For two strict neighborhoods $V \subseteq W$ of 
$]X[_{\cP_0}$ in $]\ol{X}[_{\cP_0}$ (resp. $]X[^{\log}_{\cP}$ in 
$]\ol{X}[^{\log}_{\cP}$), we denote the admissible open immersion 
$V \hra W$ by $\alpha_{VW}$ (resp. $\alpha^{\log}_{VW}$) 
and we denote $\alpha_{V]\ol{X}[_{\cP_0}}$ 
(resp. $\alpha^{\log}_{V]\ol{X}[^{\log}_{\cP}}$) simply by $\alpha_V$ 
(resp. $\alpha^{\log}_V$). 
Then the admissible open immersions $j_X, j_X^{\log}$ induce the 
functors $j_X^{\dag}, j_X^{\log,\dag}$ 
which are defined by $j_X^{\dag} := \varinjlim_{V} \alpha_{V,*}\alpha_V^{-1}, 
j_X^{\dag} := \varinjlim_{V} \alpha^{\log}_{V,*}\alpha^{\log,-1}_V$. 
It is known that 
$j_X^{\dag}$ (resp. $j_X^{\log,\dag}$) sends $\cO_{]\ol{X}[_{\cP_0}}$-modules 
(resp. $\cO_{]\ol{X}[^{\log}_{\cP}}$-modules) to 
$j_X^{\dag}\cO_{]\ol{X}[_{\cP_0}}$-modules 
(resp. $j_X^{\dag}\cO_{]\ol{X}[^{\log}_{\cP}}$-modules) and 
coherent $\cO_{]\ol{X}[_{\cP_0}}$-modules 
(resp. coherent $\cO_{]\ol{X}[^{\log}_{\cP}}$-modules) to 
coherent $j_X^{\dag}\cO_{]\ol{X}[_{\cP_0}}$-modules 
(resp. coherent $j_X^{\dag}\cO_{]\ol{X}[^{\log}_{\cP}}$-modules). 
Then we have the following: 

\begin{lem}\label{dagdagdaglem} 
The functor 
$j_X^{\d} \circ \varphi_{X,*}: \Mod(\cO_{]\ol{X}[^{\log}_{\cP}}) \lra 
\Mod(j_X^{\dag}\cO_{]\ol{X}[_{\cP_0}})$ is exact and 
send coherent modules to coherent modules. 
We have the equalities of functors 
$$j_X^{\d} \circ \varphi_{X,*} = R(j_X^{\d} \circ \varphi_{X,*}) = 
j_X^{\d} \circ 
R\varphi_{X,*}: \Mod(\cO_{]\ol{X}[^{\log}_{\cP}}) \lra 
D^+(\Mod(j_X^{\d}\cO_{]\ol{X}[_{\cP_0}})), $$
$$j_X^{\dag} = j_X^{\dag} \circ \varphi_{X,*} 
\circ \varphi_X^*: \Mod(\cO_{]\ol{X}[_{\cP_0}}) \lra 
\Mod(j_X^{\d} \cO_{]\ol{X}[_{\cP_0}}), $$
$$ \varphi_{X,*} \circ j^{\log,\d}_X = 
j_X^{\d} \circ \varphi_{X,*}: 
\Coh(\cO_{]\ol{X}[^{\log}_{\cP}}) \lra 
\Coh(j_X^{\d}\cO_{]\ol{X}[_{\cP_0}}). $$
\end{lem} 

\begin{pf} 
It is easy to see that the functor $j_X^{\d} \circ \varphi_{X,*}$ 
is left exact. 
To show the right exactness, let us take a strict neighborhood $V_0$ of 
$]X[_{\cP_0}$ in $]\ol{X}[_{\cP_0}$ such that $\varphi^{-1}_X(V_0)$ is 
isomorphic to $V_0$ via $\varphi_X$. Then, for any surjection 
$\cF_1 \lra \cF_2$ of $\cO_{]\ol{X}[^{\log}_{\cP}}$-modules, 
$(\varphi_{X,*}\cF_1)|_{V_0} \lra (\varphi_{X,*}\cF_2)|_{V_0}$ is 
surjective. Now recall that the functor 
$\varinjlim_{V}\alpha_{V,*}\alpha_{VV_0}^{-1}$ (where $V$ runs through 
strict neighborhoods contained in $V_0$) is exact 
(\cite[2.1.3]{berthelot2}). By 
applying this functor, we obtain the surjection 
$(j_X^{\d} \circ \varphi_{X,*})\cF_1 \lra 
(j_X^{\d} \circ \varphi_{X,*})\cF_2$. 
So $j_X^{\d} \circ \varphi_{X,*}$ is exact. 
Moreover, for a coherent $\cO_{]\ol{X}[^{\log}_{\cP}}$-module $\cF$, 
$(\varphi_{X,*}\cF)|_{V_0}$ is a coherent $\cO_{V_0}$-modules and so 
$(j_X^{\d} \circ \varphi_{X,*})\cF = 
\varinjlim_{V}\alpha_{V,*}\alpha_{VV_0}^{-1} ((\varphi_{X,*}\cF)|_{V_0})$ 
is a coherent $j_X^{\dag}\cO_{]\ol{X}[^{\log}_{\cP}}$-module. So 
$j_X^{\d} 
\circ \varphi_{X,*}$ sends coherent modules to coherent modules. \par 
The first 
equalities in the lemma follows from the exactness of $j_X^{\dag}$ and 
$j^{\dag}_X \circ \varphi_{X,*}$. 
The second equality of the lemma is proved as 
\begin{align*}
j_X^{\dag} \circ \varphi_{X,*} 
\circ \varphi_X^* & = 
\varinjlim_V \alpha_{V,*}\alpha^{-1}_V\varphi_{X,*}\varphi^*_X \\
& = 
\varinjlim_V \alpha_{V,*}(\varphi_X|_{\varphi_X^{-1}(V)})_* 
\alpha^{\log,-1}_{\varphi_X^{-1}(V)}\varphi_X^* \\ 
& = 
\varinjlim_V \alpha_{V,*}(\varphi_X|_{\varphi_X^{-1}(V)})_* 
(\varphi_X|_{\varphi_X^{-1}(V)})^*\alpha_V^{-1} 
= \varinjlim_V \alpha_{V,*}\alpha_V^{-1} = j_X^{\d}, 
\end{align*}
where the fourth equality follows from Lemma \ref{daglem1}. \par 
Finally we prove the last equality of the lemma. The map of functors 
$j^{\d}_X \circ \varphi_{X,*} \lra \varphi_{X,*} \circ j^{\log,\d}_X$ 
is defined by 
\begin{align*}
j^{\d}_X \circ \varphi_{X,*} & = 
\varinjlim_V\alpha_{V,*}\alpha_V^{-1}\varphi_{X,*} = 
\varinjlim_V\alpha_{V,*}(\varphi_X|_{\varphi_X^{-1}(V)})_*
\alpha^{\log,-1}_{\varphi_X^{-1}(V)} \\ 
& = 
\varinjlim_V\varphi_{X,*}\alpha^{\log}_{\varphi_X^{-1}(V),*}
\alpha^{\log,-1}_{\varphi_X^{-1}(V)} 
\lra 
\varphi_{X,*}\varinjlim_V \alpha^{\log}_{\varphi^{-1}(V),*}
\alpha^{\log,-1}_{\varphi^{-1}(V)} = 
\varphi_{X,*} \circ j^{\log,\d}_X. 
\end{align*}
We prove the exactness of 
$\varphi_{X,*} \circ j^{\log,\d}_X$ on $\Coh(\cO_{]\ol{X}[^{\log}_{\cP}})$. 
Since $\varphi_X$ is isomorphic on $\varphi_X^{-1}(V)$, 
the restriction functor 
$\varphi_X^*: \Coh(j_X^{\d}\cO_{]\ol{X}[_{\cP}}) \lra 
\Coh(j_X^{\log,\d}\cO_{]\ol{X}[^{\log}_{\cP}})$ is an equivalence of 
categories by \cite[2.1.4, 2.1.10]{berthelot2}. Moreover, 
since $\cP' \lra \cP$ (where $\cP'$ is as in the diagram $(*)$ before Lemma 
\ref{daglem1}) is formally etale on a neighborhood of $X$, 
we have, for any $E \in \Coh(j_X^{\d}\cO_{]\ol{X}[_{\cP}})$, 
the equality $R\varphi_{X,*}\varphi_X^*E = E$ by 
a special case of \cite[8.3.5]{chts}. Hence, for any 
$F \in \Coh(j_X^{\log,\d}\cO_{]\ol{X}[^{\log}_{\cP}})$, we have 
$R\varphi_{X,*}F = \varphi_{X,*}F$. Since $j_X^{\log,\d}$ is exact, 
we can conclude that the functor $\varphi_{X,*}\circ j_X^{\log,\d}$ 
is exact on $\Coh(\cO_{]\ol{X}[^{\log}_{\cP}})$. So, to prove that 
the map 
$j^{\d}_X \circ \varphi_{X,*} \lra \varphi_{X,*} \circ j^{\log,\d}_X$ 
is isomorphic on $\Coh(\cO_{]\ol{X}[^{\log}_{\cP}})$, it suffices to check 
the isomorphism 
$$j^{\d}_X \circ \varphi_{X,*}\cO_{]\ol{X}[^{\log}_{\cP}} 
\os{\simeq}{\lra} 
\varphi_{X,*} \circ j^{\log,\d}_X\cO_{]\ol{X}[^{\log}_{\cP}}.$$ 
It is true since the left hand side is calculated as 
$$ 
j^{\d}_X \circ \varphi_{X,*}\cO_{]\ol{X}[^{\log}_{\cP}} = 
\varinjlim_V\alpha_{V,*}\cO_V = j^{\d}_X\cO_{]\ol{X}[_{\cP_0}} $$ 
and the right hand side is calculated as 
$$ 
\varphi_{X,*} \circ j^{\log,\d}_X\cO_{]\ol{X}[^{\log}_{\cP}} = 
\varphi_{X,*}\varphi_X^*j_X^{\d}\cO_{]\ol{X}[_{\cP_0}}= 
j_X^{\d}\cO_{]\ol{X}[_{\cP_0}}. 
$$ 
So we are done. 
\end{pf} 

\begin{defn}
With the above notation, we define the functor 
$j_X^{\dd}: \Mod(\cO_{]\ol{X}[^{\log}_{\cP}}) \lra 
\Mod(j_X^{\d}\cO_{]\ol{X}[_{\cP_0}})$ 
by $j_X^{\dd} := j_X^{\d} \circ \varphi_{X,*}$. 
\end{defn}

Next let us consider the case where the closed immersion 
$i:(\ol{X},M_{\ol{X}}) \hra (\cP,M_{\cP})$ does not necessarily 
admit a factorization $(*)$. Also in this case, we can define the 
functor 
$j_X^{\dd}: \Mod(\cO_{]\ol{X}[^{\log}_{\cP}}) \lra 
\Mod(j_X^{\d}\cO_{]\ol{X}[_{\cP_0}})$ 
by $j_X^{\dd} := j_X^{\d} \circ \varphi_{X,*}$. 
Then we have the following: 

\begin{prop}\label{dagdagprop}
With the above assumption, the functor $j_X^{\dd}$ sends 
coherent modules to coherent modules and we have 
$Rj_X^{\dd}E = j_X^{\dd}E$ for a coherent 
$\cO_{]\ol{X}[^{\log}_{\cP}}$-module $E$. 
\end{prop} 

\begin{pf} 
First, let us factorize the map $\varphi_X: ]\ol{X}[^{\log}_{\cP}
 \lra ]\ol{X}[_{\cP_0}$ as 
$]\ol{X}[^{\log}_{\cP} \os{\varphi'_X}{\lra} ]\ol{X}[_{\cP} 
\os{\varphi''_X}{\lra} ]\ol{X}[_{\cP_0}$. Then 
there exists a strict neighborhood $V_0$ of $]X[_{\cP_0}$ in 
$]\ol{X}[_{\cP_0}$ such that ${\varphi''_X}^{-1}(V_0)$ is isomorphic to 
$V_0$ via $\varphi''_X$. Then we have 
\begin{align*}
j_X^{\dd} 
& = \varinjlim_{V} \alpha_{V]\ol{X}[_{\cP_0,*}} \alpha^{-1}_{VV_0} 
\alpha^{-1}_{V_0]\ol{X}[_{\cP_0}} \varphi''_{X,*} \varphi'_{X,*} \\ 
& = \varinjlim_{V} \alpha_{V]\ol{X}[_{\cP_0,*}} \alpha^*_{VV_0} 
(\varphi''_X\vert_{{\varphi''_X}^{-1}(V_0)})_* 
\alpha^*_{{\varphi''_X}^{-1}(V_0)]\ol{X}[_{\cP}} \varphi'_{X,*}. 
\end{align*}
Moreover, $\varphi''_X\vert_{{\varphi''_X}^{-1}(V_0)}$ is an isomorphism, 
$\varinjlim_{V} \alpha_{V]\ol{X}[_{\cP_0},*} \alpha^{-1}_{VV_0}$ 
is an exact functor sending coherent modules to coherent modules and 
the same is true if we replace $V_0$ by a smaller strict neighborhood. 
So it suffices to prove that, for some strict neighborhood 
$V$ of $]X[_{\cP}$ in $]\ol{X}[_{\cP}$ contained in ${\varphi''}^{-1}_X(V_0)$, 
the functor 
$\alpha^{-1}_{V]\ol{X}[_{\cP}} \varphi'_{X,*}$ sends coherent modules to 
coherent modules and that 
$\alpha^{-1}_{V]\ol{X}[_{\cP}} R^q \varphi'_{X,*} E = 0$ holds 
for any coherent $\cO_{]\ol{X}[^{\log}_{\cP}}$-module $E$ and $q \geq 1$. \par 
Next, let us take a Cartesian diagram 
\begin{equation*}
\begin{CD}
(\ol{X}^{(\b)},M_{\ol{X}^{(\b)}}) @>{i^{(\b)}}>> 
(\cP^{(\b)},M_{\cP^{(\b)}}) \\ 
@VVV @VgVV \\ 
(\ol{X},M_{\ol{X}}) @>i>> (\cP,M_{\cP}), 
\end{CD}
\end{equation*}
where $g$ is a strict formally etale hypercovering such that 
each $i^{(n)}$ admits a factorization like $(*)$ and that 
each $\cP^{(n)}$ is affine over $\cP$. Then $g$ induces a morphism 
$\ti{g}:]\ol{X}^{(\b)}[_{\cP^{(\b)}}^{\log} \lra 
]\ol{X}[_{\cP}^{\log}$. Then, for a coherent 
$\cO_{]\ol{X}[_{\cP}^{\log}}$-module $E$, we have 
$\alpha^{-1}_{V]\ol{X}[_{\cP}} R^q \varphi'_{X,*}E = 
\alpha^{-1}_{V]\ol{X}[_{\cP}} R^q (\varphi'_{X} \circ \ti{g})_*(\ti{g}^*E) \, 
(q \geq 0)$ by Lemma \ref{wdlem1}. 
Let us consider the following factorization of 
$\varphi'_{X} \circ \ti{g}$: 
$$ ]\ol{X}^{(\b)}[_{\cP^{(\b)}}^{\log} \os{\gamma^{(\b)}}{\lra} 
]\ol{X}^{(\b)}[_{\cP^{(\b)}} \os{\beta^{(\b)}}{\lra} ]\ol{X}[_{\cP}. $$
Then we have 
$\alpha^{-1}_{V]\ol{X}[_{\cP}} R^q(\varphi'_{X} \circ \ti{g})_*(\ti{g}^*E) 
= \alpha^{-1}_{V]\ol{X}[_{\cP}}R^q(\beta^{(\b)} \circ \gamma^{(\b)})_*
(\ti{g}^*E)$ for $q \geq 0$. \par 
Let us fix $q \geq 0$ and put $q':=\max (2,q(q+1)/2)$. 
By Lemma \ref{daglem1}, there exist a strict neighborhood 
$V^{(n)}$ of $]X^{(n)}[_{\cP^{(n)}}$ in $]\ol{X}^{(n)}[_{\cP^{(n)}}$ 
(here $X^{(n)}:= X \times_{\ol{X}} \ol{X}^{(n)}$) 
for each $n \leq q'$ 
such that $\gamma^{(n),-1}(V^{(n)})$ is isomorphic to $V^{(n)}$ via 
$\gamma^{(n)}$. Since a standard strict neighborhood of 
$]X^{(n)}[_{\cP^{(n)}}$ in $]\ol{X}^{(n)}[_{\cP^{(n)}}$ is the 
inverse image of a standard strict neighborhood of 
$]X[_{\cP}$ in $]\ol{X}[_{\cP}$, we can assume by shrinking $V^{(n)}$ that 
$\beta^{(n),-1}(V) = V^{(n)}$ holds (for $n\leq q'$) for some 
strict neighborhood $V$ of $]X[_{\cP}$ in $]\ol{X}[_{\cP}$. 
Then we put $V^{(\b)}:= \beta^{(\b),-1}(V), 
\ti{V}^{(\b)} := \gamma^{(\b),-1}(V^{(\b)})$. 
(Then we have the isomorphism $\gamma^{(n)}|_{\ti{V}^{(n)}}: 
\ti{V}^{(n)} \os{\sim}{\lra} V^{(n)}$ for 
$n \leq q'$.) If we put $F^{(\b)}:= 
\alpha^{-1}_{\ti{V}^{(\b)}]\ol{X}[^{\log}_{\cP}}\ti{g}^*E$, 
we have 
$\alpha^{-1}_{V]\ol{X}[_{\cP}}R^q(\beta^{(\b)} 
\circ \gamma^{(\b)})_*(\ti{g}^*E) 
= R^q((\beta^{(\b)}\vert_{V^{(\b)}}) \circ 
(\gamma^{(\b)}\vert_{\ti{V}^{(\b)}}))_* F^{(\b)}$. \par 
Let us consider the morphism $\beta^{(\b)}|_{V^{(\b)}}: V^{(\b)} \lra V$: 
For any sufficiently small affinoid admissible open set 
$U := \Spm (\Q \otimes_{\Z} A)\hra V$, 
the morphism $V^{(\b)} \times_V U \lra U$ induced by 
$\beta^{(\b)}|_{V^{(\b)}}$ is the one induced by some 
affine formally etale hypercovering of $\Spf A$. 
So we have the descent property and the cohomological descent property 
for coherent modules for the morphism $\beta |_{V^{(\b)}}$ 
(\cite[7.3.3]{chts}). Since 
$\gamma^{(n)}|_{\ti{V}^{(n)}}$ is an isomorphism for $n \leq q'$, 
$\{ (\gamma^{(n)}|_{\ti{V}^{(n)}})_*F^{(n)} \}_{n \leq q'}$ forms a 
compatible family of coherent $\cO_{V^{(n)}}$-modules. So there exists 
a coherent $\cO_V$-module $G$ and a homomorphism 
$(\beta^{(\b)}|_{V^{(\b)}})^* G \lra 
(\gamma^{(\b)}|_{\ti{V}^{(\b)}})_*F^{(\b)}$ which is an isomorphism for 
$\b \leq q'$. \par 
Since we have $R^t(\gamma^{(u)}|_{\ti{V}^{(u)}})_*F^{(u)}=0$ for 
$u \leq q'$ and $t>0$, we have the equality 
$R^v(\beta^{(u)}|_{V^{(u)}})_*R^t(\gamma^{(u)}|_{\ti{V}^{(u)}})_*F^{(u)}=0$ 
for $u \leq q', v \geq 0, t>0$. By the spectral sequence 
$$ 
E_1^{u,v} := R^v(\beta^{(u)}|_{V^{(u)}})_*
R^t(\gamma^{(u)}|_{\ti{V}^{(u)}})_*F^{(u)} \,\Longrightarrow\, 
R^{u+v}(\beta^{(\b)}|_{V^{(\b)}})_* R^t(\gamma^{(\b)}|_{\ti{V}^{(\b)}})_*
F^{(\b)}, 
$$ 
we obtain the equality 
$R^s(\beta^{(\b)}|_{V^{(\b)}})_* R^t(\gamma^{(\b)}|_{\ti{V}^{(\b)}})_*
F^{(\b)} = 0$ for $s \leq q', t>0$. Then, by the spectral sequence 
$$ 
E_2^{s,t} := 
R^s(\beta^{(\b)}|_{V^{(\b)}})_* R^t(\gamma^{(\b)}|_{\ti{V}^{(\b)}})_*
F^{(\b)} 
\,\Longrightarrow\, 
R^{s+t}((\beta^{(\b)}|_{V^{(\b)}}) \circ (\gamma^{(\b)}|_{\ti{V}^{(\b)}}))_*
F^{(\b)}, 
$$ 
we obtain the equality 
$R^q((\beta^{(\b)}|_{V^{(\b)}}) \circ (\gamma^{(\b)}|_{\ti{V}^{(\b)}}))_*
F^{(\b)} = 
R^q(\beta^{(\b)}|_{V^{(\b)}})_* (\gamma^{(\b)}|_{\ti{V}^{(\b)}})_*
F^{(\b)}.$ Now let us consider the commutative diagram 
\begin{equation}\label{spdiag}
\begin{CD}
E_1^{u,v} := 
R^v (\beta^{(u)}|_{V^{(u)}})_* 
(\gamma^{(u)}|_{\ti{V}^{(u)}})_*F^{(u)} @. \,\Longrightarrow\, @. 
R^{u+v}(\beta^{(\b)}|_{V^{(\b)}})_* (\gamma^{(\b)}|_{\ti{V}^{(\b)}})_*
F^{(\b)} \\ 
@AAA @. @AAA \\ 
E_1^{u,v} := 
R^v (\beta^{(u)}|_{V^{(u)}})_* 
(\beta^{(u)}|_{V^{(u)}})^*G @. \,\Longrightarrow\, @. 
R^{u+v}(\beta^{(\b)}|_{V^{(\b)}})_* (\beta^{(\b)}|_{V^{(\b)}})^* G. 
\end{CD}
\end{equation}
Since the left vertical arrow is an isomorphism for $u \leq q'$, 
the morphism 
$$ R^{q}(\beta^{(\b)}|_{V^{(\b)}})_* (\beta^{(\b)}|_{V^{(\b)}})^* G \lra 
R^{q}(\beta^{(\b)}|_{V^{(\b)}})_* (\gamma^{(\b)}|_{\ti{V}^{(\b)}})_*
F^{(\b)}$$ 
induced by the right vertical arrow in \eqref{spdiag} 
is an isomorphism and by the cohomological descent property 
for $\beta^{(\b)} |_{V^{(\b)}}$, we have 
$(\beta^{(\b)}|_{V^{(\b)}})_* (\beta^{(\b)}|_{V^{(\b)}})^* G = G$, 
$R^q(\beta^{(\b)}|_{V^{(\b)}})_* (\beta^{(\b)}|_{V^{(\b)}})^* G \allowbreak 
= 0\,(q \geq 1)$. So we have 
$(\beta^{(\b)}|_{V^{(\b)}})_* (\gamma^{(\b)}|_{\ti{V}^{(\b)}})_*
F^{(\b)} =G$ and 
$R^q(\beta^{(\b)}|_{V^{(\b)}})_* (\gamma^{(\b)}|_{\ti{V}^{(\b)}})_*
F^{(\b)} \allowbreak =\allowbreak 
0\,(q \geq 1)$. These imply the equalities 
$\alpha^{-1}_{V]\ol{X}[_{\cP}}(\beta^{(\b)} \circ 
\gamma^{(\b)})_*(\ti{g}^*E) =G$, 
$\alpha^{-1}_{V]\ol{X}[_{\cP}}R^q(\beta^{(\b)} 
\circ \gamma^{(\b)})_*(\ti{g}^*E) 
\allowbreak =\allowbreak 0 \,(q \geq 1)$. So we are done. 
\end{pf}

Now assume that we are given a diagram 
\begin{equation}\label{ddag-diag}
(\ol{X},M_{\ol{X}}) \os{f}{\lra} (\ol{Y},M_{\ol{Y}}) \os{\iota}{\hra} 
(\cY,M_{\cY}), 
\end{equation}
where $f$ is a morphism in $\L$, $(\cY,M_{\cY})$ is an object in 
$\pLF$ and $\iota$ is a closed immersion. Assume also that we are given 
open immersions $X \subseteq \ol{X}, Y \subseteq \ol{Y}$ satisfying  
$X \subseteq (\ol{X},M_{\ol{X}})_{\triv}$, 
$Y \subseteq (\cY,M_{\cY})_{\triv}$ and $f(X) \subseteq Y$. \par 
Assume for the moment that 
the above diagram fits into the diagram 
\begin{equation}\label{ddag-diag2}
\begin{CD}
(\ol{X},M_{\ol{X}}) @>{i}>> (\cP,M_{\cP}) \\ 
@VfVV @VgVV \\ 
(\ol{Y},M_{\ol{Y}}) @>{\iota}>> (\cY,M_{\cY}), 
\end{CD}
\end{equation}
where $(\cP,M_{\cP})$ is a $p$-adic fine log formal $\cB$-scheme, 
$g$ is a formally log smooth morphism and 
$i$ is a closed immersion satisfying $X \subseteq (\cP,M_{\cP})_{\triv}$. 
For $n \in \N$, let $(\cP(n),M_{\cP(n)})$ be the $(n+1)$-fold 
fiber product of $(\cP,M_{\cP})$ over $(\cY,M_{\cY})$ and let 
$\cP_0(n)$ be the $(n+1)$-fold fiber product of $\cP$ over $\cY$. 
Then we have the exact functors 
$$ j^{\dd}_X: \Coh(\cO_{]\ol{X}[^{\log}_{\cP(n)}}) \lra 
\Coh(j^{\dag}_X\cO_{]\ol{X}[_{\cP_0(n)}}) $$ 
which are compatible with pull-backs by the projections 
$$ p_i:\,]\ol{X}[^{\log}_{\cP(1)} \lra \,]\ol{X}[^{\log}_{\cP(0)} \,\,\, 
(i=1,2), 
\qquad 
p_{ij}:\,]\ol{X}[^{\log}_{\cP(2)} \lra \,]\ol{X}[^{\log}_{\cP(1)} \,\,\, 
(1 \leq i < j \leq 3) $$
and the diagonal $\Delta:\,]\ol{X}[^{\log}_{\cP(0)} \hra 
]\ol{X}[^{\log}_{\cP(1)}$. So they naturally induce the functor 
$$ 
I_{\conv}((\ol{X}/\cY)^{\log}) \os{\sim}{\lra} 
\Strat''((\ol{X} \hra \cP/\cY)^{\log}) \lra 
I^{\dag}((X,\ol{X})/\cY_K, \cP) = I^{\dag}((X,\ol{X})/\cY_K), $$
which we denote also by $j^{\dd}_X$. \par 
In the case where 
there does not necessarily exist 
the diagram \eqref{ddag-diag2}, we can construct the diagram 
\begin{equation*}
(\ol{X},M_{\ol{X}}) \os{h}{\lla} (\ol{X}^{(\b)},M_{\ol{X}^{(\b)}}) 
\os{i^{(\b)}}{\hra} (\cP^{(\b)}, M_{\cP^{(\b)}}), 
\end{equation*}
where $h$ is a strict etale hypercovering, 
$i^{(\b)}$ is a closed immersion over $\iota$ such that 
each $(\cP^{(n)}, M_{\cP^{(n)}})$ is formally log 
smooth over $(\cY,M_{\cY})$ and 
$X^{(\b)} \subseteq (\cP^{(\b)},M_{\cP^{(\b)}})_{\triv}$ holds, where 
$X^{(\b)} := X \times_{\ol{X}} \ol{X}^{(\b)}$. 
If we denote the category of descent data with respect to 
$I_{\conv}((\ol{X}^{(n)}/\cY)^{\log})$ 
(resp. $I^{\dag}((X^{(n)},\ol{X}^{(n)})/\cY)$) for $n=0,1,2$ simply by 
$\Delta_{\conv}$ (resp. $\Delta^{\dag}$), the functor $j^{\dd}_X$ in the 
previous paragraph naturally induces the functor $\Delta_{\conv} \lra 
\Delta^{\dag}$. So we have the canonical functor 
$$ 
I_{\conv}((\ol{X}/\cY)^{\log}) \os{\sim}{\lra} 
\Delta_{\conv} \lra \Delta^{\dag} \os{\sim}{\lla} 
I^{\dag}((X,\ol{X})/\cY), $$
which we denote by $j^{\dd}_X$. It is easy to check that 
the definition of the functor $j^{\dd}_X$ is independent of any choice which 
we have made. 

\begin{rem}\label{drcpx}
Let the situation be as above and let $\cE$ be an object in 
$I_{\conv}((\ol{X}/\cY)^{\log})$. Then, associated to $\cE$, we have a 
coherent $\cO_{]\ol{X}^{(\b)}[^{\log}_{\cP^{(\b)}}}$-module $\cE^{(\b)}$ 
with integrable log connection 
$$\nabla: \cE^{(\b)} \lra \cE^{(\b)} 
\otimes_{\cO_{]\ol{X}^{(\b)}[^{\log}_{\cP^{(\b)}}}}
\omega^1_{]\ol{X}^{(\b)}[^{\log}_{\cP^{(\b)}}/\cY_K}.$$ 
On the other hand, associated to $j^{\dd}_X\cE$, we have a coherent 
$j^{\d}_X\cO_{]\ol{X}^{(\b)}[_{\cP^{(\b)}}}$-module with integrable 
connection of the form 
$$\nabla': j^{\dd}_X\cE^{(\b)} \lra 
j^{\dd}_X\cE^{(\b)} \otimes_{j^{\d}_X\cO_{]\ol{X}^{(\b)}[_{\cP^{(\b)}}}} 
j^{\d}_X \Omega^1_{]\ol{X}^{(\b)}[_{\cP^{(\b)}}/\cY_K}.$$  
By construction of the functor $j^{\dd}_X$, we can see that 
$j^{\dd}_X(\nabla)=\nabla'$ holds. So we have the canonical isomorphism 
of de Rham complexes 
$j_X^{\dd}\DR(]\ol{X}^{(\b)}[^{\log}_{\cP^{(\b)}}/\cY_K, \cE)
 = 
\DR^{\d}(]\ol{X}[_{\cP}/\cY_K, \allowbreak j^{\dd}_X\cE)$. Note that, since 
each term of 
$\DR(]\ol{X}^{(\b)}[^{\log}_{\cP^{(\b)}}/\cY_K, \cE)$ is coherent, 
we have the quasi-isomorphism 
$j_X^{\dd}\DR(]\ol{X}^{(\b)}[^{\log}_{\cP^{(\b)}}/ \allowbreak \cY_K, \cE)
= Rj_X^{\dd}\DR(]\ol{X}^{(\b)}[^{\log}_{\cP^{(\b)}}/\cY_K, \cE)$ by 
Proposition \ref{dagdagprop}. 
\end{rem}

Now we prove our first main result in this section: 

\begin{thm}\label{main1}
Assume we are given a diagram 
\begin{equation}\label{maindiag}
(\ol{X},M_{\ol{X}}) \os{f}{\lra} (\ol{Y},M_{\ol{Y}}) 
\os{\iota}{\hra} (\cY,M_{\cY}) 
\end{equation}
and open immersions $j_X:X \hra \ol{X}, j_Y:Y \hra \ol{Y}$ 
satisfying $X \subseteq (\ol{X}, M_{\ol{X}})_{\triv}, 
Y \subseteq (\cY,M_{\cY})_{\triv}$ and $f^{-1}(Y)=X$, where 
$f$ is a proper log smooth integral morphism in $\LB$ having 
log smooth parameter and $\iota$ is a closed immersion in 
$\pLF$. Then, 
for a locally free isocrystal $\cE$ on $(\ol{X}/\cY)^{\log}_{\conv}$ and 
$q \geq 0$,  
we have the isomorphism 
$$ j^{\dd}_Y R^qf_{\ol{X}/\cY, \an *}\cE = 
R^qf_{(X,\ol{X})/\cY,\rig *}j_X^{\dd}\cE $$ 
of $j_Y^{\d}\cO_{]\ol{Y}[_{\cY}}$-modules. In particular, 
$R^qf_{(X,\ol{X})/\cY,\rig *}j_X^{\dd}\cE$ is a coherent 
$j_Y^{\d}\cO_{]\ol{Y}[_{\cY}}$-module. 
\end{thm} 

\begin{pf}
First we construct the map 
$j^{\dd}_Y R^qf_{\ol{X}/\cY, \an *}\cE \lra 
R^qf_{(X,\ol{X})/\cY,\rig *}j_X^{\dd}\cE$. Let us take 
an embedding system 
$$ (\ol{X},M_{\ol{X}}) \lla (\ol{X}^{(\b)},M_{\ol{X}^{(\b)}}) \hra 
(\cP^{(\b)},M_{\cP^{(\b)}}) $$ 
satisfying $X^{(\b)} \subseteq (\cP^{(\b)},M_{\cP^{(\b)}})_{\triv}$, 
where $X^{(\b)} := X \times_{\ol{X}} \ol{X}^{(\b)}$. 
Then we have the diagram of rigid analytic spaces 
\begin{equation*}
\begin{CD}
]\ol{X}^{(\b)}[^{\log}_{\cP^{(\b)}} @>{h^{\log}}>> ]\ol{Y}[^{\log}_{\cY} \\ 
@VVV @VVV \\ 
]\ol{X}^{(\b)}[_{\cP^{(\b)}} @>h>> ]\ol{Y}[_{\cY}. 
\end{CD}
\end{equation*}
Then we have 
{\small{
\begin{align*} 
(\spadesuit) \,\,\,\, Rj^{\dd}_YRf_{\ol{X}/\cY,\an *} \cE 
& = 
Rj^{\dd}_YRh^{\log}_*\DR(]\ol{X}^{(\b)}[^{\log}_{\cP^{(\b)}}/\cY_K, \cE) 
\\ 
& \lra 
Rh_*Rj^{\dd}_X \DR(]\ol{X}^{(\b)}[^{\log}_{\cP^{(\b)}}/\cY_K, \cE) \,\,\,\, 
\text{(induced by $j^{\dd}_Y \circ h^{\log}_* \ra 
h_* \circ j^{\dd}_X$)} \\ 
& \os{\sim}{\lla} 
Rh_*j^{\dd}_X \DR(]\ol{X}^{(\b)}[^{\log}_{\cP^{(\b)}}/\cY_K, \cE) \,\,\,\, 
\text{(Remark \ref{drcpx})} \\ 
& = Rh_*\DR^{\d}(]\ol{X}^{(\b)}[_{\cP^{(\b)}}/\cY_K, j^{\dd}_X\cE) \,\,\,\, 
\text{(Remark \ref{drcpx})} \\ 
& = Rf_{(X,\ol{X})/\cY,\rig *}j^{\dd}_X\cE, 
\end{align*}}}
and by taking 
the $q$-th cohomology sheaf of the above
 map, we obtain the map 
\begin{equation}\label{map}
j^{\dd}_Y R^qf_{\ol{X}/\cY, \an *}\cE \lra 
R^qf_{(X,\ol{X})/\cY,\rig *}j_X^{\dd}\cE. 
\end{equation}
(We have $\cH^q(Rj^{\dd}_YRf_{\ol{X}/\cY,\an *} \cE) = 
j^{\dd}_YR^qf_{\ol{X}/\cY,\an *} \cE$ by Corollary \ref{coherence} 
and Proposition \ref{dagdagprop}.) 
We will prove that the map \eqref{map}
is an isomorphism. \par 
Let us take a strict formally etale hypercovering 
$\epsilon: (\cY^{\lr},M_{\cY^{\lr}}) \lra (\cY,M_{\cY})$ such that, 
if we put $(Y^{\lr},M_{Y^{\lr}}) := (Y,M_Y) \times_{(\cY,M_{\cY})} 
(\cY^{\lr},M_{\cY^{\lr}})$, the closed immersion 
$(\ol{Y}^{\lmr},M_{\ol{Y}^{\lmr}}) \hra (\cY^{\lmr},M_{\cY^{\lmr}})$ admits 
a factorization 
$$ (\ol{Y}^{\lmr},M_{\ol{Y}^{\lmr}}) \hra ({\cY'}^{\lmr}, M_{{\cY'}^{\lmr}})
\lra (\cY^{\lmr},M_{\cY^{\lmr}}) $$
for each $m$, where the first map is an exact closed immersion and the 
second map is formally log etale. 
Then, since $R^qf_{\ol{X}/\cY,\an *}\cE$ is coherent, 
we have $R\epsilon_*\epsilon^*j^{\dd}_YR^qf_{\ol{X}/\cY,\an *}\cE 
\allowbreak = \allowbreak 
\epsilon_*\epsilon^*j^{\dd}_YR^qf_{\ol{X}/\cY,\an *}\cE = 
j^{\dd}_YR^qf_{\ol{X}/\cY,\an *}\cE$ by \cite[7.3.3]{chts}, where 
we denoted the morphism $]\ol{Y}^{\lr}[_{\cY^{\lr}} \lra ]\ol{Y}[_{\cY}$ 
also by $\epsilon$. 
Moreover, we have the equality of functors 
$j^{\dd}_{Y^{\lr}}\epsilon^* = \epsilon^*j^{\dd}_Y$ for 
coherent modules and the base change isomorphism 
$\epsilon^* R^qf_{\ol{X}/\cY,\an *}\cE = 
R^qf_{\ol{X}^{\lr}/\cY^{\lr},\an *}\cE^{\lr}$, 
where $\cE^{\lr}$ is the restriction of $\cE$ to $\ol{X}^{\lr}/\cY^{\lr}$. 
(This follows from the analytically flat base change theorem, 
Theorem \ref{coherence0}, Remark \ref{rem3-1} and Lemma \ref{2.1.26}.)
So we have the isomorphism 
\begin{equation}\label{maineq1}
\epsilon_*j^{\dd}_{Y^{\lr}}R^qf_{\ol{X}^{\lr}/\cY^{\lr},\an *}\cE^{\lr} = 
j_Y^{\dd} R^qf_{\ol{X}/\cY,\an *}\cE. 
\end{equation}
On the other hand, let 
{\small{
$$ 
(\cP^{(\b) \lr},M_{\cP^{(\b) \lr}}) \hookleftarrow 
(\ol{X}^{(\b) \lr},M_{\ol{X}^{(\b) \lr}}) \lra 
(\ol{X}^{\lr},M_{\ol{X}^{\lr}}) \lra (\ol{Y}^{\lr},M_{\ol{Y}^{\lr}}) \hra 
(\cY^{\lr},M_{\cY^{\lr}}) $$ }}
be the pull-back of the diagram 
$$ 
(\cP^{(\b)},M_{\cP^{(\b)}}) \hookleftarrow
(\ol{X}^{(\b)},M_{\ol{X}^{(\b)}}) \lra 
(\ol{X},M_{\ol{X}}) \lra (\ol{Y},M_{\ol{Y}}) \hra 
(\cY,M_{\cY}) $$ 
by $\epsilon$ and denote the induced morphisms 
$$ ]\ol{X}^{(\b) \lr}[_{\cP^{(\b) \lr}} \lra ]\ol{Y}[_{\cY}, 
\qquad 
]\ol{X}^{(\b)\lr}[_{\cP^{(\b)\lr}} \lra ]\ol{Y}^{\lr}[_{\cY^{\lr}} $$ 
by $\ti{h}, \breve{h}$, respectively. Then we have 
\begin{align*}
Rf_{(X,\ol{X})/\cY,\rig *}j^{\dd}_X \cE 
& = 
Rh_*\DR^{\d}(]\ol{X}^{(\b)}[_{\cP^{(\b)}}/\cY_K,j^{\dd}_X\cE) \\ 
& = 
R\ti{h}_* \DR^{\d}(]\ol{X}^{(\b) \lr}[_{\cP^{(\b) \lr}}/\cY_K, j^{\dd}_X \cE)
\qquad \text{(\cite[7.3.3, 9.1.1]{chts})} \\ 
& = 
R\ti{h}_* \DR^{\d}(]\ol{X}^{(\b) \lr}[_{\cP^{(\b) \lr}}/\cY^{\lr}_K, 
j^{\dd}_X \cE^{\lr}) \\
& = 
R\epsilon_*R\breve{h}_*\DR^{\d}
(]\ol{X}^{(\b) \lr}[_{\cP^{(\b) \lr}}/\cY^{\lr}_K, 
j^{\dd}_X \cE^{\lr}) \\
& = 
R\epsilon_*Rf_{(X^{\lr},\ol{X}^{\lr})/\cY^{\lr},\rig *}j^{\dd}_X\cE^{\lr}. 
\end{align*}
Then, if we have the isomorphism 
\begin{equation}\label{lr-isom}
j^{\dd}_{Y^{\lr}}R^qf_{\ol{X}^{\lr}/\cY^{\lr},\an *}\cE \os{\cong}{\lra} 
R^qf_{(X^{\lr},\ol{X}^{\lr})/\cY^{\lr}, \rig *}j^{\dd}_X\cE \qquad 
(q \geq 0), 
\end{equation} 
we have 
$$
R^s\epsilon_*R^qf_{(X^{\lr},\ol{X}^{\lr})/\cY^{\lr},\rig *}j^{\dd}_X\cE 
= 
R^s\epsilon_*j^{\dd}_{Y^{\lr}}R^qf_{\ol{X}^{\lr}/\cY^{\lr},\an *}\cE 
= 
R^s\epsilon_*\epsilon^*j^{\dd}_YR^qf_{\ol{X}/\cY,\an *}\cE 
= 0
$$ 
for $s>0$. Hence we have 
\begin{align*}
R^qf_{(X,\ol{X})/\cY,\rig *}j^{\dd}_X \cE & = 
\epsilon_*R^qf_{(X^{\lr},\ol{X}^{\lr})/\cY^{\lr},\rig *}j^{\dd}_X\cE \\ 
& \cong 
\epsilon_* j^{\dd}_{Y^{\lr}}R^qf_{\ol{X}^{\lr}/\cY^{\lr},\an *}\cE 
= 
j^{\dd}_YR^qf_{\ol{X}/\cY,\an *}\cE. 
\end{align*}
So the theorem is reduced to the isomorphism 
\eqref{lr-isom}. To prove the isomorphism 
\eqref{lr-isom}, we may replace $\b$ by $m \in \N$. So, to prove the 
theorem, we may replace $(Y,M_Y) \hra (\cY,M_{\cY})$ by 
$(Y^{\lmr},M_{Y^{\lmr}}) \hra (\cY^{\lmr},M_{\cY^{\lmr}})$, that is, 
we may assume that $(Y,M_Y) \hra (\cY,M_{\cY})$ admits a factorization 
$$ (Y,M_Y) \hra (\cY',M_{\cY'}) \lra (\cY,M_{\cY}), $$
where the first map is an exact closed immersion and 
the second map is formally log etale. \par 
Next, let us put 
$({\cP'}^{(\b)},M_{{\cP'}^{(\b)}}) := 
(\cP^{(\b)},M_{\cP^{(\b)}}) \times_{(\cY,M_{\cY})} (\cY',M_{\cY'})$. 
Then the second projection $({\cP'}^{(\b)},M_{{\cP'}^{(\b)}}) \lra 
(\cY',M_{\cY'})$ and the map $(\cY',M_{\cY'}) \lra (\cY,M_{\cY})$ 
induce the diagram 
\begin{equation*}
\begin{CD}
]\ol{X}^{(\b)}[^{\log}_{{\cP'}^{(\b)}} @>{{h'}^{\log}}>>
]\ol{Y}[^{\log}_{\cY'} @>{\gamma^{\log}}>>
]\ol{Y}[^{\log}_{\cY} \\ 
@VVV @V{\varphi'_Y}VV @V{\varphi_Y}VV \\ 
]\ol{X}^{(\b)}[_{{\cP'}^{(\b)}} @>{h'}>>
]\ol{Y}[_{\cY'} @>{\gamma}>> 
]\ol{Y}[_{\cY} 
\end{CD}
\end{equation*}
and by definition of log tubular neighborhood, 
$\gamma^{\log}, \varphi'_Y$ are isomorphisms. 
(Hence we have $\gamma=\varphi_Y$.) 
Then we have the following isomorphisms: 
\begin{align*}
j^{\dd}_YR^qf_{\ol{X}/\cY,\an *}\cE 
& = 
j^{\dd}_Y\gamma^{\log}_*R^q{h'}^{\log}_*
\DR(]\ol{X}^{(\b)}[^{\log}_{{\cP'}^{(\b)}}/\cY_K, \cE) \qquad 
(*) \\ 
& = 
j^{\dd}_Y\gamma^{\log}_*R^q{h'}^{\log}_*
\DR(]\ol{X}^{(\b)}[^{\log}_{{\cP'}^{(\b)}}/\cY'_K, \cE) \qquad 
(**) \\ 
& = 
\gamma_*j^{\dd}_Y
R^q{h'}^{\log}_*
\DR(]\ol{X}^{(\b)}[^{\log}_{{\cP'}^{(\b)}}/\cY'_K, \cE) \qquad 
(***) \\ 
& = 
\gamma_*j_Y^{\dd}R^qf_{\ol{X}/\cY', \an *}\cE, 
\end{align*}
where $(*), (**), (***)$ are isomorphisms by the following reasons: 
$(*)$ is an isomorphism since we can calculate the 
log analytic cohomology by using the formally log smooth morphism 
$({\cP'}^{(\b)},M_{{\cP'}^{(\b)}}) \lra (\cY,M_{\cY})$. 
$(**)$ is an isomorphism because $(\cY',M_{\cY'}) \lra (\cY,M_{\cY})$ 
is formally log etale. $(***)$ is an isomorphism 
because we have 
$$ j^{\dd}_Y \circ \gamma^{\log}_*  = 
\varphi_{Y,*} \circ j^{\d,\log}_Y \circ \gamma^{\log}_* 
= \gamma_* \circ j^{\d}_Y \circ \varphi'_{Y,*} = 
\gamma_* \circ j^{\dd}_Y $$ 
for coherent modules by Lemma \ref{dagdagdaglem}. 
On the other hand, we have 
\begin{align*}
Rf_{(X,\ol{X})/\cY, \rig *}j^{\dd}_X\cE & = 
R\gamma_*Rh'_*\DR^{\d}(]\ol{X}^{(\b)}[_{{\cP'}^{(\b)}}/\cY_K, j^{\dd}_X\cE) \\ 
& = 
R\gamma_*Rh'_*\DR^{\d}(]\ol{X}^{(\b)}[_{{\cP'}^{(\b)}}/\cY'_K, j^{\dd}_X\cE) \\
& = 
R\gamma_* Rf_{(X,\ol{X})/\cY',\rig *}j^{\dd}_X\cE 
\end{align*}
because we can calculate the relative rigid cohomology 
by using ${\cP'}^{(\b)} \lra \cY$ and $\cY' \lra \cY$ is formally etale 
on a neighborhood of $Y$. So, if we have the isomorphism 
\begin{equation}\label{dash-isom}
j^{\dd}_YR^qf_{\ol{X}/\cY',\an *}\cE \os{\cong}{\lra} 
R^qf_{(X,\ol{X})/\cY', \rig *}j^{\dd}_X\cE, 
\end{equation} 
we have $R^s\gamma_*R^qf_{(X,\ol{X})/\cY', \rig *}j^{\dd}_X\cE = 0$ 
for $s>0$ (which follows from \cite[8.3.5]{chts} and the fact that 
$\gamma^*: \Coh(j^{\d}_Y\cO_{]\ol{Y}[_{\cY}}) \lra 
\Coh(j^{\d}_Y\cO_{]\ol{Y}[_{\cY'}})$ is an equivalence of 
categories) and so we have 
\begin{align*}
R^qf_{(X,\ol{X})/\cY, \rig *}j^{\dd}_X\cE & = 
\gamma_*R^qf_{(X,\ol{X})/\cY', \rig *}j^{\dd}_X\cE \\ & \cong 
\gamma_* j^{\dd}_YR^qf_{\ol{X}/\cY',\an *}\cE = 
j_Y^{\dd}R^qf_{\ol{X}/\cY, \an *}\cE. 
\end{align*}
So the theorem is reduced to the isomorphism \eqref{dash-isom}. 
Hence, to prove the theorem, we may assume that the map 
$(\ol{Y},M_{\ol{Y}}) \hra (\cY,M_{\cY})$ is an exact closed immersion. \par 
Next, let $\{(\cY_m,M_{\cY_m})\}_m$ be the system of universal enlargements of 
$(\ol{Y},M_{\ol{Y}}) \hra (\cY,M_{\cY})$ and let 
$$ (\ol{X}_m,M_{\ol{X}_m}) \lra (\ol{Y}_m, M_{\ol{Y}_m}) \hra 
(\cY_m,M_{\cY_m}) $$ 
be the pull-back of the diagram \eqref{maindiag} by $(\cY_m,M_{\cY_m}) 
\lra (\cY,M_{\cY})$ and put $X_m := X \times_{\ol{X}} \ol{X}_m$. 
Then we have the admissible covering 
$]\ol{Y}[_{\cY} = \bigcup_m \cY_{m,K}$ and it is easy to see that 
the restriction of 
$j^{\dd}_YR^qf_{\ol{X}/\cY, \an *}\cE$ 
(resp. $R^qf_{(X,\ol{X})/\cY, \rig *}j^{\dd}_X\cE$) to $\cY_{m,K}$ 
is naturally isomorphic to 
$j^{\dd}_YR^qf_{\ol{X}_m/\cY_m, \an *}\cE$ 
(resp. $R^qf_{(X_m,\ol{X}_m)/\cY_m, \rig *}j^{\dd}_X\cE$). 
So, we may replace $(\ol{Y},M_{\ol{Y}}) \hra (\cY,M_{\cY})$ by 
$(\ol{Y}_m, M_{\ol{Y}_m}) \hra (\cY_m, M_{\cY_m})$ to prove the theorem: 
So we can reduce to the case where $(\ol{Y},M_{\ol{Y}}) \hra (\cY,M_{\cY})$ 
is a homeomorphic exact closed immersion. \par 
Now we prove that the map \eqref{map} is an isomorphism under the 
assumption that $(\ol{Y},M_{\ol{Y}}) \hra (\cY,M_{\cY})$ is a 
homeomorphic exact closed immersion. To show this, it suffices to 
prove that the map 
\begin{equation}\label{map2}
Rj^{\dd}_YRh^{\log}_*\DR(]\ol{X}^{(\b)}[^{\log}_{\cP^{(\b)}}/\cY_K, \cE) 
\lra 
Rh_*Rj^{\dd}_X \DR(]\ol{X}^{(\b)}[^{\log}_{\cP^{(\b)}}/\cY_K, \cE) 
\end{equation}
in $(\spadesuit)$ is a quasi-isomorphism. We may replace $\b$ by 
$m \in \N$ and by taking a Zariski covering of $\cP^{(m)}$ (and 
the induced covering of $\ol{X}^{(m)}$), 
we may reduce to the case where $\ol{X}^{(m)}$ is affine. Then 
there exists a formally log smooth morphism 
$(\cQ,M_{\cQ}) \lra (\cY,M_{\cY})$ with 
$(\cQ,M_{\cQ}) \times_{(\cY,M_{\cY})} (\ol{Y},M_{\ol{Y}}) = 
(\ol{X}^{(m)},M_{\ol{X}^{(m)}})$. 
Then, the both hand sides and the map in \eqref{map2} (with $\b=m$) 
are unchanged in derived category 
if we replace $(\cP^{(m)},M_{\cP^{(m)}})$ by 
$(\cQ,M_{\cQ})$. 
So we can reduce to the case that the morphism 
$(\cP^{(m)}, M_{\cP^{(m)}}) \lra (\cY,M_{\cY})$ satisfies the 
equality 
$(\cP^{(m)}, M_{\cP^{(m)}}) \times_{(\cY,M_{\cY})} 
(\ol{Y},M_{\ol{Y}}) = (\ol{X}^{(m)},M_{\ol{X}^{(m)}})$. 
In this case, we can easily see that the map \eqref{map2} (with $\b = m$) 
is an isomorphism, because we can see the equality of 
functors $j^{\dd}_Y \circ h^{\log}_* \os{=}{\lra} 
h_* \circ j^{\dd}_X$ in this case, by the quasi-compactness and 
quasi-separatedness of the morphism 
$]\ol{X}[_{\cP^{(m)}} = \cP^{(m)}_K \lra \cY_K = ]\ol{Y}[_{\cY}$. 
So we have finished the proof of the theorem. 
\end{pf}

Now let us take a triple of the form $(S,S,\cS)$ and assume that 
we are given a diagram 
$$ (\ol{X},M_{\ol{X}}) \os{f}{\lra} (\ol{Y},M_{\ol{Y}}) 
\os{g}{\lra} S $$ 
and open immersions $j_X:X \hra \ol{X}, j_Y:Y \hra \ol{Y}$ with 
$X \subseteq (\ol{X},M_{\ol{X}})_{\triv}, 
Y \subseteq (\ol{Y},M_{\ol{Y}})_{\triv}$ and $f^{-1}(Y) =X$, where 
$f$ is a proper log smooth integral 
morphism in $\LB$ having log smooth parameter 
and $g$ is a morphism in $\LB$. Let 
$\cE$ be a locally free isocrystal on 
$(\ol{X}/\cS)^{\log}_{\conv}$. Then we have the following theorem, 
which is the second main result in this section. 

\begin{thm}\label{main2}
Let the notation be as above. Then, for any $q \in \N$, 
there exists uniquely an overconvergent isocrystal $\cF$ on 
$(Y,\ol{Y})/\cS_K$ satisfying the following condition$:$ 
For any strict morphism $(\ol{Z},M_{\ol{Z}}) \lra (\ol{Y}, M_{\ol{Y}})$, 
open immersion $Z \hra \ol{Z}$ with $Z \subseteq Y \times_{\ol{Y}} \ol{Z}$ 
and a closed immersion $(\ol{Z}, M_{\ol{Z}}) \hra (\cZ,M_{\cZ})$ 
into a $p$-adic fine log formal $\cB$-scheme $(\cZ,M_{\cZ})$ 
formally log smooth over $\cS$ 
satisfying $Z \subseteq (\cZ,M_{\cZ})_{\triv}$, 
the restriction of $\cF$ to the category 
$I^{\d}((Z,\ol{Z})/\cS,\cZ)$ is given functorially by 
$(R^qf_{(X \times_Y Z, \ol{X} \times_{\ol{Y}} \ol{Z})/\cZ, \rig *} 
j_X^{\dd}\cE, 
\epsilon)$, where $\epsilon$ is given by 
\begin{align*}
p_2^*R^qf_{(X \times_Y Z, \ol{X} \times_{\ol{Y}} \ol{Z})/\cZ, \rig *} 
j_X^{\dd}\cE & \os{\sim}{\ra} 
R^qf_{(X \times_Y Z, \ol{X} 
\times_{\ol{Y}} \ol{Z})/\cZ \times_{\cS} \cZ, \rig *} j_X^{\dd}\cE \\ 
& \os{\sim}{\leftarrow}
p_1^*R^qf_{(X \times_Y Z, \ol{X} 
\times_{\ol{Y}} \ol{Z})/\cZ, \rig *} j_X^{\dd}\cE. 
\end{align*}
$($Here $p_i \,(i=1,2)$ 
denotes the $i$-th projection $]\ol{Z}[_{\cZ \times_{\cS} \cZ} \lra 
]\ol{Z}[_{\cZ}.)$ 
\end{thm} 

\begin{pf} 
By Corollary \ref{iso2}, there exists an isocrystal $\ti{\cF}$ on 
$(\ol{Y}/\cS)^{\log}_{\conv}$ such that, 
for the closed immersion 
$(\ol{Z},M_{\ol{Z}}) \hra (\cZ,M_{\cZ})$ as above, $\ti{\cF}$ induces via 
the functor 
$$ I_{\conv}((Y/\cS)^{\log}) \lra 
I_{\conv}((Z/\cS)^{\log}) \simeq \Strat''((Z \hra \cZ/\cS)^{\log}) $$
the object 
$(R^qf_{\ol{X} \times_{\ol{Y}} \ol{Z}/\cZ,\an *}\cE, \ti{\epsilon})$, where 
$\ti{\epsilon}$ is the canonical isomorphism 
$$ p_2^*R^qf_{\ol{X} \times_{\ol{Y}} \ol{Z}/\cZ,\an *}\cE  \os{\sim}{\ra} 
R^qf_{\ol{X} \times_{\ol{Y}} \ol{Z}/\cZ \times_{\cS} \cZ,\an *}
\cE \os{\sim}{\leftarrow}
p_1^*R^qf_{\ol{X} \times_{\ol{Y}} \ol{Z}/\cZ,\an *}\cE. $$ 
$($Here $p_i$ denotes the $i$-th projection 
$]\ol{Z}[^{\log}_{\cZ \times_{\cS} \cZ} \lra 
]\ol{Z}[^{\log}_{\cZ}.)$ Then, by Theorem \ref{main1}, one can see that the overconvergent isocrystal 
$\cF := j^{\dd}_Y\ti{\cF}$ is the one which satisfies the condition of 
the theorem. \par 
We prove the uniqueness of $\cF$. Take an embedding system 
\begin{equation}\label{mainembsys}
(\ol{Y},M_{\ol{Y}}) \lla (\ol{Y}^{(\b)},M_{\ol{Y}^{(\b)}}) 
\hra (\cY^{(\b)},M_{\cY^{(\b)}}) 
\end{equation} 
satisfying $Y^{(\b)} \subseteq (\cY^{(\b)},M_{\cY^{(\b)}})_{\triv}$, 
where $Y^{(\b)}:=Y \times_{\ol{Y}} \ol{Y}^{(\b)}$. 
Let us denote the category of descent data with respect to 
$I^{\d}((Y^{(n)},\ol{Y}^{(n)})/\cS_K, \cY^{(n)}) \,(n=0,1,2)$ by 
$I^{\d}((Y^{(\b)},\ol{Y}^{(\b)})/\cS_K, \allowbreak \cY^{(\b)})$. Then, by the 
condition required for $\cF$ in the case where 
$Z, (\ol{Z},M_{\ol{Z}}), (\cZ,M_{\cZ})$ are equal to 
$Y^{(n)}, (\ol{Y}^{(n)},M_{\ol{Y}^{(n)}}), (\cY^{(n)},M_{\cY^{(n)}})$ 
respectively (for $n=0,1,2$), the restriction of $\cF$ by the functor 
\begin{equation}\label{equivalence}
I^{\d}((Y,\ol{Y})/\cS_K) \lra 
I^{\d}((Y^{(\b)},\ol{Y}^{(\b)})/\cS_K, \cY^{(\b)})
\end{equation} 
is uniquely determined. Since the functor \eqref{equivalence} is 
an equivalence of categories by Proposition \ref{etaledescent}, 
we have the uniqueness of $\cF$. So the proof of the theorem is finished. 
\end{pf} 

\begin{cor}\label{main3}
Let $(X,\ol{X}) \os{f}{\lra} (Y,\ol{Y}) \os{g}{\lra} (S,S)$ be the morphism 
of pairs induced by the situation in Theorem \ref{main2} and let 
$\cE$ be as in Theorem \ref{main2}. Then Conjecture \ref{bconj2} 
$($the version without Frobenius structure$)$ is 
true for the overconvergent isocrystal $j_X^{\dd}\cE$. 
\end{cor} 

\begin{pf} 
Let us take a diagram 
\begin{equation}\label{maincov}
(\ol{Y},M_{\ol{Y}}) \os{g^{(0)}}{\lla} (\ol{Y}^{(0)},M_{\ol{Y}^{(0)}}) 
\os{i^{(0)}}{\hra} (\cY^{(0)},M_{\cY^{(0)}}), 
\end{equation}
where $g^{(0)}$ is a strict etale surjective morphism, 
$i^{(0)}$ is a closed immersion in $\pLF$ with 
$Y^{(0)} \subseteq (\cY^{(0)},M_{\cY^{(0)}})_{\triv}$ (where 
$Y^{(0)}:= Y \times_{\ol{Y}} \ol{Y}^{(0)}$), 
$(\cY^{(0)},M_{\cY^{(0)}})$ is formally log smooth over $\cS$ and 
$\cY^{(0)}$ is formally smooth over $\cS$. Using this diagram, 
we define the category $\cC$ as the category of triples 
$(Z,\ol{Z},\cZ)$ over $(Y^{(0)},\ol{Y}^{(0)},\cY^{(0)})$ such that 
$\cZ$ is formally smooth over $\cY^{(0)}$. \par 
Let us take an object $(Z,\ol{Z},\cZ)$ in $\cC$ and let 
$M_{\ol{Z}}, M_{\cZ}$ be the log structure on $\ol{Z},\cZ$ defined 
as the pull-back of $M_{\ol{Y}^{(0)}}, M_{\cY^{(0)}}$, 
respectively. Then the open immersion $Z \hra \ol{Z}$
 and the closed immersion 
$(\ol{Z},M_{\ol{Z}}) \hra (\cZ,M_{\cZ})$ satisfies the condition 
required in the statement of Theorem \ref{main2}. So, by 
Theorem \ref{main2}, the restriction of the overconvergent isocrystal 
$\cF$ by $I^{\d}((Y,\ol{Y})/\cS_K) \lra I^{\d}((Z,\ol{Z})/\cS_K,\cZ)$ 
is given by 
$(R^qf_{(X \times_Y Z, \ol{X} \times_{\ol{Y}} \ol{Z})/\cZ, \rig *} 
j_X^{\dd}\cE, \epsilon)$ as in the statement of Theorem \ref{main2}. 
Moreover, if we define 
$(\ol{Y}^{(n)},M_{\ol{Y}^{(n)}})$ (resp. $(\cY^{(n)},M_{\ol{Y}^{(n)}})$) 
to be the $(n+1)$-fold fiber product of 
$(\ol{Y}^{(0)},M_{\ol{Y}^{(0)}})$ (resp. $(\cY^{(0)},M_{\ol{Y}^{(0)}})$) 
over $(\ol{Y},M_{\ol{Y}})$ (resp. $\cS$) and if we put 
$Y^{(n)}:= Y \times_{\ol{Y}} \ol{Y}^{(n)}$, we have the embedding 
system \eqref{mainembsys} and the triple $(Y^{(n)},\ol{Y}^{(n)},\cY^{(n)})$ 
is contained in the category $\cC$. From these facts, we can deduce the 
uniqueness of the overconvergent isocrystal $\cF$ in the same way as 
the proof of Theorem \ref{main2}. So we are done. 
\end{pf} 

Finally, we give a result on Frobenius structure on the overconvergent 
isocrystal which is constructed in Theorem \ref{main2}. \par 
In the following, let us assume that $k$ is perfect and that 
$(\cB, M_{\cB})$ is the formal scheme $\Spf W$ endowed with trivial 
log structure. (So $(B,M_B):=(\cB,M_{\cB}) \otimes_{\Z_p} \F_p$ is 
the scheme $\Spec k$ endowed with trivial log sturucture.) 
Let $V$ be a finite totally ramified extension of $W$, let 
$\pi$ be a uniformizer of $V$ and put 
$\cS := \Spf V$, $S := \Spec V/\pi V = B$. \par 
Let us fix an integer $q$ which is a power of $p$. For a fine log scheme 
$(X,M_X)$ over $\F_p$, let $F_X$ be the $q$-th power Frobenius 
(the endomorphism $(X,M_X) \lra (X,M_X)$ induced by the $q$-th power 
endomorphism of the structure sheaf). 
Assume that we have an endomorphism $\sigma: \cS \lra \cS$ which 
lifts $F_S:S \lra S$ and fix it. Then, if we have a pair $(X,\ol{X})$ 
over $(S,S)$, 
we have the canonical functor
$$ 
F^*_X: I_{\conv}((X/\cS)^{\log}) \lra I_{\conv}((X/\cS)^{\log}) $$
induced by the morphisms $F_X$ and $\sigma$. An 
overconvergent $F$-isocrystal on 
$(X,\ol{X})/\cS_K$ is defined to be a pair $(\cG,\alpha)$, where 
$\cG$ is an overconvergent isocrystal on $(X,\ol{X})/\cS_K$ and $\alpha$ is 
an isomorphism $F^*_X\cG \os{\sim}{\lra} \cG$. $\alpha$ is called a 
Frobenius structure on $\cG$. \par 
Now assume that 
we are given a diagram 
$$ (\ol{X},M_{\ol{X}}) \os{f}{\lra} (\ol{Y},M_{\ol{Y}}) 
\os{g}{\lra} S $$ 
and open immersions $j_X:X \hra \ol{X}, j_Y:Y \hra \ol{Y}$ with 
$X \subseteq (\ol{X},M_{\ol{X}})_{\triv}, 
Y \subseteq (\ol{Y},M_{\ol{Y}})_{\triv}$ and $f^{-1}(Y)=X$, 
where $f$ is a proper log smooth integral 
morphism in $\LB$ having log smooth parameter 
and $g$ is a morphism in $\LB$. 
Let $\cE$ be a locally free isocrystal on $(\ol{X}/\cS)^{\log}_{\conv}$ 
such that $j^{\dd}\cE$ (which is an overconvergent isocrystal on 
$(X,\ol{X})/\cS_K$) admits a Frobenius structure $\alpha: F_X^*j_X^{\dd}\cE 
\os{\simeq}{\lra} j_X^{\dd}\cE$. 
By Corollary \ref{main3}, 
we have the $q$-th rigid cohomology overconvergent isocrystal 
$\cF$ of $j_X^{\dd}\cE$, which is an overconvergent isocrystal on 
$(Y,\ol{Y})/\cS_K$. 
Then we have the 
following theorem: 

\begin{thm} 
With the above notation, there exists the canonical Frobenius 
structure on $\cF$ which is induced by the Frobenius 
structure $\alpha$ on $j_X^{\dd}\cE$. 
\end{thm} 

\begin{pf} 
Let us take a diagram 
\begin{equation*}
(\ol{Y},M_{\ol{Y}}) \os{g^{(0)}}{\lla} (\ol{Y}^{(0)},M_{\ol{Y}^{(0)}}) 
\os{i^{(0)}}{\hra} (\cY^{(0)},M_{\cY^{(0)}}) 
\end{equation*}
as in the proof of Corollary \ref{main3} endowed with 
the endomorphism $\sigma_{\cY^{(0)}}:
(\cY^{(0)},M_{\cY^{(0)}}) \lra (\cY^{(0)},M_{\cY^{(0)}})$ 
which is compatible with $F_{\ol{Y}^{(0)}}$ and $\sigma$, and let 
us define the category $\ti{\cC}$ as the category of triples 
$(Z,\ol{Z},\cZ)$ in $\cC$ endowed with an endomorphism 
$\sigma_{\cZ}: \cZ \lra \cZ$ compatible with $F_{\ol{Z}}, \sigma$ and 
$\sigma_{\cY^{(0)}}$. For $(Z,\ol{Z},\cZ,\sigma_{\cZ})$ in 
$\ti{\cC}$, 
let $\cF_{\cZ}$ be the restriction of $\cF$ to 
$I^{\d}((Z,\ol{Z})/\cS_K,\cZ)$. Then it suffices to introduce 
a Frobenius structure on $\cF_{\cZ}$ in functorial way: Indeed, 
if we define 
$(\ol{Y}^{(n)},M_{\ol{Y}^{(n)}})$, $(\cY^{(n)},M_{\ol{Y}^{(n)}})$, 
$Y^{(n)}$ as in the proof of Corollary \ref{main3} and 
if we define $\sigma_{\cY^{(n)}}: (\cY^{(n)},M_{\ol{Y}^{(n)}}) \lra 
(\cY^{(n)},M_{\ol{Y}^{(n)}})$ as the morphism induced by $\sigma_{\cY^{(0)}}$, 
we have the embedding system \eqref{mainembsys} endowed with 
$\sigma_{\cY^{(\b)}}$ and we have 
$(Y^{(n)},\ol{Y}^{(n)},\cY^{(n)},\sigma_{\cY^{(n)}}) \in \ti{\cC}$. 
So we can define the Frobenius structure on the restriction of $\cF$ by 
the equivalence of categories \eqref{equivalence}. \par 
So let us take $(Z,\ol{Z},\cZ,\sigma_{\cZ})$ in $\ti{\cC}$ and 
rewrite $Z \times_Y X, \ol{Z} \times_{\ol{Y}} \ol{X}, Z,\ol{Z}, 
\cZ, \sigma_{\cZ}$ as $X,\ol{X},Y,\ol{Y},\cY,\allowbreak \sigma_{\cY}$, 
respectively. Then we are in the situation where 
there exists a closed immersion 
$(\ol{Y},M_{\ol{Y}}) \hra (\cY,M_{\cY})$ into 
a fine log formal $\cB$-scheme $(\cY,M_{\cY})$ formally log smooth 
over $\cS$ with $Y \subseteq (\cY,M_{\cY})_{\triv}$, 
endowed with an endomorphism $\sigma_{\cY}:(\cY,M_{\cY}) \lra 
(\cY,M_{\cY})$ compatible with $F_{\ol{Y}}$ and $\sigma$. 
In this situation, 
we have the canonical homomorphism 
of $j_Y^{\d}\cO_{]\ol{Y}[_{\cY}}$-modules 
\begin{equation}\label{mapmap}
\sigma_{\cY}^*R^qf_{(X,\ol{X})/\cY,\rig *}j^{\dd}_X\cE 
\os{F^*}{\lra} 
R^qf_{(X,\ol{X})/\cY, \rig *}F_X^*j^{\dd}_X\cE 
\os{\alpha}{\lra}
R^qf_{(X,\ol{X})/\cY \rig *}j^{\dd}_X\cE 
\end{equation} 
(where $F^*$ is induced by $F_X^*$ and $\sigma_{\cY}$) 
and it suffices to prove that this homomorphism is an isomorphism. \par 
We prove the above claim, following partly \cite[3.3.3, 4.1.4]{tsuzuki3}. 
By using the faithfulness of the restriction functor 
$\Coh(j^{\d}_Y\cO_{]\ol{Y}[_{\cY}}) \lra 
\Coh(\cO_{]Y[_{\cY}})$, we can reduce to the case $Y=\ol{Y}$ 
(and $X=\ol{X}, j_X^{\dd}\cE=\cE$) to prove the 
claim. Moreover, by using \cite[9.4.7]{bgr}, we may reduce to the 
case where $Y =: y$ is a closed point. Let $k'$ be a finite Galois extension 
of $k$ containing the function field of $y$, let $S' := \Spec k', 
\cS':= \Spf V'$ (where $V'$ is the unramified extension of $V$ with 
residue field $k'$) and 
let us denote the base change of the diagram 
$X \os{f}{\lra} y \os{\iota}{\hra} \cY$ by $\cS' \lra \cS$ by 
$X' \os{f}{\lra} y' \os{\iota}{\hra} \cY'$. Then, since $\cS' \lra \cS$ 
is finite etale, we may replace $X,y,\cY,S,\cS$ by $X',y',\cY',S',\cS'$ 
to prove 
the claim. Since $y'$ is a disjoint union of finite number of the schemes 
isomorphic to $S'$, we may assume that the function field of $y$ is 
equal to $k$ to prove the claim (for $X \lra y \hra \cY$). Then we can form 
the following diagram: 
\begin{equation*}
\begin{CD}
X @>f>> y @>{\iota}>> \cY \\
@\vert @\vert @VVV \\ 
X @>f>> y @>>> \cS, 
\end{CD}
\end{equation*}
where the lower right horiozontal arrow is the closed immersion induced by 
the surjection $V \lra k$. Then, by the analytically flat base change 
theorem, it suffices to prove the claim for the lower horizontal line 
in the above diagram, that is, we may assume $\cY = \cS$. 
In this case, the homomorphism \eqref{mapmap} is nothing but the 
homomorphism of convergent cohomologies 
\begin{equation}\label{convfrob}
H^q((X/\cS)_{\conv},\cE) \os{F^*}{\lra} 
H^q((X/\cS)_{\conv},F_X^*\cE) \os{\alpha}{\cong} H^q((X/\cS)_{\conv},\cE)
\end{equation}
induced by $(\cE,\alpha)$. Moreover, 
by \cite[5.1]{tsuzuki2}, we have the direct image functor 
$\tau_*: I_{\conv}(X/\cS) \allowbreak \lra I_{\conv}(X/\cB)$ which sends 
convergent $F$-isocrystals $(\cE,\alpha)$ on $(X/\cS)_{\conv}$
to convergent $F$-isocrystals $(\tau_*\cE,\tau_*\alpha)$ on 
$(X/\cB)_{\conv}$, and the homomorphism \eqref{convfrob} is rewritten as 
\begin{equation*}
H^q((X/\cB)_{\conv},\tau_*\cE) \os{F^*}{\lra} 
H^q((X/\cB)_{\conv},F_X^*\tau_*\cE) \os{\tau_*\alpha}{\cong} 
H^q((X/\cB)_{\conv},\cE). 
\end{equation*}
So we are reduced to the case $\cS=\cB$. In this case, 
the homomorphism \eqref{convfrob} is 
identical with the following homomorphism of crystalline 
cohomologies
\begin{equation}\label{crysmap}
H^q((X/\cB)_{\crys},\Phi(\cE)) \os{F^*}{\lra} 
H^q((X/\cB)_{\crys},F^*\Phi(\cE)) \os{\alpha}{\cong} 
H^q((X/\cB)_{\crys},\Phi(\cE)). 
\end{equation}
So we are reduced to showing that the first map in the diagram \eqref{crysmap}
is an isomorphism. To prove this assertion, we may reduce to the case 
$q=p$. 
Denote the fiber product $B \times_{F_B,B} X$ by $X'$ and factorize 
the Frobenius map $F_X:X \lra X$ 
as $X \os{F_{X/B}}{\lra} X' \os{\pr}{\lra} X$, 
where $F_{X/B}$ is relative Frobenius and $\pr$ is the projection. 
Then, if we denote $\pr^*\Phi(\cE)$ by $\cE'$, 
the first map in \eqref{crysmap} factors as 
\begin{equation}\label{crysmap2}
H^q((X/\cB)_{\crys},\Phi(\cE)) \os{\pr^*}{\lra} 
H^q((X'/\cB)_{\crys},\cE') \os{F^*_{X/B}}{\lra} 
H^q((X/\cB)_{\crys},F^*_{X/B}\cE'). 
\end{equation}
The first map in \eqref{crysmap2} is an 
isomorphism by base change theorem (Theorem \ref{crysbc}). 
So it suffices to prove that the second map is an isomorphism. 
In this case, note that $\cE'$ has the form $\cE' = \Q \otimes \cF$, 
where $\cF$ is a $p$-torsion free crystal on $(X/\cB)_{\crys}$ in 
the sense in \cite[5.1]{ogus4}. (It follows from the fact that 
$L''$ in Lemma \ref{crys-hpdilem} is $p$-torsion free. 
See also \cite[0.7.5]{ogus2}.) So the triple 
$(\cF,F^*_{X/B}\cF, \id:F^*_{X/B}\cF \ra F^*_{X/B}\cF)$ forms an 
$F$-span on $X/\cB$ in the sense in \cite[5.2.1]{ogus4}. 
Moreover, by \cite[5.2.9]{ogus4}, this triple is automatically an 
admissible $F$-span, because we have 
$R_{X/\cB}=0$ in the notation there since 
the log structures on $X,\cB$ are trivial. Then, by \cite[7.3.1]{ogus4}, 
we can conclude that the second map in \eqref{crysmap2} is an isomorphism. 
So we have finished the proof of the theorem. 
\end{pf} 


\end{document}